\documentclass{amsart}
\usepackage{amsmath}
\usepackage{amsxtra}
\usepackage{amstext}
\usepackage{amssymb}
\usepackage{latexsym}
\usepackage{dsfont} 
\usepackage{verbatim}
\usepackage{tabls}
\usepackage{rotating}
\usepackage{color}

\newtheorem{thm}{Theorem}[section]

\newtheorem{cor}[thm]{Corollary}
\newtheorem{lem}[thm]{Lemma}
\newtheorem{prop}[thm]{Proposition}
\newtheorem{ques}[thm]{Problem}
\newcommand{\abs}[1]{\left\vert#1\right\vert}
\newcommand{\norm}[1]{\left\Vert#1\right\Vert}
\theoremstyle{definition}
\newtheorem{defn}[thm]{Definition}

\newtheorem{rem}[thm]{Remark}

\numberwithin{equation}{section}
\def\XXint#1#2#3{{\setbox0=\hbox{$#1{#2#3}{\int}$}
     \vcenter{\hbox{$#2#3$}}\kern-.5\wd0}}

\begin{document}
\title[Nonlinear equations with natural growth terms]{Local and global behaviour of nonlinear equations with natural growth terms}

\author{Benjamin J. Jaye}
\address{Department of Mathematics, Kent State University, Kent, Ohio 44240, USA}
\email{bjaye@kent.edu}

\author{Igor E. Verbitsky}
\address{Department of Mathematics, University of Missouri, Columbia, Missouri 65211, USA}
\email{verbitskyi@missouri.edu}

\thanks{Supported in part
by NSF grant  DMS-0901550.}

\dedicatory{In memory of Professor Nigel Kalton}

\begin{abstract}In this paper we  study the Dirichlet problem 
\begin{equation*}\begin{cases}\label{introdir}-\Delta_p u = \sigma |u|^{p-2}u + \omega \text{ in } \Omega,\\
\; u = 0\text{ on }\partial\Omega,\end{cases}\end{equation*} where $\sigma$ and $\omega$ are nonnegative Borel measures, and 
$\Delta_p u = \nabla \cdot (\nabla u \, |\nabla u|^{p-2})$ is the $p$-Laplacian.   Here $\Omega\subseteq \mathbf{R}^n$ is either a bounded domain, or the entire space.  Our main estimates concern optimal pointwise bounds of solutions  in terms of two local Wolff's potentials, under minimal regularity assumed on $\sigma$ and $\omega$.  
In addition, analogous results for equations modeled by the $k$-Hessian in place of the $p$-Laplacian will be discussed.
\end{abstract}

\keywords{Quasilinear equations, natural growth terms, Wolff potentials, discrete Carleson measures}

\subjclass[2000]{Primary 35J60,  42B37. Secondary 31C45, 35J92, 42B25}


\maketitle

\section{Introduction}

\subsection{}In this paper we develop an approach to studying the local and global pointwise behaviour of solutions to equations with natural growth terms, under minimal regularity assumptions.  Let $\Omega\subseteq\mathbf{R}^n$ be an open set, with $n\geq 2$, and let $1<p<n$.  The model problem under consideration is the following Dirichlet problem:
\begin{equation}\label{quas2w}
\begin{cases} \,-\Delta_p u = \sigma u^{p-1}+\omega \; \text{ in } \Omega, \\
\, u>0 \, \text{ in } \Omega, \\
\,  u=0 \, \text{ on } \partial \Omega. \end{cases}
\end{equation}
Throughout this paper we will assume that the \textit{potential} $\sigma$ is a locally finite nonnegative measure, 
and the \textit{inhomogeneous term} $\omega$ is a finite nonnegative measure.  Here $\Delta_p u = \nabla \cdot (\abs{\nabla u}^{p-2}\nabla u)$ is the $p$-Laplacian operator.  Since our aim is to study (\ref{quas2w}) in a low regularity setting, in all our results the $p$-Laplacian operator can be replaced by a general second order quasilinear operator with standard structural assumptions, for instance the $\mathcal{A}$-Laplacian operator (see e.g. \cite{HKM}).

The equation (\ref{quas2w}) is a very natural perturbation of the $p$-Laplacian operator, and the local behaviour of solutions to equations of the type (\ref{quas2w}) is a heavily studied topic, beginning with the classic works of Serrin \cite{Ser2}; Trudinger \cite{T67}; and Ladyzhenskaya and Ural'tseva \cite{LU68}, where suitable $L^q$ assumptions are imposed on $\sigma$ and $\omega$.
The purpose of this article is to study the \textit{pointwise} behaviour of solutions to (\ref{quas2w}), including the cases where the potential and data are too rough to fall within the studies previously cited.  In particular classical tools such as Harnack's inequality are no longer valid in general for positive solutions of (\ref{quas2w}).

In recent papers \cite{AHBV} and  \cite{Gre02} it is pointed out that the existence problem for (\ref{quas2w}) is non-trivial for general measure right hand side $\omega$, even under the assumption that $\sigma \in L^q(\Omega)$ for $q>n/p$.  In these papers, the existence problem for (\ref{quas2w}) for measure data $\omega$ is solved under the assumption that $\sigma \in L^q(\Omega)$ for $q>n/p$ with small norm.  If one  avoids   the phenomenon of interaction between $\sigma$ and $\omega$, then 
a simple analysis (see Remark 6.1 of \cite{AHBV}) shows that this $L^q$ class of potentials $\sigma$ is optimal on the Lebesgue scale
in order to solve the equation (\ref{quas2w}) for all finite measures $\omega$.

Here we investigate solutions of (\ref{quas2w}) taking into account the interaction between the two terms $\sigma$ and $\omega$.  The problem turns out to be less robust than the \textit{super-critical}  case studied earlier in \cite{PV0, PV, PV2}, where the $\sigma u^{p-1}$ term in (\ref{quas2w}) is replaced by $\sigma u^q$ with $q>p-1$.  The equations with natural growth terms $q=p-1$ have all the hallmarks of the end-point case where more subtle methods of analysis are in order.  For example, in what follows we will make extensive use of John-Nirenberg type BMO estimates in weighted spaces where the underlying measure is non-doubling.

As a result of our study, existence results are extended to more general classes of measures $\sigma$ which could be singular with respect to Lebesgue measure.  More salient is that our approach reveals  pointwise behaviour of solutions. The latter reduces questions of finer regularity of solutions $u$ of (\ref{quas2w}) to merely checking norm mapping properties of certain nonlinear integral operators.  By now this is a well developed approach to deducing fine properties of nonlinear equations, see e.g. \cite{DM10, DM11, Min07, Min11, PV,  PV2} and references therein. 

The results of this paper are bound to be complicated due to the two-weight nature of the problem (the interaction between $\sigma$ and $\omega$).  To compensate for this, we provide several examples of  classes of both $\sigma$ and $\omega$ where our main theorems can be applied and where the pointwise expressions we obtain for solutions of (\ref{quas2w}) simplify.

\subsection{}In order to motivate our main results, we will first discuss the problem in the entire space $\mathbf{R}^n$. In this case the Dirichlet problem (\ref{quas2w}) reads:
\begin{equation}\label{quas2wentire}
\begin{cases} \,-\Delta_p u = \sigma u^{p-1}+\omega \; \text{ in } \mathbf{R}^n, \\
\, \inf_{x\in \mathbf{R}^n} u(x) = 0. \end{cases}
\end{equation}
In equation (\ref{quas2wentire}), solutions are considered in the sense of $p$-superharmonic functions, see Section \ref{renormsoln}.  When dealing with Dirichlet boundary conditions in bounded domains we will work with the stronger notion of renormalized solutions.

In our previous work \cite{JVFund}, we studied an important special case of (\ref{quas2w}), namely, the \textit{fundamental solution}:
\begin{equation}\label{fundsoln}
-\Delta_p u = \sigma u^{p-1} + \delta_{x_0} \text{ in } \mathbf{R}^n,\; \text{ and } \inf_{x\in \mathbf{R}^n}u(x) = 0.
\end{equation}
Here $\delta_{x_0}$ is the Dirac delta measure with pole at $x_0$.  By producing sharp global pointwise bounds for solutions (see (\ref{fundsolnrecap}) below), we showed that the problem is controlled by two local potentials; the local nonlinear Wolff's potential, 
and a local fractional linear potential (defined in (\ref{Riesz}) and (\ref{Wolff}) respectively).
These two potentials will play a prominent role in what follows.  

Furthermore, it was shown that a necessary condition for the existence of a positive supersolution (in any reasonable sense), i.e., the integral inequality:
$$-\Delta_p u \geq \sigma u^{p-1} \text{ in }\mathbf{R}^n,
$$
is that the potential $\sigma$ should satisfy the condition:
\begin{equation}\label{capstart}
\sigma(E) \leq C(\sigma) \, \text{cap}_p(E)\text{ for all compact sets }E\subset\mathbf{R}^n,
\end{equation}
with $C(\sigma)=1$.  Here $\text{cap}_p$ is the standard $p$-capacity relative to $\mathbf{R}^n$ (see (\ref{entirecap}) below).  It is known that in many nonlinear elliptic problems with measure data working with this capacity is very natural, see e.g. \cite{DMMOP99}.  With this in mind, a primary result of \cite{JVFund} can be summarized as follows:

\textit{There exists a positive constant $C=C(n,p)$ so that if $\sigma$ satisfies (\ref{capstart}) with $C(\sigma)< C$, then there exists a solution of (\ref{fundsoln}) so that}:
\begin{equation}\begin{split}\label{fundsolnrecap}
u(x,x_0) \approx  c\abs{x-x_0}^{\frac{p-n}{p-1}}& \exp\Big( c\int_{0}^{\abs{x-x_0}}\Bigl(\frac{\sigma(B(x,r))}{r^{n-p}}\Bigl)^{1/(p-1)}\frac{dr}{r}\Bigl)\\
&\cdot \exp\Bigl(c\int_{0}^{\abs{x-x_0}} \frac{\sigma(B(x_0,r))}{r^{n-p}}\frac{dr}{r}\Bigl).
\end{split}\end{equation}

Here and elsewhere in the paper, the symbol $\approx$ denotes a bilateral pointwise bound, so that the constant $c=c(n,p)>0$ in (\ref{fundsolnrecap}) may differ on each side of the bound.  In other words, $u(x)\approx cf(x)e^{cg(x)}$, if there exist constants $c_1, c_2>0$ 
which depend only on $n$ and $p$  so that:
$$c_1 f(x)e^{c_1g(x)}\leq u(x) \leq c_2f(x)e^{c_2g(x)},
$$
where $u, f, g$ are nonnegative functions. 
Our principal aim is to extend such a result to when $\delta_{x_0}$ is replaced by a general measure $\omega$.  We will see that this generalisation is by no means straightforward.

The bound (\ref{fundsolnrecap}) leads to a natural candidate for a global pointwise bound for solutions of (\ref{quas2wentire}) for general $\omega$.  Indeed, one might expect to be able to find a solution of (\ref{quas2wentire}) with the bilateral estimate:\begin{equation}\begin{split}\label{falseconj1}
u(x) \approx c \int_0^{\infty}&\Bigl[\frac{1}{r^{n-p}}\exp\Big( c\int_{0}^{r}\Bigl(\frac{\sigma(B(x,s))}{s^{n-p}}\Bigl)^{1/(p-1)}\frac{ds}{s}\Bigl)\\
& \cdot \int_{B(x,r)}\exp\Bigl(c\int_{0}^{r} \frac{\sigma(B(z,s))}{s^{n-p}}\frac{ds}{s}\Bigl)d\omega(z)\Bigl]^{1/(p-1)}\frac{dr}{r}.
\end{split}\end{equation}

There is substantial evidence to support (\ref{falseconj1}).  The reader can check that it coincides with (\ref{fundsolnrecap}) when $\omega = \delta_{x_0}$.  Second, it would recover the pointwise bounds already found very recently in the linear case $p=2$ in \cite{FV1, FNV}.  However, in general (\ref{falseconj1}) turns out to be false.

In order to see that (\ref{falseconj1}) fails in general, we introduce another class of examples, of interest in their own right.  Let us consider solutions $u$ of:
\begin{equation}\begin{cases}\label{1rhsentire}
-\Delta_p u = \sigma u^{p-1} \text{ in } \mathbf{R}^n,\\
\; \inf_{x\in \mathbf{R}^n} u=1.
\end{cases}\end{equation}
This equation has been heavily studied (in both the entire space and in domains) in the case $p=2$, where it is related to the so-called \textit{gauge}, or the Feyman-Kac functional, see e.g. \cite{CZ95, FNV}.  We will see that for the purposes of the pointwise bounds in this paper, this equation is essentially equivalent to the problem:
\begin{equation}\begin{cases}\label{1bventire}
-\Delta_p u = \sigma u^{p-1} + \sigma \text{ in } \mathbf{R}^n,\\
\; \inf_{x\in \mathbf{R}^n} u=0,
\end{cases}\end{equation}
which is of the form (\ref{quas2w}).
Both equations (\ref{1rhsentire}) and (\ref{1bventire}) will be studied in Section \ref{1bvsection}, where it will be shown that \textit{there exists a constant $C=C(n,p)$ so that if $\sigma$ satisfies (\ref{capstart}) with $C(\sigma) < C$, then there exists a solution $u$ of (\ref{1bventire}) so that}:
\begin{equation}\label{1bventirebd}
u(x)\approx \Bigl[\exp\Bigl(c\int_0^{\infty}\Bigl(\frac{\sigma(B(x,r))}{r^{n-p}}\Bigl)^{1/(p-1)} \frac{dr}{r}\Bigl) - 1\Bigl].
\end{equation}

The main observation regarding the bound (\ref{1bventirebd}) is that the linear potential $\mathbf{I}^r_{\alpha}(d\sigma)(x)$ does not appear at all.  Consequently, if $1<p< 2$, one can find examples of $\sigma$ so that quantity appearing in (\ref{falseconj1}) is identically infinite (when $\omega = \sigma$), but also so that the quantity in (\ref{1bventirebd}) is uniformly bounded\footnote{One way to do this is to pick a set $E\subset B(0,1)$ so that the Riesz capacity $\text{cap}_{1,p}(E)>0$, but $\text{cap}_{p/2,2}(E)=0$. This is possible if and only if $1<p<2$, see Theorem 5.5.1 of \cite{AH}.  Then choose $\sigma$ to be the capacitary measure associated with $E$. By the dual definition of capacity (Theorem 2.5.2 of \cite{AH}) the result follows.}.  It follows that the bound (\ref{falseconj1}) is not sharp in general.

From the bound (\ref{1bventirebd}), along with a simple summation by parts argument, one is lead to the another potential bound for solutions of (\ref{quas2wentire}):
\begin{equation}\begin{split}\label{falseconj2}
u(x) \approx  \int_0^{\infty}&\Bigl[\frac{1}{r^{n-p}}\exp\Big( c\int_{0}^{r}\Bigl(\frac{\sigma(B(x,s))}{s^{n-p}}\Bigl)^{1/(p-1)}\frac{ds}{s}\Bigl)\\
& \cdot \int_{B(x,r)}\exp\Bigl(c\int_{0}^{r} \Bigl(\frac{\sigma(B(z,s))}{s^{n-p}}\Bigl)^{1/(p-1)}\frac{ds}{s}\Bigl)d\omega(z)\Bigl]^{1/(p-1)}\frac{dr}{r}.
\end{split}\end{equation}

The bound (\ref{falseconj2}) coincides (up to multiplicative constants) with (\ref{1bventirebd}) if $\omega = \sigma$.  However, when $\omega = \delta_{x_0}$, (\ref{falseconj2}) clearly does not match (\ref{fundsolnrecap}).

\subsection{An example theorem} It turns out that a combination of the two bounds (\ref{falseconj1}) and (\ref{falseconj2}) yields optimal pointwise estimates for solutions of (\ref{quas2wentire}).  As the discussion above shows, such a result should depend on whether $1<p\leq 2$ or $p\geq 2$.  Analogous results will be proved in a bounded domain $\Omega$.  These results will be stated in Section \ref{statement}. For the purpose of this introduction we content ourselves with a statement of our main result in the case $p\geq 2$ 
and $\Omega = \mathbf{R}^n$:

\begin{thm}\label{introthm}Let $p\geq 2$.  Suppose that there exists a solution of (\ref{quas2wentire}). Then there exists a constant $c>0$, depending on $n$ and $p$, so that:
\begin{equation}\begin{split}\label{pgeq2entirelowbd}
u(x) \geq c  \int_0^{\infty}&\Bigl[\frac{1}{r^{n-p}}\exp\Big( c\int_{0}^{r}\Bigl(\frac{\sigma(B(x,s))}{s^{n-p}}\Bigl)^{1/(p-1)}\frac{ds}{s}\Bigl)\\
& \cdot \int_{B(x,r)}\exp\Bigl(c\int_{0}^{r} \frac{\sigma(B(z,s))}{s^{n-p}}\frac{ds}{s}\Bigl)d\omega(z)\Bigl]^{1/(p-1)}\frac{dr}{r}.
\end{split}\end{equation}
Conversely,  there exists a positive constant $C(n,p,c)>0$ so that if $\sigma$ satisfies (\ref{capstart}) with $C(\sigma)< C$, there is a solution of (\ref{quas2wentire}) such that:
\begin{equation}\begin{split}\label{pgeq2entireupbd}
u(x) \leq c_1  \int_0^{\infty}&\Bigl[\frac{1}{r^{n-p}}\exp\Big( c\int_{0}^{r}\Bigl(\frac{\sigma(B(x,s))}{s^{n-p}}\Bigl)^{1/(p-1)}\frac{ds}{s}\Bigl)\\
& \cdot \int_{B(x,r)}\exp\Bigl(c\int_{0}^{r} \Bigl(\frac{\sigma(B(z,s))}{s^{n-p}}\Bigl)^{1/(p-1)}\frac{ds}{s}\Bigl)d\omega(z)\Bigl]^{1/(p-1)}\frac{dr}{r},
\end{split}\end{equation}
for a positive constant $c_1 = c_1(n,p,c)>0$, provided the right-hand side of the preceding inequality is finite at a single point $x\in \mathbf{R}^n$ (for some choice of $c>0$). 
\end{thm}

\begin{rem}[Concerning the optimality of (\ref{pgeq2entirelowbd}) and (\ref{pgeq2entireupbd})]\label{pgeq2sharprem} For certain choices of $\omega$, either (\ref{pgeq2entirelowbd}) or (\ref{pgeq2entireupbd}) are sharp.  This was shown in discussion in the previous paragraph; indeed, the display  (\ref{pgeq2entirelowbd}) is sharp if $\omega = \delta_{x_0}$, and (\ref{pgeq2entireupbd}) is sharp if $\omega = \sigma$.
\end{rem}

There are many classes of $\omega$ and $\sigma$ which satisfy the theorem above.  In Remark \ref{wolffentire} below we consider a simple condition on $\sigma$ which ensures the existence of a solution $u$ for any measure $\omega$ with corresponding pointwise bound.  This condition in particular covers the work of \cite{AHBV} cited above.

In Section \ref{examples}, we consider examples of conditions on $\omega$ so that one can deduce the existence of a solution for any potential $\sigma$ satisfying the condition (\ref{capstart}) with small constant.  In particular we will focus on the three cases:
\begin{enumerate}
\item if $\omega$ is a weak $A_{\infty}$-weight (this includes `power weight' right hand sides),
\item if $\omega$ lies in $L^q$ for some $q>1$,
\item if $\omega$ lies in a suitable Morrey space.
\end{enumerate}
In these three cases, we will see that the bound (\ref{pgeq2entireupbd}) simplifies. 

\begin{rem}\label{wolffentire}  From Theorem \ref{introthm}, it follows that there exists a constant $C=C(n,p)>0$ so that if $p \ge 2$ and 
\begin{equation}\label{wolffbd}\int_0^{\infty}\Bigl(\frac{\sigma(B(x,r))}{r^{n-p}}\Bigl)^{1/(p-1)}\frac{dr}{r} \leq C, \text{ for all }x\in \mathbf{R}^n,
\end{equation}
then there exists a positive constant $c=c(n,p)>0$, along with a solution $u$ of (\ref{quas2wentire}) such that:
$$\frac{1}{c}\int_0^{\infty}\Bigl(\frac{\omega(B(x,r))}{r^{n-p}}\Bigl)^{1/(p-1)}\frac{dr}{r} \leq u(x)  \leq c\int_0^{\infty}\Bigl(\frac{\omega(B(x,r))}{r^{n-p}}\Bigl)^{1/(p-1)}\frac{dr}{r}. 
$$
The reader should note that the condition (\ref{wolffbd}) is satisfied whenever $\sigma \in L^q(\mathbf{R}^n)$ for $q>n/p$ (with small $L^q$ norm), and so the theorems presented in this paper recover the relevant results of \cite{AHBV} mentioned above.  An analogous statement holds in the case $1<p<2$, and also when $\mathbf{R}^n$ is replaced by a bounded domain $\Omega$, as we will see in Section \ref{statement}.
\end{rem}

The preceding remarks explain that in a certain sense, the bounds of this paper are optimal.  However, the following question remains:
\begin{ques}\label{openques} Find a matching bilateral pointwise bound for solutions of (\ref{quas2w}) which is sharp for each measure $\omega$.  \end{ques}

Answering this question would be tantamount to inverting the nonlinear operator $-\Delta_p u -\sigma u^{p-1}$ pointwise.  Such a bound must necessarily have a much more nonlinear dependence on $\omega$.

The bounds (\ref{pgeq2entirelowbd}) and (\ref{pgeq2entireupbd}) are proved by studying certain \textit{local} function spaces whose underlying measure is $\sigma$, the measure appearing in the lower order term in (\ref{quas2w}).
The proof of the lower bound (\ref{pgeq2entirelowbd}) relies on a localisation procedure, see Section \ref{localise} below.  The lower bounds are then proved in Proposition \ref{lowbdprop}.

Our approach to proving the existence of a solution to (\ref{quas2w}) \textit{with the corresponding global upper bound} (\ref{pgeq2entireupbd}) goes via the construction of solutions to certain nonlinear integral inequalities, see Section \ref{construct} below.  It had been already seen in \cite{JVFund}, that such integral inequalities are intimately linked with solutions of (\ref{quas2w}).  It is these constructions which are the deepest portion of this paper.  With this integral supersolution in hand, we complete the proofs of our main results with an iterative argument, which is carried out in Section \ref{existence}.

\subsection{}  In Section \ref{1bvsection}, we study the equation (\ref{1bventire}), along with its counterpart in a bounded domain.  Our interest here is primarily to assert the bounds alluded to in display (\ref{1bventirebd}) above and the surrounding discussion.  It will be convenient to utilize a well known substitution (see e.g. \cite{Hi48, KK79, MP02, ADP06, AHBV} and Proposition \ref{logsubgen} below), to study the closely related quasilinear Riccati type equation:
\begin{equation}\begin{cases}\label{riccatiintro}
-\Delta_p v = (p-1)\abs{\nabla v}^p + \sigma \text{ in }\Omega,\\
v = 0\, \text{ on } \partial \Omega.
\end{cases}\end{equation}
There have been many recent papers devoted to studying such equations under a variety of assumptions on $\sigma$ and  $v$, see e.g. \cite{FM1, AHBV, Por02, MP02, GT03, ADP06, PS06} and references therein. In the process of asserting (\ref{1bventirebd}), we obtain the existence of solutions of (\ref{riccatiintro}) with pointwise bounds for general measures $\sigma$ (Theorem \ref{ricthm} below).  We thereby obtain an extension to quasilinear operators of work of Hansson, Maz'ya and Verbitsky \cite{HMV}, which complements the results in the aforementioned papers.  

\subsection{}  The plan of the paper is as follows.  In Section \ref{statement} we precisely state our main results. Sections \ref{background}--\ref{existence} are then devoted to proving our main results: Section \ref{background} introduces the required notation and background; Section \ref{construct} is then concerned with the construction of supersolutions to integral equations.  In Section \ref{lowbdsec}, we obtain lower bounds to solutions of (\ref{quas2w}).   The proofs are then concluded in Section \ref{existence}, where the constructions of Section \ref{construct} are used to prove the existence of solutions of (\ref{quas2w}) with corresponding bounds.

The final sections of the paper deal with applications and auxiliary results.  A study of solutions of the equation (\ref{1rhsentire}), along with their relationship to the equation (\ref{riccatiintro}) with natural growth in the gradient  is carried out in Section \ref{1bvsection}.  Section \ref{examples} is then devoted to special cases where our theorems are applicable.  Finally, in Section \ref{Hessian}, we consider fully nonlinear analogues of the Dirichlet problem (\ref{quas2w}) for Hessian equations with natural growth terms.

\section{Main results}\label{statement}

\subsection{} In this section we state our results.  We will use two notions of solution to study (\ref{quas2w}); the local notion of $p$-superharmonicity and the stronger global notion of solutions in the renormalizad sense, see Section \ref{renormsoln} below for a brief discussion.  

Let us first state the capacity condition on $\sigma$ that will appear throughout the paper.  For an open set $\Omega\subset \mathbf{R}^n$, if $u$ is a positive $p$-superharmonic solution of the inequality:
\begin{equation}\label{diffineq}-\Delta_p u \geq \sigma u^{p-1} \text{ in } \Omega,\end{equation} then, $\sigma$ obeys the following capacity condition:
\begin{equation}\label{capintro}
\sigma(E) \leq C \, \rm{cap}_{p}(E, \Omega) \;\text{ for any compact set } E\subset \Omega,
\end{equation}
with $C=1$.  This was proved as Lemma 4.3 in \cite{JVFund}\footnote{In \cite{JVFund}, two proofs of the necessity of (\ref{capintro}) are given.}.  Here $\text{cap}_p$ is the standard $p$-capacity associated to the Sobolev space $W^{1,p}(\Omega)$:
\begin{equation}\label{pcapintro}
\text{cap}_{p}(E, \Omega) = \inf\{ \;\norm{\nabla f}^p_{\text{L}^p} \;: \; f\geq 1 \; \text{on} \; E, \, \, f \in C^\infty_0(\Omega)\;\}.\end{equation}
Since any solution of (\ref{quas2w}) trivially satisfies (\ref{diffineq}), it follows that $\sigma$ satisfies (\ref{capintro}) whenever there exists a solution of (\ref{quas2w}).  \textit{Therefore, without loss of generality we impose that $\sigma$ satisfies (\ref{capintro}) throughout the paper}.

To prove the existence of solutions of (\ref{quas2w}) with corresponding upper bounds, we introduce a stronger condition, namely that:
\begin{equation}\label{strongcapcond}\sigma(E)\leq C(\sigma) \, \text{cap}_p(E)\, \text{ for all compact sets } E\subset \mathbf{R}^n,\end{equation}  
for a positive constant $C(\sigma)>0$.  Here $\text{cap}_p(E) = \text{cap}_p(E, \mathbf{R}^n)$ is the $p$-capacity in the entire space, i.e.
\begin{equation}\label{entirecap}\text{cap}_{p}(E) = \inf\{ \;\norm{\nabla f}^p_{\text{L}^p} \;: \; f\geq 1 \; \text{on} \; E, \, \, f \in C^\infty_0(\mathbf{R}^n)\;\}.
\end{equation}
It is immediate that $\text{cap}_p(E)\leq \text{cap}_p(E, \Omega)$, whenever $E\subset \Omega\subset \mathbf{R}^n$.   It is well known that (as a result of the Sobolev inequality) if $\sigma \in \textrm{L}^{\frac{n}{p}, \infty}(\mathbf{R}^n)$, then $\sigma$ satisfies (\ref{strongcapcond}).  However, much more general $\sigma$ are admissible for (\ref{strongcapcond}), possibly singular with respect to Lebesgue measure.  Note further that (\ref{capintro}) and (\ref{strongcapcond}) coincide when $\Omega = \mathbf{R}^n$, and so our results are sharpest in the entire space.

In order to be concise, we will use standard notation for the local potentials.  The fractional linear Riesz potential $\mathbf{I}^r_\alpha  (d \sigma)$, and the nonlinear Wolff potential $\mathbf{W}_{\beta,s}^r(d \sigma)$ are defined by:
\begin{equation}\label{Riesz}
\mathbf{I}^r_{\alpha}(d\sigma)(x) = \int_0^r \frac{\sigma(B(x,t))}{t^{n-\alpha}} \frac{dt}{t}, \text{ and}
\end{equation}
\begin{equation}\label{Wolff}
\mathbf{W}_{\beta ,s}^r(d\sigma)(x) =  \int_0^r \Bigl(\frac{\sigma(B(x,t))}{t^{n-\beta s}}\Bigl)^{1/(s-1)}\frac{dt}{t},
\end{equation}
respectively, where $0<r<\infty$, $0<\alpha<n$, $0<\beta < n/s$, and $1<s<\infty$. If $r=\infty$ then the superscript $r$ in the 
notations above will be dropped. 

 In the quasilinear 
case for equations of the $p$-Laplacian type we set: $\alpha=p$, $\beta=1$ and $s=p$. Recall that Theorem \ref{introthm} from the introduction concerned the equation (\ref{quas2w}) when $\Omega = \mathbf{R}^n$ in the case $2<p<n$, so we will next state our result when $\Omega = \mathbf{R}^n$ and $1<p<2$. 

\begin{thm}\label{mainthmleq2entire}  Let $1<p<2$, and suppose that $u$ is a solution of (\ref{quas2wentire}) in the $p$-superharmonic sense, then $\sigma$ satisfies (\ref{strongcapcond}).  In addition there is a constant $c=c(n,p)>0$ so that for all $x\in \mathbf{R}^n$,
\begin{equation}\begin{split}\label{pleq2entirelow}
u(x) \geq & c \int_{0}^{\infty}\Bigl(\frac{e^{c\mathbf{W}^r_{1,p}( d\sigma)(x)}}{r^{n-p}} \int_{B(x,r)} e^{c\mathbf{W}^r_{1,p}(d\sigma)(z)}d\omega(z)\Bigl)^{\frac{1}{p-1}}\frac{dr}{r}.
\end{split}
\end{equation}
On the other hand, under the assumption that the right hand side of (\ref{pleq2entireup}) is finite for some $x\in \mathbf{R}^n$ and $c>0$, there is a positive constant $C_0 = C_0(n,p,c)>0$, so that if $\sigma$ satisfies (\ref{strongcapcond}) with constant $C(\sigma)<C_0$, then there exists a solution $u$ of (\ref{quas2wentire}).  Furthermore, there is a positive constant $c_1=c_1(n,p,c)$ such that the constructed solution $u$ satisfies
\begin{equation}\label{pleq2entireup}
\begin{split}
u(x) \leq &c_1  \int_{0}^{\infty}\Bigl(\frac{e^{c\mathbf{W}^r_{1,p}(d\sigma)(x)}}{r^{n-p}} \int_{B(x,r)}e^{c\mathbf{I}^r_{p}(d\sigma)(z)}d\omega(z)\Bigl)^{\frac{1}{p-1}}\frac{dr}{r},
\end{split}
\end{equation}
for all $x\in \mathbf{R}^n$.
\end{thm}

The discussion in the introduction shows that bounds in Theorem \ref{mainthmleq2entire} are again optimal (cf. Remark \ref{pgeq2sharprem}).  In particular, display (\ref{pleq2entireup}) is sharp if $\omega = \delta_{x_0}$, and (\ref{pleq2entirelow}) is sharp if $\omega = \sigma$.

As in Remark \ref{wolffentire}, it follows that if $1<p<2$, and there exists a constant $C>0$ so that:
$$\int_0^{\infty}\frac{\sigma(B(x,r))}{r^{n-p}} \frac{dr}{r}\leq C \text{ for all }x\in \mathbf{R}^n,
$$
then there exists a positive constant $c=c(n,p)>0$ and an $p$-superharmonic solution of (\ref{quas2wentire}) so that:
$$\frac{1}{c} \int_0^{\infty} \Bigl(\frac{\omega(B(x,r))}{r^{n-p}}\Bigl)^{1/(p-1)} \frac{dr}{r}\leq u(x) \leq c \int_0^{\infty} \Bigl(\frac{\omega(B(x,r))}{r^{n-p}}\Bigl)^{1/(p-1)} \frac{dr}{r}.
$$
The corresponding statement continues to hold in bounded domains.  For examples of well known classes of $\omega$ where our theorems apply, see Section \ref{examples}.

\begin{rem}\label{tailrem}The condition used in the existence result above that (\ref{pleq2entireup}) is finite at a single point $x\in \mathbf{R}^n$ is equivalent to the same expression (\ref{pleq2entireup}) being finite almost everywhere in $\mathbf{R}^n$.  In fact, either statement follows from the following weaker \textit{tail estimate}:  \textit{there exists $x_0\in \mathbf{R}^n$ and $R>0$ so that}:
\begin{equation}\begin{split}\label{finitetailpgeq2}
 \int_R^{\infty}&\Bigl[\frac{1}{r^{n-p}}\exp\Big( c\int_{0}^{r}\Bigl(\frac{\sigma(B(x_0,s)\backslash B(x_0,R))}{s^{n-p}}\Bigl)^{1/(p-1)}\frac{ds}{s}\Bigl)\\
& \cdot \int_{B(x_0,r)}\exp\Bigl(c\int_{0}^{r} \Bigl(\frac{\sigma(B(z,s))}{s^{n-p}}\Bigl)\frac{ds}{s}\Bigl)d\omega(z)\Bigl]^{1/(p-1)}\frac{dr}{r}<\infty.
\end{split}\end{equation}
We discuss this further in Section \ref{finiteappendix}.  The analogous result is true when $p\geq 2$, and also in the case of bounded domains treated below.
\end{rem}

Let us now turn to our main results for the equation (\ref{quas2w}) in bounded domains $\Omega$.  Define $d(x) = \inf_{y\in \partial\Omega}\abs{x-y}$, to be the distance to the boundary of $\Omega$, and let $d_{\Omega}$ be the diameter of $\Omega$. 

\begin{thm}[Lower bounds]\label{bddlowbd} Let $1<p<n$.  Suppose that $u$ is a $p$-superharmonic solution of (\ref{quas2w}) in a bounded domain $\Omega$.  Then there is a constant $c=c(n,p)>0$, such that for all $x\in \Omega$: 

(i) if $1<p\leq 2$ then:
\begin{equation}\begin{split}\label{pleq2lowintrobd}
u(x) \geq & \,c \int_{0}^{\frac{d(x)}{5}}\Bigl(\frac{e^{c\mathbf{W}^r_{1,p}( d\sigma)(x)}}{r^{n-p}} \int_{B(x,r)} e^{c\mathbf{W}^r_{1,p}(d\sigma)(z)}d\omega(z)\Bigl)^{\frac{1}{p-1}}\frac{dr}{r}.
\end{split}
\end{equation}

(ii) if $2\leq p<n$ then:
\begin{equation}\begin{split}\label{pgeqlowbdst}
u(x) \geq c & \int_{0}^{\frac{d(x)}{5}}\Bigl(\frac{e^{c\mathbf{W}_{1,p}^r(d\sigma)(x)}}{r^{n-p}} \int_{B(x,r)} e^{c\mathbf{I}_{p}^r(d\sigma)(z)}d\omega(z)\Bigl)^{\frac{1}{p-1}}\frac{dr}{r}.
\end{split}
\end{equation}
\end{thm}

This theorem, along with the lower bounds in our previous stated results in the entire space, follow from Proposition \ref{lowbdprop} below.  Let us now turn to the existence of solutions with corresponding global upper bounds.

\begin{thm}[Existence and upper bounds]\label{bddupbd}Let $1<p<n$. Suppose that:
\begin{itemize}\item $1<p\leq 2$ and (\ref{pleq2upintrobd}) below is finite for some $x\in \Omega$ and $c>0$, or
\item $2\leq p<n$ and (\ref{pgequpbdst}) below is finite for some $x\in \Omega$ and $c>0$.
\end{itemize}
Then there is a positive constant $C_0=C_0(n,p,c)$ so that if $\sigma$ satisfies (\ref{capstart}) with $C(\sigma)<C_0$, then there exists a renormalized solution $u$ of (\ref{quas2w}) in $\Omega$ satisfying the following pointwise estimate for $x\in \Omega$:

(i) if $1<p\leq 2$:\begin{equation}\label{pleq2upintrobd}
u(x) \leq  c_1 \int_{0}^{2 d_{\Omega}}\Bigl(\frac{e^{c\mathbf{W}^r_{1,p}(\chi_{\Omega} d\sigma)(x)}}{r^{n-p}} \int_{B(x,r) \cap \Omega}e^{c\mathbf{I}^r_{p}(\chi_{\Omega}d\sigma)(z)}d\omega(z)\Bigl)^{\frac{1}{p-1}}\frac{dr}{r},
\end{equation}

(ii) if $2\leq p<n$:
\begin{equation}
\begin{split}\label{pgequpbdst}
u(x) \leq c_1 & \int_{0}^{2d_{\Omega}}\Bigl(\frac{e^{c\mathbf{W}_{1,p}^r(\chi_{\Omega} d\sigma)(x)}}{r^{n-p}} \int_{B(x,r) \cap \Omega}e^{c\mathbf{W}_{1,p}^r(\chi_{\Omega}d\sigma)(z)}d\omega(z)\Bigl)^{\frac{1}{p-1}}\frac{dr}{r}.
\end{split}
\end{equation}
Here $c_1 = c_1(c,n,p)>0$.
\end{thm}

These theorems remain to be optimal by looking at the special cases where $\omega = \delta_{x_0}$ or $\omega = \sigma$, as in our theorems stated in the entire space.

\section{Preliminaries}\label{background}

\subsection{Notation}\label{notation}
For an open set $\Omega$, and a measure $\sigma$ defined on $\Omega$, we let $L^p(\Omega, d\sigma)$ (or $L^p_{\text{loc}}(\Omega, d\sigma)$) to be the space of functions integrable (or locally integrable) to the $p$-th power with respect to the measure $\sigma$.  When $\sigma$ is Lebesgue measure, we instead write $L^p(\Omega)$ (or $L^p_{\text{loc}}(\Omega)$).

For a measure defined on $\mathbf{R}^n$, the mixed norm space $L^p(\ell^q,d\sigma)$ is defined as the space of sequences of functions $\{f_Q\}_{Q\in \mathcal{Q}}$ so that:
\begin{equation}\label{mixednormdef}\big|\big| f\big|\big|_{\displaystyle  L^p(\ell^q, d\sigma)} = \Bigl(\int_{\Omega}\Bigl[ \sum_{Q\in \mathcal{Q}} |f_Q(x)|^q \Bigl]^{p/q}d\sigma(x) \Bigl)^{1/p}<\infty.
\end{equation}
Here $\mathcal{Q}$ is the lattice of dyadic cubes in $\mathbf{R}^n$, see Section \ref{disccarlsec}.

We define the Sobolev space $W^{1,p}(\Omega)$ (respectively $W^{1,p}_{\text{loc}}(\Omega)$) to be the space of functions $u$ so that $u\in L^p(\Omega)$ and $|\nabla u| \in L^p(\Omega)$ (respectively $u\in L^p_{\text{loc}}(\Omega)$ and $|\nabla u|\in L^p_{\text{loc}}(\Omega)$).

Throughout this paper we use the display $A\lesssim B$, to mean $A\leq CB$, with $C$ a positive constant depending on the relevant allowed parameters of the particular theorem or lemma being proved.  For a $\sigma$ measurable set $E$, we will often denote by $|E|_{\sigma} = \sigma(E)$, the $\sigma$ measure of $E$.  Finally, for a dyadic cube $P\in \mathcal{Q}$ (see Section \ref{disccarlsec}), the display:
$$\sum_{Q\subset P}\; \text{  reads as ``\textit{the sum over all dyadic cubes $Q$ which are contained in $P$}".}$$

\subsection{Notions of solution} \label{renormsoln} Let $\mu$ be a nonnegative measure defined on a (possibly unbounded) domain $\Omega$, and extend $\mu$ to be $0$ outside $\Omega$ so that the resulting measure is defined on $\mathbf{R}^n$.  In this section we introduce two of the notions of solution for quasilinear equations with measure data, i.e. the Dirichlet problem:
\begin{equation}\label{renorm}
\begin{cases}\,-\Delta_p u = \mu\; \text{ in }\Omega,\\
\, u=0\;\; \text{ on } \partial\Omega.
\end{cases}
\end{equation}
It is well known, see e.g. \cite{K02}, that there are several notions of solution independently developed to study (\ref{renorm}).  We will briefly discuss two of them; renormalized solutions, and $p$-superharmonic solutions.

\textbf{Superharmonic solutions.}  We say that $u: \Omega \rightarrow (-\infty, \infty]$ is $p$-superharmonic if $u$ is a lower semicontinuous function, not identically infinite in any component of $\Omega$, and satisfying the following comparison principle:  whenever $D\subset\subset \Omega$ and $h \in C( \bar{D})$ is $p$-harmonic in $D$, with $h \leq u$ on $\partial D$, then $h \leq u$ in $D$.

It is well known (see, for instance \cite{HKM}), that for each $p$-superharmonic function in $\Omega$ we can associate a measure $\mu[u]$ on $\Omega$.  We then say that $-\Delta_p u = \mu$ in $\Omega$ \textit{in the $p$-superharmonic sense}, if $u$ is $p$-superharmonic in $\Omega$, and $\mu[u] = \mu$ in the sense of distributions.  In particular:
\begin{defn} We define $u\geq 0$ to be a solution of $-\Delta_pu = \sigma 
u^{p-1}+\omega$ in the $p$-superharmonic sense if $u\in L^{p-1}_{\text{loc}}(\Omega, d\sigma)$ and $d\mu[u] = u^{p-1}d\sigma +d\omega$.
\end{defn}

\textbf{Renormalized solutions.}  When working in a bounded domain $\Omega$, the more refined notion of \textit{renormalized solutions}, introduced by P. L. Lions and F. Murat, is often most convenient, see \cite{DMMOP99} for a comprehensive introduction.  Given a \textit{finite} nonnegative measure $\mu$, it is well known that we can decompose it as $\mu = \mu_0 + \mu_s$, where $\mu_0$ is absolutely continuous with respect to $\text{cap}_p$, and $\mu_s^{+}$ is singular with respect to $\text{cap}_p$.  We say that $u$ is a \textit{renormalized solution} of (\ref{renorm}) if:
\begin{equation}\label{intasump}\begin{split}T_k(u) \in & W_{0}^{1,p}(\Omega), \text{ for all } k>0; \; u \in \textrm{L}^{(p-1)\frac{n}{n-p}, \infty}(\Omega); \\
& \text{ and } \nabla u \in \textrm{L}^{(p-1)\frac{n}{n-1}, \infty}(\Omega).\end{split}\end{equation}
In addition, for all Lipschitz functions $h\in W^{1, \infty}(\mathbf{R})$ so that its derivative $h^\prime$ has compact support we have:
\begin{equation}\begin{split}
\int_{\Omega}\abs{\nabla u}^{p}h^\prime(u)\phi dx \,+ & \int_{\Omega}\abs{\nabla u}^{p-2}\nabla u\cdot \nabla \phi h(u) dx  =  \int_{\omega} h(u) \phi d\mu_0 \\
& + h(\infty)\int_{\Omega} \phi d\mu_s,
\end{split}
\end{equation}
whenever $\phi \in W^{1,r}(\Omega)\cap L^{\infty}(\Omega)$ with $r>n$, and so that $h(u)\phi \in W^{1,p}_0(\Omega)$.  Here $h(\infty) = \lim_{t\rightarrow\infty}h(t)$.  In particular:
\begin{defn}  A function $u$ satisfying (\ref{intasump}) is a renormalized solution of (\ref{quas2w}) if $u\in L^{p-1}_{\text{loc}}(\Omega, d\sigma)$, and for any $h$ as above:
\begin{equation}\begin{split}
\int_{\Omega}\abs{\nabla u}^{p}h^\prime(u)\phi dx \,+ & \int_{\Omega}\abs{\nabla u}^{p-2}\nabla u\cdot \nabla \phi h(u) dx  =  \int_{\Omega} |u|^{p-2}uh(u)\phi d\sigma\\
&+\int_{\omega} h(u) \phi d\mu_0 + h(\infty)\int_{\Omega} \phi d\mu_s,
\end{split}
\end{equation}
whenever $\phi \in W^{1,r}(\Omega)\cap L^{\infty}(\Omega)$ with $r>n$, and so that $h(u)\phi \in W^{1,p}_0(\Omega)$.
\end{defn}

The class of renormalized solutions is narrower than $p$-superharmonic solutions.  The wider class of test functions prevents  the existence of solutions which are influenced by a measure supported \textit{on the boundary of the domain}, as in some well known counterexamples for uniqueness of $p$-superharmonic functions (see \cite{K02}).  This additional property yields a \textit{global} potential estimate.  Recall the definition of the Wolff potential from (\ref{Wolff}) with $\beta=1$ and $s=p$. 
\begin{thm}[Potential estimates]\label{thmpotest}  Let $\mu$ be a nonnegative finite measure in $\Omega$.  There exists a positive constant $C_1 = C_1(n,p)$ such that the following two statements hold:

\indent a).  Suppose $B(x, 2r)\subset \Omega$, and suppose that $u$ satisfies $-\Delta_p u =\mu$ in $B(x, 2r)$ in the $p$-superharmonic sense, then:
\begin{equation}\label{potestloc}
u(x) \geq \frac{1}{C_1}\mathbf{W}_{1,p}^{r}(d\mu)(x).
\end{equation}
\indent b).  Suppose in addition that $\Omega$ is a bounded domain, and suppose that $u$ is a \textit{renormalized solution} of (\ref{renorm}), then there is a positive constant $C_1 = C_1(n,p)$ such that:
\begin{equation}\label{potest}
\frac{1}{C_1}\mathbf{W}_{1,p}^{d(x)/2}(d\mu)(x) \leq u(x) \leq C_1 \mathbf{W}_{1,p}^{2\text{diam}(\Omega)}(d\mu)(x).
\end{equation}
\end{thm}
\begin{proof}Part a) and the lower bound of part b) in Theorem \ref{thmpotest} are due to Kilpel\"{a}inen and Mal\'y \cite{KM92}. The upper bound is a global version of a local estimate of Kilpel\"ainen and 
Mal\'y \cite{KM1} obtained in \cite{PV}.\end{proof}

\subsection{Consequences of the capacity condition.} \label{necconds} 


Let $\sigma$ be a measure defined on $\mathbf{R}^n$.  In what follows we will need several consequences of the capacity condition (\ref{strongcapcond}), i.e. that there exists a constant $C(\sigma)>0$ so that:
$$\sigma(E)\leq C(\sigma)\text{cap}_p(E) \text{ for all compact sets }E\subset \mathbf{R}^n.
$$

The following result is a well-known theorem, due to work of Maz'ya, D. Adams and B. Dahlberg, which recasts the capacity condition as a multiplier condition:

\begin{thm}\cite{MazSob}\label{multcond}  Suppose that a nonnegative measure satisfies (\ref{strongcapcond}) with constant $C(\sigma)$, then the following inequality holds:
\begin{equation}\label{multcondst}
\int_{\mathbf{R}^n} |h|^p d\sigma \leq C(\sigma)\Bigl(\frac{p}{p-1}\Bigl)^p \int_{\mathbf{R}^n} |\nabla h|^p dx, \text{ for all } h\in C^{\infty}_0(\mathbf{R}^n). 
\end{equation}
\end{thm}
Let us next note that there exists a positive constant $c_{n,p}$, depending on $n$ and $p$, so that if $\sigma$ satisfies (\ref{strongcapcond}), then for each ball $B(x,r)\subset\mathbf{R}^n$:
\begin{equation}\label{ballest}
\sigma(B(x,r)) \leq C(\sigma)c_{n,p} r^{n-p}.
\end{equation}
Display (\ref{ballest}) is a special case of (\ref{capintro}).  It follows from the elementary fact that:
$$\text{cap}_p(B(x,r)) = \text{cap}_p(B(0,1))\cdot r^{n-p}.
$$

In what follows we will need several technical lemmas from \cite{JVFund}.  We first quote Corollary 4.11 in \cite{JVFund}, which is a John-Nirenberg type result:
\begin{lem} \label{corexpest} Suppose that $\sigma$ satisfies (\ref{strongcapcond}).  Then whenever $\beta\cdot C(\sigma)< 1$, there is a constant $C(n,p)>0$ such that:
\begin{equation}\label{expest}\int_{E}e^{\;\beta \mathbf{W}_{1,p}(\chi_{E}d\sigma)(y)} d\sigma(y) \leq \frac{C(n,p)}{1-C(\sigma)\beta} \sigma(E), \text{ for all compact sets } E\subset \mathbf{R}^n.
\end{equation}
\end{lem}
The second result shows that the Hausdorff measure condition (\ref{ballest}) gives us control over the tail of the Wolff potential.  We present the proof in Appendix \ref{appendix1} below, in slightly more generality, since the proof had been deferred from our previous paper \cite{JVFund}.
\begin{lem}\label{lemtailest}  Let $\sigma$ be satisfy the ball condition (\ref{ballest}). Then there is a positive constant $C=C(n, p, C(\sigma))>0$, so that for all $x \in \mathbb{R}^n$ and $y \in B(x,t)$,  $t>0$, it follows: 
\begin{equation}\label{tailest} \abs{\int_t^{\infty} \left [ \Bigl(\frac{\sigma(B(x,r))}{r^{n-p}}\Bigr)^{\frac{1}{p-1}} - \Bigl(\frac{\sigma(B(y,r))}{r^{n-p}}\Bigr)^{\frac{1}{p-1}} \right ]\frac{dr}{r}} \leq C .
\end{equation}
\end{lem}

Our second result using the ball estimate is a weighted exponential integrability result from the \cite{HJ11}.  The class of weights  is the so-called \textit{weak $A_{\infty}$} class, which we define now:
\begin{defn}  A nonnegative function $w$ is a weak $A_{\infty}$ weight, if there are constants $C_w>0$ and $\theta>0$ so that, for all balls $B$, and measurable subsets $E\subset B$:
\begin{equation}\label{weakainfdef}
\frac{|E|_w}{|2B|_w} \leq C_w\Big(\frac{|E|}{|Q|}\Bigl)^{\theta}.
\end{equation}
\end{defn}

\begin{thm}\label{morremb}\cite{HJ11}  Let $\omega$ be a weak $A_{\infty}$ weight, and suppose $\sigma$ is a measure such that (\ref{ballest}) holds.  Then for any $0<q<\infty$, there exist finite positive constants $C, c>0$ depending on $n,p, C(\sigma), q$; along with the constants $\theta$ and $C_{\omega}$ from (\ref{weakainfdef}), so that:
\begin{equation}\label{expint}
\frac{1}{|2B|_{\omega}}\int_{2B} \exp\Bigl[c\int_0^{\infty}\Bigl(\frac{\sigma(B(x,s)\cap B)}{s^{n-p}}\Bigl)^q\frac{ds}{s}\Bigl] d\omega \leq C,
\end{equation}
for all balls $B\subset\mathbf{R}^n$.
\end{thm}

\subsection{Dyadic Carleson embedding theorem}\label{disccarlsec}  In this section we briefly discuss the dyadic Carleson measure theorem, which we employ several times.  First, recall that a cube $Q$ in $\mathbf{R}^n$ is a \textit{dyadic cube} (we will write $Q\in \mathcal{Q}$), if $Q$ can be written $Q = [2^k m,2^{k+1}m)^n$ for some $k, m \in \mathds{Z}$.  We denote by $\ell(Q)$ the sidelength of the cube $Q$.  The reader should note that if $P$ and $Q$ are two dyadic cubes with non-empty intersection, then either $Q\subset P$, or $P\subset Q$.

It is known, see \cite{JVFund} Lemma 4.7, that condition (\ref{strongcapcond}) is equivalent to the existence of a constant $\tilde{C}>0$, so that $\sigma$ satisfies the following \textit{discrete Carleson measure} condition:
 \begin{equation} \label{discarla}
 \sum_{Q\subset P}c_Q \abs{Q}_{\sigma}^{p'} \leq \tilde{C} \, \abs{P}_{\sigma}\;\; \text{ for all } P\in \mathcal{Q}.
\end{equation}
Here, and throughout this paper, the sequence $c_Q$ is defined:
\begin{equation}\label{cQ} c_Q = \ell(Q)^{\frac{p-n}{p-1}}
.\end{equation}
Furthermore, there exists a constant $A=A(n,p)>0$ such that: $\tilde{C}/A \leq C(\sigma) \leq A \tilde{C}$, with $C(\sigma)$  as in (\ref{strongcapcond}).
The \textit{dyadic Carleson embedding theorem} (see, e.g.,  \cite{NTV99, COV}), is then: 
\begin{thm}\label{Carlemb}  Suppose that $\sigma$ satisfies (\ref{discarla}) with constant $\tilde{C}>0$.  Then, for any $s>1$,
 \begin{equation} \label{discarlb}
\sum_{Q \in \mathcal{Q}} c_Q\abs{Q}_{\sigma}^{p'} \left \vert \frac{1}{\abs{Q}_{\sigma}} \int_{Q} f d \sigma \right \vert^{s} \le  \tilde{C}\Bigl(\frac{s}{s-1}\Bigl)^s \, ||f||_{L^{s}(d \sigma)}^{s},
\end{equation}
for every $f \in L^{s}(d \sigma)$.
\end{thm}



\subsection{Dyadic shifting lemma}  Here we describe a tool to transfer results for dyadic potential operators to their continuous analogues.  This technique goes back to the seminal papers \cite{FS, GJ}.   Define $\mathcal{Q}_t$ to be the shifted dyadic lattice by $t\in \mathbf{R}^n$, i.e. $\mathcal{Q}_t = \{ \; Q + t \; : \, Q\in \mathcal{Q}\}$.
\begin{lem}\label{dyshift}  Let $\phi$, and $\psi$ be two functions mapping measurable sets into nonnegative measurable functions, so that whenever $A \subset B$ are two measurable sets, and $x\in \mathbf{R}^n$, it follows that $\phi(A)(x) \leq \phi(B)(x)$, and $\psi(A)(x)\leq \psi(B)(x)$.

Then, there exists $j_0 = j_0(n) \in \mathds{N}$, and $C=C(n,p)>0$, such that for all $k\in \mathds{Z}$, and $x\in \mathbf{R}^n$:
\begin{equation}\begin{split}\nonumber
\int_0^{2^k} & \Bigl(\frac{\phi(B(x,r))(x))}{r^{n-p}}\int_{B(x,r)}\psi(B(x,r))(z)d\omega(z)\Bigl)^{1/(p-1)}\frac{dr}{r}\\
& \leq \frac{C}{\abs{B(0, 2^{k+j_0})}}\int_{B(0, 2^{k+j_0})}\sum_{\substack{x\in Q_t, Q_t\in \mathcal{Q}_t \\ \ell(Q_t) \leq 2^{k+j_0}}}c_{Q_t}\Bigl(\phi(Q_t)(x)\int_{Q_t}\psi(Q_t)(z)d\omega(z)\Bigl)^{1/(p-1)}dt
\end{split}\end{equation}
\end{lem}
This lemma has the same proof as the standard dyadic shift argument, for instance see \cite{COV1}, p. 399.  We will use this result with $\psi$ and $\phi$ certain exponential weights, which will clearly satisfy the hypothesis of the Lemma.

\section{A nonlinear integral obstacle problem}\label{construct}

In this section we will construct solutions to certain nonlinear integral inequalities. It is the principal analytic argument in our existence theorems with corresponding upper bounds.

Let $\sigma$ be a nonnegative measure satisfying:
\begin{equation}\label{suffcapcond}
\sigma(E) \leq C(\sigma) \,  \text{cap}_p(E), \text{ for all compact sets } E\subset \mathbf{R}^n.
\end{equation}
where the capacity $\text{cap}_p(E)$ is defined as in (\ref{entirecap}).  Here the constant $C(\sigma)$ is reserved to be the least constant so that (\ref{suffcapcond}) holds.

Consider the nonlinear integral operator $\mathcal{T}$, acting on nonnegative functions $f\geq 0$, $f\in \textrm{L}^{p-1}_{\text{loc}}(\mathbf{R}^n, d\sigma)$ by:
\begin{equation}
\label{Ndefnsuper}
\mathcal{T}(f)(x) = \mathbf{W}_{1,p}(f^{p-1}d\sigma)(x) = \int_0^{\infty}\Bigl(\frac{1}{r^{n-p}}\int_{B(x,r)} f^{p-1}(z)d\sigma(z)\Bigl)^{1/(p-1)}\frac{dr}{r}.
\end{equation}

This section is devoted to the following problem: \textit{for a finite positive measure $\omega$, find a positive function $v$ so that $v\in \textrm{L}^{p-1}_{\text{loc}}(d\sigma)$ so that:}
\begin{equation}\label{obstacle}
\begin{cases} \, v(x) \geq \mathbf{W}_{1,p}(\omega)(x) \text{ and,}\\
\text{ there exists } C>0 \text{ so that } \mathcal{T}(v) \leq C (v - \mathbf{W}_{1,p}(\omega)).
 \end{cases}
\end{equation}
Solutions of (\ref{obstacle}) are solutions of a nonlinear obstacle problem for the integral operator $\mathcal{T}$, with obstacle $\mathbf{W}_{1,p}(\omega)$.  Given a solution of (\ref{obstacle}), a simple weak continuity argument shows the existence of solutions to (\ref{quas2w}), as we will show in Section \ref{existence}.

We will present a solution of (\ref{obstacle}) under the assumption that $C(\sigma)$ is small enough. The function $v$, as well as the argument to prove (\ref{obstacle}), will differ in the cases $1<p\leq2$, and $p\geq 2$.  Consider the function $v$, defined by:
\begin{equation}\label{vdefn}
v(x) = \int_0^{\infty} \Bigl(\frac{e^{\beta\mathbf{W}_{1,p}(\chi_{B(x,r)}d\sigma)(x)}}{r^{n-p}}\int_{B(x,r)}e^{\beta\mathbf{V}_{B(x,r)}(z)}d\omega(z)\Bigl)^{1/(p-1)}\frac{dr}{r}, 
\end{equation}
here:
\begin{equation}\label{localpot}
\mathbf{V}_{B(x,r)}(y) = \begin{cases}
\displaystyle \sum_{j=-\infty}^{\infty}\Bigl(\frac{\sigma(B(y, r/2^j)\cap B(x,r))}{(r/2^j)^{n-p}}\Bigl)^{1/(p-1)}, \, \text{ if } p\geq2,\\
\displaystyle \sum_{j=0}^{\infty}\frac{\sigma(B(y, r/2^j)\cap B(x,r))}{(r/2^j)^{n-p}},\, \text{ if } 1<p\leq2.
\end{cases}
\end{equation}

Using (\ref{ballest}), the following inequalities hold:
\begin{equation}\label{easyVbd}\mathbf{V}_{B(x,r)}(y) \lesssim \begin{cases}
\; \mathbf{W}^r_{1,p}(\chi_{B(x,r)}d\sigma)(y) + 1 \; \text{ if }p\geq 2,\\
\;\, \mathbf{I}^r_{p}(\chi_{B(x,r)}d\sigma)(y) + 1\; \text{ if }1<p\leq 2.\\
\end{cases}\end{equation}
As a result of (\ref{easyVbd}),  $v$ is less than a constant multiple of the right hand side of the bounds (\ref{pleq2upintrobd}) and (\ref{pgequpbdst}), for a suitably chosen constant $c>0$. 
Our primary result is the following:
\begin{thm}\label{superthm}  Suppose that there exists $\beta>0$ so that $v$ is finite almost everywhere.  Then there is a constant $C_0 = C_0(n,p)>0$, such that if $\beta C(\sigma)<C_0$, then $v$ is a solution of (\ref{obstacle}).\end{thm}

We will show in Section \ref{finiteappendix} below that $v<\infty$ almost everywhere for some $\beta>0$ if and only if it is finite at a single point.

The particular (and slightly cumbersome) structure of $\mathbf{V}$ in (\ref{localpot}) is in order to obtain a clean iteration argument.  One should keep in mind in what follows that the exponent $\beta$ in the exponential weights appearing in (\ref{vdefn}) have to remain constant in the iteration.

\begin{proof}[Proof of Theorem \ref{superthm}]To begin the proof, write:
$$\mathbf{W}_{1,p}(v^{p-1}d\sigma)(x) \lesssim I_{\text{out}} + I_{\text{in}}
$$
where:
\begin{equation}\begin{split}I_{\text{out}} = \int_0^{\infty} \Bigl(\frac{1}{r^{n-p}}\int_{B(x,r)}\Bigl\{& \int_r^{\infty} \Bigl(\frac{e^{\beta\mathbf{W}_{1,p}(\chi_{B(y,s)}d\sigma)(y)}}{s^{n-p}}\\
&\cdot \mu_{B(y,s)}\Bigl)^{1/(p-1)}\frac{ds}{s}\Bigl\}^{p-1}d\sigma(y)\Bigl)^{1/(p-1)}\frac{ds}{s},\end{split}\end{equation}
and
\begin{equation}\begin{split}I_{\text{in}}= \int_0^{\infty} \Bigl(\frac{1}{r^{n-p}}\int_{B(x,r)}\Bigl\{ & \int_0^{r} \Bigl(\frac{e^{\beta\mathbf{W}_{1,p}(\chi_{B(y,s)}d\sigma)(y)}}{s^{n-p}}\\
&\cdot \mu_{B(y,s)}\Bigl)^{1/(p-1)}\frac{ds}{s}\Bigl\}^{p-1}d\sigma(y)\Bigl)^{1/(p-1)}\frac{ds}{s},\end{split}\end{equation}
here:
\begin{equation}\label{muball}\mu_{B(y,s)} = \int_{B(y,s)}e^{\beta\mathbf{V}_{B(y,s)}(z)}d\omega(z).
\end{equation}
To prove Theorem \ref{superthm}, it clearly suffices to prove that there exists $C_0>0$ such that if $\beta C(\sigma) <C_0$ then it follows:
\begin{equation}\label{supsol1}
I_{\text{out}} \lesssim v - \mathbf{W}_{1,p}(\omega)\, \text{ and,}\end{equation}
\begin{equation}\label{supsol2}
I_{\text{in}} \lesssim v - \mathbf{W}_{1,p}(\omega).
\end{equation} 
We will first prove the inequality (\ref{supsol1}).  This inequality is \textit{responsible for the build-up of the tails of the potentials}, it is relatively simple and the proof is valid for all $1<p<n$.  The difficulty thus lies in inequality (\ref{supsol2}).  By inspection, (\ref{supsol2}) will immediately follow from the inequality:
\begin{equation}\begin{split}\label{localAinf} 
\int_{B(x,r)}\Bigl\{ \int_0^{r} \Bigl(\frac{e^{\beta\mathbf{W}_{1,p}(\chi_{B(y,s)}d\sigma)(y)}}{s^{n-p}}\cdot & \mu_{B(y,s)}\Bigl)^{1/(p-1)}\frac{ds}{s}\Bigl\}^{p-1}d\sigma(y) \\
& \lesssim \mu_{B(x,Ar)} - \omega(B(x, Ar)),
\end{split}
\end{equation}
for a constant $A>0$ depending on $n$ and $p$ (in our arguments $A$ will be at most $4$).  Here $\mu_{B(x,r)}$ is as in (\ref{muball}).

In particular, once (\ref{supsol1}) is proved, the problem of finding a solution of (\ref{obstacle}) is reduced to a local integral estimate.   We will prove that (\ref{supsol1}) and (\ref{supsol2}) in the following series of lemmas.\end{proof}
\begin{lem} There exists $C_0=C_0(n,p)>0$, so that if $C(\sigma)<C_0$, then (\ref{supsol1}) holds
\end{lem}
\begin{proof} For $y\in B(x,r)$ with $s>r$, it is clear that $B(x, 2s) \supset B(y,s)$.  Therefore:
\begin{equation}\begin{split}\nonumber I_{\text{out}} \leq \int_0^{\infty} \Bigl(\frac{1}{r^{n-p}}\int_{B(x,r)}\Bigl\{ \int_r^{\infty} \Bigl(&\frac{e^{\beta\mathbf{W}_{1,p}(\chi_{B(x,2s)}d\sigma)(y)}}{s^{n-p}}\\
& \cdot \mu_{B(x,2s)}\Bigl)^{1/(p-1)}\frac{ds}{s}\Bigl\}^{p-1}d\sigma(y)\Bigl)^{1/(p-1)}\frac{ds}{s}. \end{split}\end{equation}
Let us re-write this as:
$$I_{\text{out}} \leq \int_0^{\infty} \Bigl(\frac{1}{r^{n-p}}\Bigl)^{1/(p-1)}III\frac{dr}{r}
$$
with $III = III(x, r)$ defined by:
$$III = \Bigl(\int_{B(x,r)}\Bigl\{ \int_r^{\infty} \Bigl(\frac{e^{\beta\mathbf{W}_{1,p}(\chi_{B(x,2s)}d\sigma)(y)}}{s^{n-p}}\cdot \mu_{B(x,2s)}\Bigl)^{1/(p-1)}\frac{ds}{s}\Bigl\}^{p-1}d\sigma(y)\Bigl)^{1/(p-1)}.
$$
To estimate $III$, note that for any $y\in B(x,r)$, an application of Lemma \ref{lemtailest} yields:
\begin{equation}\label{tailestlem1} \begin{split} \exp\bigl(\beta\mathbf{W}_{1,p}(\chi_{B(x,2s)}&d\sigma)(y)\bigl) \leq C \exp\bigl(\beta \mathbf{W}_{1,p}(\chi_{B(x, 2r)}d\sigma)(y)\bigl)\\
& \cdot\exp\Bigl(\beta \int_r^{\infty} \Bigl(\frac{\sigma(B(x,t)\cap B(x, 2s))}{t^{n-p}}\Bigl)^{1/(p-1)}\frac{dt}{t}\Bigl).
\end{split}\end{equation}
Substituting display (\ref{tailestlem1}) into $III$ and integrating through, we estimate:
\begin{equation}\begin{split}III \lesssim  \int_r^{\infty} \Bigl[&\frac{\mu_{B(x, 2s)}}{s^{n-p}}\exp\Bigl(\beta \int_r^{\infty} \Bigl(\frac{\sigma(B(x,t)\cap B(x, 2s))}{t^{n-p}}\Bigl)^{1/(p-1)}\frac{dt}{t}\Bigl)\\
&\Bigl]^{1/(p-1)}\frac{ds}{s}\cdot\Bigl(\int_{B(x,2r)}e^{\beta(p-1) \mathbf{W}_{1,p}(\chi_{B(x, 2r)}d\sigma)(y)}d\sigma(y)\Bigl)^{1/(p-1)}.
\end{split}\end{equation}
Next, by the exponential integrability lemma (Lemma \ref{corexpest}), we assert that there is a constant $C_0$ such that if $\beta C(\sigma)<C_0$, then:
\begin{equation}\begin{split}\label{innerintest}III \lesssim  \int_r^{\infty} \Bigl[ &\frac{\mu_{B(x, 2s)}}{s^{n-p}}\exp\Bigl(\beta \int_r^{\infty} \Bigl(\frac{\sigma(B(x,t)\cap B(x, 2s))}{t^{n-p}}\Bigl)^{1/(p-1)}\frac{dt}{t}\Bigl)\\
&\Bigl]^{1/(p-1)}\frac{ds}{s}\cdot\sigma(B(x, 2r))^{1/(p-1)}.
\end{split}\end{equation}
Substituting the estimate (\ref{innerintest}) for $III$ into the definition of $I_{\text{out}}$, we find that:
\begin{equation}
\begin{split} I_{\text{out}} \lesssim & \int_0^{\infty}\Bigl(\frac{\sigma(B(x, 2r))}{r^{n-p}}\Bigl)^{1/(p-1)}\int_r^{\infty}\Bigl(\frac{\mu_{B(x, 2s)}}{s^{n-p}}\\
& \cdot \exp\Bigl(\beta \int_r^{\infty} \Bigl(\frac{\sigma(B(x,t)\cap B(x, 2s))}{t^{n-p}}\Bigl)^{1/(p-1)}\frac{dt}{t}\Bigl)\Bigl)^{1/(p-1)}\frac{ds}{s}\frac{dr}{r}.
\end{split}\end{equation}
Applying Fubini's theorem yields:
\begin{equation}\label{afterfubini1}
\begin{split} I_{\text{out}} \lesssim & \int_0^{\infty}\Bigl(\frac{\mu_{B(x, 2s)}}{s^{n-p}}\Bigl)^{1/(p-1)} \Bigl[\int_0^s\Bigl(\frac{\sigma(B(x, 2r))}{r^{n-p}}\Bigl)^{1/(p-1)}\\
& \cdot \exp\Bigl(\beta/(p-1) \int_r^{\infty} \Bigl(\frac{\sigma(B(x,t)\cap B(x, 2s))}{t^{n-p}}\Bigl)^{1/(p-1)}\frac{dt}{t}\Bigl)\frac{ds}{s}\Bigl]\frac{dr}{r}.
\end{split}\end{equation}
On the other hand, employing (\ref{suffcapcond}) in the form of (\ref{ballest}), it follows:
\begin{equation}\label{easycapexpest}\exp \Bigl(\beta/(p-1) \int_s^{\infty} \Bigl(\frac{\sigma(B(x,t)\cap B(x, 2s))}{t^{n-p}}\Bigl)^{1/(p-1)}\frac{dt}{t}\Bigl)\leq C(n,p,C(\sigma), \beta),
\end{equation}
note here the constant is independent of $r$.  Changing variables, and applying (\ref{easycapexpest}) in (\ref{afterfubini1}), it follows:
\begin{equation}\nonumber
\begin{split} I_{\text{out}} \lesssim & \int_0^{\infty}\Bigl(\frac{\mu_{B(x, s)}}{s^{n-p}}\Bigl)^{1/(p-1)} \int_0^s\Bigl(\frac{\sigma(B(x, r))}{r^{n-p}}\Bigl)^{1/(p-1)}\\
& \cdot \exp\Bigl(\beta/(p-1) \int_r^{s} \Bigl(\frac{\sigma(B(x,t)}{t^{n-p}}\Bigl)^{1/(p-1)}\frac{dt}{t}\Bigl)\frac{dr}{r}\frac{ds}{s}.
\end{split}\end{equation}
An elementary integration by parts argument now yields:
\begin{equation}\begin{split}\nonumber \int_0^s\Bigl(\frac{\sigma(B(x, r))}{r^{n-p}}\Bigl)^{1/(p-1)}&\exp\Bigl(\beta/(p-1) \int_r^{s} \Bigl(\frac{\sigma(B(x,t)}{t^{n-p}}\Bigl)^{1/(p-1)}\frac{dt}{t}\Bigl)\frac{dr}{r}\\
& \lesssim \Bigl\{\exp\Bigl(\beta/(p-1) \int_0^{s} \Bigl(\frac{\sigma(B(x,t)}{t^{n-p}}\Bigl)^{1/(p-1)}\frac{dt}{t}\Bigl) - 1\Bigl\}.
\end{split}\end{equation}
In conclusion:
\begin{equation}\nonumber
\begin{split} I_{\text{out}} \lesssim  \int_0^{\infty}&\Bigl(\frac{\mu_{B(x, s)}}{s^{n-p}}\Bigl)^{1/(p-1)}\cdot \Bigl\{\exp\Bigl[\beta/(p-1) \int_0^{\infty} \\
&\Bigl(\frac{\sigma(B(x,t)\cap B(x, s))}{t^{n-p}}\Bigl)^{1/(p-1)}\frac{dt}{t}\Bigl]-1\Bigl\}\frac{ds}{s}.
\end{split}\end{equation}
Since $\mu_{B(x,r)} \geq \omega(B(x,r))$, the inequality (\ref{supsol1}) follows.  \end{proof}

\subsection{The proof of (\ref{localAinf})}  It now follows from the discussion after display (\ref{supsol2}), that Theorem \ref{superthm} will be proved once we show that we can choose $C_0>0$ depending on $n$ and $p$, such that if $\beta C(\sigma)<C_0$, then (\ref{localAinf}) holds for all balls $B(x,r)$.  We first prove (\ref{localAinf}) when $1<p\leq 2$, indeed:

\begin{lem}\label{sleq2upbdlem}  Let $1<p\leq2$.  There exists $C_0 = C_0(n,p)>0$ such that if $\beta C(\sigma)<C_0$, with $C(\sigma)$ as in (\ref{suffcapcond}), then there is a constant $C=C(n,p)>0$ such that for all $B(x,r)$:
\begin{equation}\begin{split}\label{slequpbdlemst}
\int_{B(x,r)}\Bigl\{\int_0^r & e^{\beta\mathbf{W}_{1,p}(\chi_{B(x,r)}(y))}\Bigl(\frac{1}{t^{n-p}}\int_{B(y,t)}e^{\beta \mathbf{V}_{B(y,t)})(z)}d\omega(z)\Bigl)^{1/(p-1)}\\
&\Bigl\}^{p-1}d\sigma(y) \leq C \int_{B(x,2r)}\Bigl( e^{\beta \mathbf{V}_{B(x,2r)}(z)} -1\Bigl)d\omega(z).
\end{split}\end{equation}
\end{lem}
\begin{proof}  Without loss of generality, let $\omega\equiv 0 $ and $\sigma\equiv 0 $ on $\mathbf{R}^n \backslash B(x, 2r)$.  Then, note that by definition:
\begin{equation}\label{restricrem}\mathbf{V}_{B(x,2r)}(z) = \sum_{j=0}^{\infty}\frac{\sigma(B(z, 2r/2^j))}{(2r/2^j)^{n-p}}.
\end{equation}
We will display the local potential in (\ref{slequpbdlemst}) as a sum and then use a sequence space imbedding.  Indeed, the left hand side of (\ref{slequpbdlemst}) is less than a constant multiple of:
\begin{equation}\begin{split}\int_{B(x,r)} \Bigl\{\sum_{j=0}^{\infty} & e^{\beta\mathbf{W}_{1,p}(\chi_{B(y, r/2^j)}d\sigma)(y)}\Bigl[\bigl(\frac{r}{2^j}\bigl)^{p-n}\\
& \cdot \int_{B(y, r/2^j)} e^{\beta\mathbf{V}_{B(y, r/2^j)}(z)}d\omega(z) \Bigl]^{1/(p-1)}\Bigl\}^{p-1}d\sigma(y).
\end{split}\end{equation}
Since $1<p\leq 2$, it follows that $||\cdot||_{\ell^1}\leq||\cdot||_{\ell^{p-1}}$.  Hence the previous display is less than:
\begin{equation}\label{boringstep}\begin{split}\int_{B(x,r)} \Bigl[\sum_{j=0}^{\infty}& e^{(p-1)\beta\mathbf{W}_{1,p}(\chi_{B(y, r/2^j)}d\sigma)(y)}\bigl(\frac{r}{2^j}\bigl)^{p-n} \\
& \int_{B(y, r/2^j)} e^{\beta\mathbf{V}_{B(y, r/2^j)}(z)}d\omega(z)\Bigl]d\sigma(y).
\end{split}\end{equation}
By Fubini's theorem, and since $B(z, 2r/2^j)\supset B(y, r/2^j)$ for $z\in B(y, r/2^j)$, the display (\ref{boringstep}) is less than:
$$\int_{B(x, 2r)} \sum_{j=0}^{\infty} \bigl(\frac{r}{2^j}\bigl)^{p-n} \int_{B(z, r/2^j)}e^{\beta\mathbf{W}_{1,p}(\chi_{B(y, r/2^j)}d\sigma)(y)}d\sigma(y) e^{\beta \mathbf{V}_{B(z,2r/2^j)}(z)}d\omega(z).
$$
By the exponential integrability lemma (Lemma \ref{corexpest}), there exists $C_0>0$ such that if $\beta C(\sigma)< C_0$, this last line is less than a constant multiple of:
\begin{equation}\label{afterexpintpleq2}\int_{B(x, 2r)} \sum_{j=0}^{\infty} \bigl(\frac{r}{2^j}\bigl)^{p-n}\abs{B(z, r/2^{j-1})}_{\sigma} e^{\beta \mathbf{V}_{B(z,2r/2^j)}(z)}d\omega(z).
\end{equation}
In order to apply Lemma \ref{corexpest} to obtain (\ref{afterexpintpleq2}), we observe:
$$\int_{B(z, r/2^j)}e^{\beta\mathbf{W}_{1,p}(\chi_{B(y, r/2^j)}d\sigma)(y)}d\sigma(y)\leq \int_{B(z, r/2^{j-1})}e^{\beta\mathbf{W}_{1,p}(\chi_{B(z, r/2^{j-1})}d\sigma)(y)}d\sigma(y),
$$
from which a direct application of the exponential integrability lemma yields (\ref{afterexpintpleq2}).

The proof will be completed by summation by parts.  Indeed, the display (\ref{afterexpintpleq2}) is equal to:
\begin{equation}\label{b4sumbyparts2}\int_{B(x, 2r)} \sum_{j=0}^{\infty} \bigl(\frac{2r}{2^j}\bigl)^{p-n}\abs{B(z, 2r/2^{j})}_{\sigma}\exp\Bigl(\beta  \sum_{k=j}^{\infty}\frac{\sigma(B(z, 2r/2^k))}{(2r/2^k)^{n-p}}\Bigl)d\omega(z).
\end{equation}
Recall the following elementary summation by parts result (see e.g. \cite{FV2}):  Suppose that $\{\lambda_j\}_j$ is a nonnegative sequence such that $0\leq \lambda_j \leq 1$.  Then:
\begin{equation}\label{sumparts2}\sum_{j=0}^{\infty}\lambda_j e^{ \sum_{k=j}^{\infty}\lambda_k} \leq 2\,\Bigl( e^{\sum_{j=0}^{\infty}\lambda_j} -1\Bigl). 
\end{equation}
Note that (\ref{sumparts2}) can be applied in (\ref{b4sumbyparts2}) provided $C_0<1$, and hence the lemma follows by recalling the definition of $\mathbf{V}_{B(x,2r)}$ from (\ref{restricrem}).
\end{proof}
We shall now move onto $p>2$, which requires a more involved argument based on Theorem \ref{Carlemb}.  Recall the definition of the local potential $\mathbf{V}_{B(x,r)}$ from (\ref{localpot}).
\begin{lem}\label{sgeq2upbdlem2}  Let $ 2<p<n$.   There exists $C_0 = C_0(n,p)$, such that if $\beta C(\sigma)<C_0$, then there is a constant $C=C(n,p)>0$ and:
\begin{equation}\begin{split}\label{sgequpbdlem2st}
\int_{B(x,2^k)}\Bigl\{\int_0^{2^k} & e^{\beta\mathbf{W}_{1,p}(\chi_{B(y,t)}(y))}\Bigl(\frac{1}{t^{n-p}}\int_{B(y,t)}e^{\beta\mathbf{V}_{B(y,t)}(z)}d\omega(z)\Bigl)^{1/(p-1)}\Bigl\}^{p-1}d\sigma(y)\\
&  \leq C \int_{B(x,2^{k+1})} \Bigl(e^{\beta\mathbf{V}_{B(x,2^{k+1})}(z)}-1\Bigl)d\omega(z),  \text{ for all }B(x,2^k), \, k\in \mathds{Z} .
\end{split}\end{equation}
\end{lem}


\begin{proof}  Without loss of generality, suppose $\sigma\equiv 0 $ and $\omega \equiv 0$ on $\mathbf{R}^n\backslash B(x, 2^{k+1})$.  By the dyadic shifting lemma (Lemma \ref{dyshift}); if $\ell \geq k+k_0$, with $k_0>0$ depending on $n$, the left hand side of (\ref{sgequpbdlem2st}) is less than a constant multiple of:
\begin{equation}\begin{split}\label{upbdave1}
\int_{B(x,2^k)} \Bigl\{\frac{1}{2^{\ell n}}& \int_{B(0, 2^{\ell})} \sum_{\substack{y \in Q_t, Q_t\in \mathcal{Q}_t \\ \ell(Q_t) \leq 2^{j+j_0}} }  c_Q \, e^{\beta\mathbf{W}_{1,p}(\chi_{Q_t} d\sigma)(y)}\\
&\cdot \Bigl( \int_{Q_t}e^{\beta\mathbf{V}_{Q_t}(z)}d\omega(z)\Bigl)^{1/(p-1)}dt\Bigl\}^{p-1} d\sigma(y)
\end{split}
\end{equation}
here $c_Q$ is as in (\ref{cQ}), and $j_0$ depends only on $n$. 

Since $p>2$, it follows from Jensen's inequality that  (\ref{upbdave1}) is less than:
\begin{equation}\begin{split}\label{upbdave2}
\frac{1}{2^{\ell n}} \int_{B(0, 2^{\ell})}\Bigl[\int_{B(x,2^k)} \Bigl\{ & \sum_{\substack{y \in Q_t, Q_t\in \mathcal{Q}_t \\ \ell(Q_t) \leq 2^{j+j_0}} }  c_Q \, e^{\beta\mathbf{W}_{1,p}(\chi_{Q_t} d\sigma)(y)}\\
&\cdot \Bigl( \int_{Q_t}e^{\beta\mathbf{V}_{Q_t}(z)}d\omega(z)\Bigl)^{1/(p-1)}\Bigl\}^{p-1} \Bigl] d\sigma(y)dt
\end{split}
\end{equation}
It suffices to be able to estimate the inner integral in (\ref{upbdave2}) (the expression in square brackets), for a fixed $t$, with constant independent on $t$.  We will therefore assume $t=0$.  As a result, it suffices to prove the estimate:
\begin{equation}\label{beforeduality}\begin{split}\int_{B(x,2^k)}\Bigl[ & \sum_{\substack{y \in Q\\ \ell(Q) \leq 2^{k+j_0}}}  c_Q   e^{\beta\mathbf{W}_{1,p}(\chi_{Q} d\sigma)(y)}\Bigl( \int_{Q}  e^{\beta\mathbf{V}_{Q} (z)}d\omega(z)\Bigl)^{1/(p-1)}dt\\
&\Bigl]^{p-1} d\sigma(y) \leq C \int_{B(x,2^{k+1})} \Bigl(e^{\beta\mathbf{V}_{B(x,2^{k+1})}(z)}-1\Bigl)d\omega(z).
\end{split}\end{equation}
We first claim that the left hand side of (\ref{beforeduality}) is less than or equal a constant multiple of:
\begin{equation}\label{afterduality}I = \sum_{\substack{Q\in \mathcal{Q} \\ \ell(Q) \leq 2^{k+j_0}}}c_Q\abs{Q}_{\sigma}^{1/(p-1)} \int_Q e^{\beta\mathbf{V}_Q(z)} d\omega(z).
\end{equation}
Here we will use the Carleson measure theorem.  First, let $\lambda_Q = \int_Q e^{\beta\mathbf{V}_Q(z)} d\omega(z)$, then (\ref{beforeduality}) is equal to:
$$\int_{B(x,2^k)}\Bigl\{ \sum_{\substack{y \in Q\\ \ell(Q) \leq 2^{k+j_0}}}  c_Q  e^{\beta\mathbf{W}_{1,p}(\chi_{Q} d\sigma)(y)}\lambda_Q^{1/(p-1)}dt\Bigl\}^{p-1} d\sigma(y),
$$
Let $q = (p-1)' = (p-1)/(p-2)$.  Applying duality in $L^{p-1}(d\sigma)$, the left hand side of (\ref{beforeduality}) is equal to:
\begin{equation}\label{pgeq2upstep1}\sup_{\norm{g}_{L^{q}(d\sigma)}=1} \Bigl(\sum_{\substack{Q \in \mathcal{Q}\\ \ell(Q) \leq 2^{k+j_0}}}  c_Q \int_Q |g(y)| e^{\beta\mathbf{W}_{1,p}(\chi_{Q} d\sigma)(y)}d\sigma(y)\lambda_Q^{1/(p-1)}\Bigl)^{p-1}.
\end{equation}
Fix such a admissible $g\geq 0$, then it follows H\"{o}lder's inequality that the quantity in the supremum appearing in display (\ref{pgeq2upstep1}) is less than:
$$I\cdot\Bigl(\sum_{\substack{Q \in \mathcal{Q}\\ \ell(Q) \leq 2^{k+j_0}}}  c_Q \abs{Q}_{\sigma}^{p'}\Bigl(\frac{1}{\abs{Q}_{\sigma}}\int_Q |g(y)| e^{\beta\mathbf{W}_{1,p}(\chi_{Q} d\sigma)(y)}d\sigma(y)\Bigl)^{\frac{p-1}{p-2}}\Bigl)^{p-2},
$$
with $I$ as in (\ref{afterduality}).  To prove the claim, it therefore remains to prove that:
\begin{equation}\label{pgeq2upstep2}\begin{split}\sum_{\substack{Q \in \mathcal{Q}\\ \ell(Q) \leq 2^{k+j_0}}}  & c_Q \abs{Q}_{\sigma}^{p'}\Bigl(\frac{1}{\abs{Q}_{\sigma}}\int_Q |g(y)| e^{\beta\mathbf{W}_{1,p}(\chi_{Q} d\sigma)(y)}d\sigma(y)\Bigl)^{\frac{p-1}{p-2}}\\
& \lesssim \int |g|^{\frac{p-1}{p-2}}d\sigma, \text{ for any } g\in L^{(p-1)/(p-2)}(d\sigma).
\end{split}\end{equation}
To this end, let $\epsilon>0$ such that $(p-2)(1+\epsilon)<p-1$, and apply H\"{o}lder's inequality in the following way:
\begin{equation}\begin{split}\Bigl(\frac{1}{\abs{Q}_{\sigma}}\int_Q |g(y)| & e^{\beta\mathbf{W}_{1,p}(\chi_{Q} d\sigma)(y)}d\sigma(y)\Bigl)^{\frac{p-1}{p-2}}\leq \Bigl(\frac{1}{\abs{Q}_{\sigma}}\int_Q |g(y)|^{1+\epsilon}d\sigma(y)\Bigl)^{\frac{p-1}{(p-2)(1+\epsilon)}}\\
& \cdot \Bigl(\frac{1}{\abs{Q}_{\sigma}}\int_Qe^{\frac{1+\epsilon}{\epsilon} \beta\mathbf{W}_{1,p}(\chi_{Q} d\sigma)(y)}d\sigma(y)\Bigl)^{\frac{(p-1)\epsilon}{(p-2)(1+\epsilon)}}.
\end{split}\end{equation}
Now, note that, if $\frac{1+ \epsilon}{\epsilon}\beta C(\sigma) <C_0$, we may apply Lemma \ref{corexpest} to estimate:
$$\Bigl(\frac{1}{\abs{Q}_{\sigma}}\int_Qe^{\frac{1+\epsilon}{\epsilon} \beta\mathbf{W}_{1,p}(\chi_{Q} d\sigma)(y)}d\sigma(y)\Bigl)^{\frac{(p-1)\epsilon}{(p-2)(1+\epsilon)}}\leq C.
$$
The left hand side of (\ref{pgeq2upstep2}) is therefore less than a constant multiple of:
$$\sum_{\substack{Q \in \mathcal{Q}\\ \ell(Q) \leq 2^{k+j_0}}}  c_Q \abs{Q}_{\sigma}^{p'}\Bigl(\frac{1}{\abs{Q}_{\sigma}}\int_Q |g(y)|^{1+\epsilon}d\sigma(y)\Bigl)^{\frac{p-1}{(p-2)(1+\epsilon)}}.
$$
Since $\frac{p-1}{(p-2)(1+\epsilon)}>1$, applying the Carleson measure theorem (Theorem \ref{Carlemb}), with $s =\frac{p-1}{(p-2)(1+\epsilon)}$, it follows that conclude that (\ref{pgeq2upstep2}) holds.  Hence the claim is proved.

When specialised to a cube $Q$, the condition (\ref{suffcapcond}) is $\sigma(Q) \lesssim C(\sigma)\ell(Q)^{n-p}$ (this follows analogously to the condition (\ref{ballest})).   It follows for $z\in Q$ with $\ell(Q)\leq 2^{k+j_0}$, that there is a constant $C=C(n,p, C(\sigma))>0$, so that:
$$\mathbf{V}_Q(z) \leq \sum_{j=-\infty}^{\ell(Q)}\Bigl(\frac{\sigma(Q\cap B(z, 2^j)\cap B(x, 2^{k+1}))}{2^{j(n-p)}}\Bigl)^{1/(p-1)} + \, C
$$
Recall $\sigma \equiv 0$ on $\mathbf{R}^n\backslash\{0\}$.  After applying Fubini's Theorem in (\ref{afterduality}), we deduce from the previous display that (\ref{afterduality}) is less than a constant multiple of:
\begin{equation}\label{sumpartsprelim1}\begin{split}\int_{B(x,2^{k+1})}&\sum_{\substack{z\in Q\\ \ell(Q) \leq 2^{k+j_0}}} c_Q \abs{Q\cap B(x, 2^{k+1})}_{\sigma}^{1/(p-1)}\\
& \cdot \exp\Bigl(\beta\sum_{j=-\infty}^{\ell(Q)}\Bigl(\frac{\sigma(Q\cap B(z, 2^j)\cap B(x, 2^{k+1}))}{2^{j(n-p)}}\Bigl)^{1/(p-1)}\Bigl)d\omega(z).
\end{split}\end{equation}
By the summation by parts result (\ref{sumparts2}), valid under assumption that $\beta C(\sigma) <1$, it follows that the integrand in (\ref{sumpartsprelim1}) is less than a constant multiple of:
\begin{equation}\label{pgeq2aftersumparts}\exp\Bigl(\beta\sum_{j=-\infty}^{k+j_0+j_1}\Bigl(\frac{\sigma(B(z, 2^j)\cap B(x, 2^{k+1}))}{2^{j(n-p)}}\Bigl)^{1/(p-1)}\Bigl)-1.
\end{equation}
Indeed, each cube $Q$ so that $\ell(Q) = 2^{j}$, with $z\in Q$, is contained in a ball $B(z, 2^{j+j_1})$, where $j_1$ is a dimensional constant.  Thus, the integrand in (\ref{sumpartsprelim1}) is less than a constant multiple of:
\begin{equation}\begin{split}\nonumber2^{j_1 \frac{n-p}{p-1}}\sum_{\ell=-\infty}^{k+j_0+j_1}&\Bigl(\frac{\abs{B(z, 2^{\ell}) \cap B(x, 2^{k+1})}_{\sigma}}{2^{\ell(n-p)}}\Bigl)^{1/(p-1)}\\
&\cdot\exp\Bigl(\beta\sum_{j=-\infty}^{\ell}\Bigl(\frac{\sigma(B(z, 2^j)\cap B(x, 2^{k+1}))}{2^{j(n-p)}}\Bigl)^{1/(p-1)}\Bigl),
\end{split}\end{equation} 
from which we may apply (\ref{sumparts2}) to conclude (\ref{pgeq2aftersumparts}).  Here we have used the definition of $c_Q$ from (\ref{cQ}).

Recalling the definition of $\mathbf{V}_{B(x, 2^{k+1})}$, we have asserted that there exists a constant $C=C(n,p,C(\sigma))>0$ such that (\ref{beforeduality}) holds.  This concludes the proof of the lemma, and with it Theorem \ref{superthm}.
\end{proof}

\subsection{On the finiteness of (\ref{vdefn})}\label{finiteappendix}

In this subsection we discuss the finiteness of the construction (\ref{vdefn}).  Recall the definition of $\mathbf{V}$ from (\ref{localpot}).  In particular, we prove that:

\begin{lem}\label{lemfinitesuff}
Suppose that, for some $x_0\in \mathbf{R}^n$ and $R>0$:
\begin{equation}\label{tailestsuf}
\int_R^{\infty}\Bigl(\frac{e^{\beta\mathbf{W}_{1,p}^r(d\sigma)(x_0)}}{r^{n-p}} \int_{B(x_0, r)} e^{\beta \mathbf{V}_{B(x,r)}}d\omega(z)\Bigl)^{1/(p-1)} \frac{dr}{r}=C_{\text{tail}}<\infty.
\end{equation}
There exists a constant $C=C(n,p)>0$ so that if $C(\sigma)<C$, then the function $v$ defined in (\ref{vdefn}) is finite almost everywhere.
\end{lem}

This lemma proves the assertion made in Remark \ref{tailrem}.

\begin{proof}
It suffices to prove that, one can choose $C(\sigma)>0$ so that:
$$\int_{B(x_0,R)} v^{\min(1,p-1)}(x) dx<\infty.
$$
To this end, first note that:
\begin{equation}\begin{split}\nonumber\int_{B(x_0, R)} \int_R^{\infty} & \Bigl(\frac{e^{\beta\mathbf{W}_{1,p}^r(d\sigma)(x)}}{r^{n-p}} \int_{B(x, r)} e^{\beta \mathbf{V}_{B(x,r)}}d\omega(z)\Bigl)^{1/(p-1)} \frac{dr}{r}dx\\
& \leq \int_{B(x_0, R)} \int_R^{\infty} \Bigl(\frac{e^{\beta\mathbf{W}_{1,p}^r(d\sigma)(x_0)}}{r^{n-p}} \int_{B(x_0, 2r)} e^{\beta \mathbf{V}_{B(x_0,2r)}}d\omega(z)\Bigl)^{1/(p-1)} \frac{dr}{r}\\
&\qquad \cdot \exp\Bigl(\beta\Bigl[\mathbf{W}_{1,p}^r(d\sigma)(x) -\mathbf{W}_{1,p}^r(d\sigma)(x_0)\Bigl]\Bigl)dx.
\end{split}\end{equation}
Let us now fix $x\in B(x_0,R)$.   Then for any $r>R$, an application of Lemma \ref{lemtailest} yields:
\begin{equation}\begin{split}|\mathbf{W}_{1,p}^r(d\sigma)(x)& -\mathbf{W}_{1,p}^r(d\sigma)(x_0)| \\
& \leq C(n, \beta, C(\sigma)) + \Big[\mathbf{W}_{1,p}^R(\sigma)(x_0) + \mathbf{W}_{1,p}^R(\sigma)(x) \Bigl].
\end{split}\end{equation}
With (\ref{tailestsuf}) in mind, we therefore estimate:
\begin{equation}\begin{split}\nonumber\int_{B(x_0, R)} \int_R^{\infty} & \Bigl(\frac{e^{\beta\mathbf{W}_{1,p}^r(d\sigma)(x)}}{r^{n-p}} \int_{B(x, r)} e^{\beta \mathbf{V}_{B(x,r)}}d\omega(z)\Bigl)^{1/(p-1)} \frac{dr}{r}dx\\
& \lesssim C_{\text{tail}} \int_{B(x_0, R)} \exp[\beta\mathbf{W}_{1,p}^R(d\sigma)(x)] dx \lesssim C_{\text{tail}},
\end{split}\end{equation}
in the last line, we have applied the exponential integrability result Theorem \ref{morremb} (which is valid provided $C(\sigma)<C(n,p)$ for some positive constant $C>0$).

To handle the remaining part of the integral of $v$, let us first suppose $p\geq 2$.   It remains to show that:
\begin{equation}\label{innerfinitepgeq2}\int_{B(x_0, R)}  \int_0^{R} \Bigl(\frac{e^{\beta\mathbf{W}_{1,p}^r(d\sigma)(x)}}{r^{n-p}} \int_{B(x, r)} e^{\beta \mathbf{V}_{B(x,r)}}d\omega(z)\Bigl)^{1/(p-1)} \frac{dr}{r}dx \lesssim C_{\text{tail}}.
\end{equation}
First, by Fubini's theorem and H\"{o}lder's inequality, the left hand side of display (\ref{innerfinitepgeq2}) is less than a constant multiple (depending on $R$) of:
\begin{equation}\label{pgeq2innerR}
\int_0^{R} \Bigl(\frac{1}{r^{n-p}}\int_{B(x_0, R)} e^{\beta\mathbf{W}_{1,p}^r(d\sigma)(x)}\int_{B(x, r)} e^{\beta \mathbf{V}_{B(x,r)}}d\omega(z)dx\Bigl)^{1/(p-1)} \frac{dr}{r}.
\end{equation}
Let us now examine the integrand in (\ref{pgeq2innerR}).  Applying Fubini's theorem, we estimate:
\begin{equation}\begin{split}\int_{B(x_0, R)} & e^{\beta\mathbf{W}_{1,p}^r(d\sigma)(x)}\int_{B(x, r)} e^{\beta \mathbf{V}_{B(x,r)}}d\omega(z)dx\\
& \leq \int_{B(x_0, 2R)}e^{\beta \mathbf{V}_{B(x_0,2R)}}\int_{B(z,r)}e^{\beta\mathbf{W}_{1,p}^r(d\sigma)(x)}dxd\omega(z).
\end{split}\end{equation}
On the other hand, on account of Theorem \ref{morremb},
$$\int_{B(z,r)}e^{\beta\mathbf{W}_{1,p}^r(d\sigma)(x)}dx\leq \int_{B(z,4r)} e^{\beta\mathbf{W}_{1,p}(\chi_{B(z,4r)}d\sigma)(x)}dx\lesssim r^n.
$$
Applying these two observations into (\ref{pgeq2innerR}), we deduce that the left hand side of (\ref{innerfinitepgeq2}) is less than a constant multiple (depending on $R$, $\beta$, $n$, $p$ and $C(\sigma)$) of:
\begin{equation}\begin{split}\nonumber\int_0^{R} \Bigl(\frac{1}{r^{n-p}}r^n  & \int_{B(x_0, 2R)}e^{\beta \mathbf{V}_{B(x_0,2R)}}d\omega(z)\Bigl)^{1/(p-1)}\frac{dr}{r} \\
&\lesssim R^{p/(p-1)} \Bigl(\int_{B(x_0, 2R)}e^{\beta \mathbf{V}_{B(x,2R)}}d\omega(z)\Bigl)^{1/(p-1)}.
\end{split}\end{equation}
This last display is finite as a result of (\ref{tailestsuf}), as required.

The case when $1<p<2$ is similar so we will brief.  We instead consider:
\begin{equation}\label{pleq2finiteinner}\int_{B(x_0, R)} \Bigl\{\int_R^{\infty} \Bigl(\frac{e^{\beta\mathbf{W}_{1,p}^r(d\sigma)(x)}}{r^{n-p}} \int_{B(x, r)} e^{\beta \mathbf{V}_{B(x,r)}}d\omega(z)\Bigl)^{1/(p-1)} \frac{dr}{r}\Bigl\}^{p-1}dx.
\end{equation}
Since $1<p<2$, the previous display is less than:
$$\int_{B(x_0, R)} \int_R^{\infty} \frac{e^{\beta\mathbf{W}_{1,p}^r(d\sigma)(x)}}{r^{n-p}} \int_{B(x, r)} e^{\beta \mathbf{V}_{B(x,r)}}d\omega(z) \frac{dr}{r}dx.
$$
From this point, we use Fubini and Theorem \ref{morremb} as in the case $p\geq2$ to deduce that display (\ref{pleq2finiteinner}) is finite.  From these estimates the lemma follows.
\end{proof}

\section{Lower bounds: the proofs of (\ref{pgeq2entireupbd}), (\ref{pleq2entireup}), (\ref{pleq2lowintrobd}) and (\ref{pgeqlowbdst}).}\label{lowbdsec}  In this section we derive lower bounds for solutions of (\ref{quas2w}).  For a ball $B(x_0, 5R)\subset \Omega$, will be concerned with positive solutions $u$ of:
\begin{equation}\label{localsoln}-\Delta_p u = \sigma u^{p-1}+\omega  \text{ in } B(x_0, 5R), \text{ in the }p\text{-superharmonic sense}.
\end{equation}  
In particular we will prove the following proposition:
\begin{prop} \label{lowbdprop} Let $\Omega$ be an open set.  Suppose $B(x_0, 5R)\subset \Omega$, and suppose that $u$ is a positive solution of (\ref{localsoln}). Then there is a constant $c=c(n,p)>0$, such that:

(i) if $1<p\leq 2$, then
\begin{equation}\begin{split}\label{mainlowestpleq2}
u(x_0) \geq c\int_0^R\Bigl(\frac{e^{c\mathbf{W}_{1,p}^r(\sigma)(x_0)}}{r^{n-p}}\int_{B(x_0, r)} e^{c\mathbf{W}^{r}_{1,p}(\sigma)(z)}d\omega(z)\Bigl)^{1/(p-1)}\frac{dr}{r},
\end{split}\end{equation}

(ii) if $2<p<n$, then
\begin{equation}\begin{split}\label{mainlowestpggeq2}
u(x_0) \geq c\int_0^R\Bigl(\frac{e^{c\mathbf{W}_{1,p}^r(\sigma)(x_0)}}{r^{n-p}}\int_{B(x_0, r)} e^{c\mathbf{I}^{r}_{1,p}(\sigma)(z)}d\omega(z)\Bigl)^{1/(p-1)}\frac{dr}{r}.
\end{split}\end{equation}
\end{prop}
Note that the lower bounds for solutions of (\ref{quas2w}) appearing in Theorems \ref{introthm}, \ref{mainthmleq2entire}, along with Theorem \ref{bddlowbd}, follow from this proposition.  We begin with a some observations regarding the localisation of integral operators:

\subsection{Localisation of operators}\label{localise}  
Let us fix $x_0 \in \Omega$, and let $R>0$ be such that $B(x_0, 5R) \subset \Omega$.  Suppose that $u$ is a positive solution of:
\begin{equation}\label{localsuper}
-\Delta_p u \geq \sigma u^{p-1} \text{ in } B(x_0, 5R), \text{ in the }p\text{-superharmonic sense}.
\end{equation} Then by the potential estimate, Theorem \ref{thmpotest}, it follows that for all $x\in B(x_0, R)$:
\begin{equation}\label{localwolffest1}
u(x) \geq C \int_0^{2R} \Bigl(\frac{1}{r^{n-p}} \int_{B(x,r)} u^{p-1}(z)d\sigma(z)\Bigl)^{1/(p-1)}\frac{dr}{r}.
\end{equation}
First, let us restrict the integration $\sigma$ to $B(x_0, R)$, then clearly:
\begin{equation}\label{localwolffest2}
u(x) \geq C \int_0^{2R} \Bigl(\frac{1}{r^{n-p}} \int_{B(x,r)\cap B(x_0, R)} u^{p-1}(z)d\sigma(z)\Bigl)^{1/(p-1)}\frac{dr}{r}.
\end{equation}
Now, still under the assumption that $x\in B(x_0, R)$, we note that if $r\geq2R$, then $B(x,r) \cap B(x_0, R) = B(x_0, R)$.  Hence we may extend the integration in (\ref{localwolffest2}), i.e.:
\begin{equation}\nonumber
u(x) \geq C \int_0^{\infty} \Bigl(\frac{1}{r^{n-p}} \int_{B(x,r)\cap B(x_0, R)} u^{p-1}(z)d\sigma(z)\Bigl)^{1/(p-1)}\frac{dr}{r}.
\end{equation}
Note here the constant $C$ has changed, but is still a positive constant depending on $n$ and $p$.  Let us now define $ \tilde\sigma$ by:
\begin{equation}\label{tildesigma}d\tilde\sigma= \chi_{B(x_0, R)} d\sigma,\end{equation}
and a nonlinear integral operator $\mathcal{N}$ by:
\begin{equation}\label{defnlocalN}
\mathcal{N}(f)(x) = \mathcal{N}_{\tilde\sigma}(f)(x) = \int_0^{\infty} \Bigl(\frac{1}{r^{n-p}} \int_{B(x,r)} f^{p-1}(z)d\tilde\sigma(y)\Bigl)^{1/(p-1)}\frac{dr}{r},
\end{equation}
The iterates of $\mathcal{N}$ are denoted by $\mathcal{N}^j(f) = \mathcal{N}(\mathcal{N}^{j-1}(f))$.

Using the definition of (\ref{tildesigma}), it follows that, for all $x\in B(x_0, R)$:
\begin{equation}\label{localwolffest4}
u(x) \geq C \mathcal{N}(u)(x).
\end{equation}

Now, suppose that $u$ is a positive solution of (\ref{localsoln}).  By the interior potential estimate, Theorem \ref{thmpotest}, it again follows there exists a constant $C=C(n,p)>0$ so that that for all $x\in B(x_0, R)$:
\begin{equation}\label{localomega}\begin{split}
u(x) 
 \geq C\mathcal{N}(u)(x) + C \mathbf{W}_{1,p}\tilde\omega(x),
\end{split}\end{equation}
where:
\begin{equation}\label{omegatilde}
d \tilde\omega = \chi_{B(x_0, R)}d\omega.
\end{equation}
Mimicking the discussion in Section 5.1 of  \cite{JVFund}, we iterate (\ref{localomega}) to obtain:
\begin{lem}\label{lowbdsum}  Suppose $u$ is a positive solution of (\ref{localsoln}). Let $x_0\in \Omega$ so that $B(x_0, 5R)\subset \Omega$.  Then, with $\mathcal{N}$ and $\tilde\omega$ as in (\ref{defnlocalN}) and (\ref{omegatilde}) repectively,  there is a constant $C=C(n,p)>0$ so that:

\indent a).  If $1<p \leq 2$:
\begin{equation}\label{Nlowbdsleq2} u(x) \geq C \sum_{j=0}^{\infty} C^j \mathcal{N}^j(\mathbf{W}_{1,p}\tilde\omega)(x), \;\text{ for all } x\in B(x_0, R) .\end{equation}
\indent b). If $2< p <n$, then for any $q>1$, 
\begin{equation}\label{Nlowbdsgeq2}
u(x) \geq C(q) \sum_{j=0}^{\infty} j^{(q\frac{2-p}{p-1})}C^j \mathcal{N}^j(\mathbf{W}_{1,p}\tilde\omega)(x), \;\text{ for all } x\in B(x_0, R),
\end{equation}
where $C(q) =C(q,n,p)>0$.
\end{lem}

Proposition \ref{lowbdprop} will follow from careful estimation of the sums in Lemma \ref{lowbdsum}.  In order to estimate the sums, we use the array of tools described in Section \ref{necconds}.  First we note the localized measure $\tilde\sigma$ satisfies the \textit{strong} capacity condition (\ref{strongcapcond}), indeed:
\begin{lem}\label{simplesigmalem}
Under the assumption that $u$ is a solution of (\ref{localsuper}), the measure $\tilde\sigma$ satisfies the strong capacity condition (\ref{strongcapcond}).  More precisely, there is a constant $C=C(n,p)>0$ so that:
$$\tilde\sigma(E)\leq C\text{cap}_p(E)\text{ for all compact sets }E\subset\mathbf{R}^n.
$$
\end{lem}
\begin{proof} Let $E$ be a compact set.  From Lemma 4.3 of \cite{JVFund} (see (\ref{capintro}) above), since $u$ is a solution of (\ref{localsuper}), we have the following estimate:
$$\tilde\sigma(E) = \sigma(E\cap B(x_0,R))\leq \text{cap}_p(E\cap B(x_0,R), B(x_0, 5R)).$$  However, by the separation of $E\cap B(x_0,R)$ and $B(x_0, 5R)$, and since $1<p<n$, we have that (as a consequence of the Sobolev inequality): 
$$\text{cap}_p(E\cap B(x_0,R), B(x_0, 5R))\leq c(n,p)\text{cap}_p(E\cap B(x_0.R)).$$  The lemma follows.
\end{proof}

From Lemma \ref{simplesigmalem} it follows that $\tilde\sigma$ is a Carleson measure, which we state as a lemma.

\begin{lem} \label{carllem} Suppose that $u$ is a positive solution of the inequality $u \geq \mathcal{N}(u)$, with $\mathcal{N}$ as in (\ref{defnlocalN}).  Then the measure $\tilde\sigma$, defined in (\ref{tildesigma}), is a discrete Carleson measure, that is there is a positive constant $C=C(n, p)$ such that for each dyadic cube $P \in \mathcal{Q}$ and  every compact set $E\subset \mathbf{R}^n$,
 \begin{equation} \label{discarl}\sum_{\substack{Q
\subset P\\Q\in \mathcal{Q}}}c_Q \abs{Q\cap E}_{\tilde\sigma}^{p'} \leq C\abs{P \cap E}_{\tilde\sigma}.
\end{equation}
Here $c_Q$ is the sequence defined in (\ref{cQ}).
\end{lem}

\begin{proof}  From Lemma \ref{simplesigmalem}, it follows that $\tilde\sigma$ satisfies (\ref{strongcapcond}).  Therefore, as described before display (\ref{discarla}), it follows that (\ref{discarl}) holds.
\end{proof}

A useful corollary of Lemma \ref{carllem}, together with Lemma \ref{lpdual} from Appendix \ref{duality}, is the following.

\begin{cor}  \label{lowbdduality} Let $1<p\leq2$, and let $\{\lambda_Q\}_{Q\in \mathcal{Q}}$ be any nonnegative sequence indexed over the dyadic cubes.  Suppose that $\tilde\sigma$ satisfies (\ref{discarl}), then there is a constant $C=C(n,p)$ such that:
\begin{equation}
\int_P\Bigl\{\sum_{x\in Q, \, Q\subset P} c_Q \lambda_Q^{1/(p-1)} \Big\}^{p-1}d\tilde\sigma(x) \geq C \sum_{Q\subset P}c_Q\abs{Q}_{\tilde\sigma}^{1/(p-1)}\lambda_Q
\end{equation}
for each dyadic cube $P\in \mathcal{Q}$.
\end{cor}

\begin{proof}  This corollary follows from Lemma \ref{lpdual} in Appendix \ref{duality}.  Letting $s=1/(p-1)$ and relabelling, the lemma boils down to the observation that that $\mu_Q = c_Q^{2-p}|Q|_{\tilde\sigma}^{\frac{2-p}{p-1}}$ is admissible for (\ref{appendspace}) in Lemma \ref{lpdual}.  But this is precisely the statement (\ref{discarl}).
\end{proof}

Having completed our discussion on localisation, we turn to proving Proposition \ref{lowbdprop}.  First up is a lemma regarding the estimation of the `tails' of sums appearing in Lemma \ref{lowbdsum}.

 \subsection{A tail estimate for $1<p<n$}  The following Lemma can be proved by mimicking the proof of Lemma 5.4 in \cite{JVFund}.
\begin{lem}\label{lowerouterbuildup} Let $1<p<n$.  There is a constant $c=c(n,p)>0$ so that for each $m\geq 1$, and every $x\in B(x_0, R)$:
\begin{equation}\begin{split}
\mathcal{N}^m (\mathbf{W}_{1,p})(\tilde\omega)(x)  \geq & \frac{c^m}{m!} \int_0^{\infty} \Bigl(\int_0^r\Bigl( \frac{\tilde\sigma (B(x, t/2))}{t^{n-p}}\Bigl)^{1/(p-1)}\frac{dt}{t}\Bigl)^m(x) \\
&\cdot \int_0^r\Bigl(\frac{\tilde\omega(B(x,s/2))}{s^{n-p}}\Bigl)^{1/(p-1)}\frac{ds}{s}\frac{dr}{r}.
\end{split}\end{equation}
\end{lem}

\subsection{The proof of Proposition \ref{lowbdprop} in case $p\geq 2$}\label{pgeq2lowsec}  Let us begin with the case when $p\geq2$, and we introduce an auxiliary function $\mathbf{B}^t(\tilde\sigma)$, defined by:
\begin{equation}\label{Blowbd}\mathbf{B}^t(\tilde\sigma)(z) = \sum_{j=0}^{\infty} \bigl(\frac{t}{2^j}\bigl)^{p - n} \abs{B(z, t/2^j)}_{\tilde\sigma}
\end{equation}
\begin{lem}\label{sgeq2lowbdlem}Let $p\geq 2$, then there is a constant $C=C(n,p)>0$ such that for all $m\in \mathds{N}\cup\{0\}$, and any $x\in B(x_0, R)$ and $r>0$:
\begin{equation}\begin{split}\label{sgeqlowbdlemst}
\int_{B(x,2r)}\Bigl\{\int_0^{2r} & \Bigl(\frac{1}{t^{n-p}}\int_{B(y,t)}\bigl(\mathbf{B}^t(\tilde\sigma)\bigl)^m(z)d\tilde\omega(z)\Bigl)^{1/(p-1)}\Bigl\}^{p-1}d\tilde\sigma(y)\\
&  \geq \frac{C}{m+1} \int_{B(x,r)} (\mathbf{B}^r(\tilde\sigma)\bigl)^{m+1}(z)d\tilde\omega(z),
\end{split}\end{equation}
here $\mathbf{B}^t(\tilde\sigma)$ is as in (\ref{Blowbd}).
\end{lem}

\begin{proof}  First note that the left hand side of (\ref{sgeqlowbdlemst}) is greater than a constant multiple of:
\begin{equation}\label{pgeq2lowbdstep1}\int_{B(x, 2r)}\Bigl\{\sum_{j=0}^{\infty}\Bigl(\bigl(\frac{r}{2^j}\bigl)^{p-n}\lambda_{B(y,\frac{r}{2^j})}\Bigl)^{1/(p-1)}\Bigl\}^{p-1}d\tilde\sigma(y),
\end{equation}
where:
$$\lambda_{B(y,\frac{r}{2^j})} = \int_{B(y,\frac{r}{2^j})\cap B(x,r)}\bigl(\mathbf{B}^{\frac{r}{2^j}}(\tilde\sigma)\bigl)^m(z)d\tilde\omega(z).$$
Since $p-1\geq 1$, it follows that display (\ref{pgeq2lowbdstep1}) is greater than:
$$\int_{B(x, 2r)}\sum_{j=0}^{\infty}\bigl(\frac{r}{2^j}\bigl)^{p-n} \int_{B(y,\frac{r}{2^j})\cap B(x,r)}\bigl(\mathbf{B}^{\frac{r}{2^j}}(\tilde\sigma)\bigl)^m(z)d\tilde\omega(z) d\tilde\sigma(y).
$$
By Fubini's theorem, the previous display is in turn equal to:
\begin{equation}\label{pgeq2lowest1}\int_{B(x, r)}\sum_{j=0}^{\infty}\bigl(\frac{r}{2^j}\bigl)^{p-n}\tilde\sigma(\{y\in B(x, 2r)\, : \, y\in B(z,\frac{r}{2^j})\})\bigl(\mathbf{B}^{\frac{r}{2^j}}(\tilde\sigma)\bigl)^m(z)d\tilde\omega(z).
\end{equation}
But, for $z \in B(x,r)$ and $t\leq r$, it follows that $B(z,t)\subset B(x,2r).$
Applying this observation in (\ref{pgeq2lowest1}) yields:
\begin{equation}\label{pgeq2lowest2}\int_{B(x, r)}\sum_{j=0}^{\infty}\bigl(\frac{r}{2^j}\bigl)^{p-n}\tilde\sigma(B(z,\frac{r}{2^j}))\bigl(\mathbf{B}^{\frac{r}{2^j}}(\tilde\sigma)\bigl)^m(z)d\tilde\omega(z).
\end{equation}
Let us now recall the following elementary summation by parts inequality.  Let $\{\lambda_j\}_j$ be a nonnegative sequence, and let $m\geq 1$, then:
\begin{equation}\label{sumparts1}\frac{1}{m}\Bigl(\sum_{j=0}^{\infty}\lambda_j\Bigl)^m\leq \sum_{j=0}^{\infty}\lambda_j \Bigl(\sum_{j=k}^{\infty}\lambda_k\Bigl)^{m-1} .\end{equation}
Applying (\ref{sumparts1}) in (\ref{pgeq2lowest2}) 
the lemma follows.
\end{proof}

To complete the proof of Proposition \ref{lowbdprop} when $p\geq 2$, recall the formal Neumann series expansions in Lemma \ref{lowbdsum}.  Applying Lemma \ref{lowerouterbuildup} into Lemma \ref{lowbdsum} results in a constant $c>0$ such that:
\begin{equation}\label{pgeq2lowbd1}
u(x_0) \geq c\int_0^{\infty} \exp\Bigl( c \int_0^r\Bigl( \frac{\tilde\sigma (B(x_0, t/2))}{t^{n-p}}\Bigl)^{1/(p-1)}\frac{dt}{t}\Bigl)\Bigl(\frac{\omega(B(x_0,r))}{r^{n-p}}\Bigl)^{1/(p-1)}\frac{dr}{r}.
\end{equation}
On the other hand, from Lemma \ref{sgeq2lowbdlem}  it follows that there exists a constant $C=C(n,p)>0$ with:
$$\mathcal{N}^m(\mathbf{W}_{1,p}(\tilde\omega)(x_0) \geq C^m\int_0^{\infty}\Bigl(\frac{1}{r^{n-p}}\int_{B(x_0,r)}(\mathbf{B}^r(\tilde\sigma)z))^m d\tilde\omega(z)\Bigl)^{1/(p-1)}\frac{dr}{r},
$$
hence Lemma \ref{lowbdsum} yields a constant $c>0$ such that:
\begin{equation}\begin{split}\label{pgeq2lowbd2}
u(x_0) \geq c\int_0^{\infty}\Bigl(\frac{1}{r^{n-p}}\int_{B(x_0,r)}e^{c\mathbf{B}^r(\tilde\sigma)(z)}d\tilde\omega(z)\Bigl)^{1/(p-1)}\frac{dr}{r}
\end{split}\end{equation}
Averaging (\ref{pgeq2lowbd1}) and (\ref{pgeq2lowbd2}) with the inequality of arithmetic and geometric means, we assert that there is a positive constant $c>0$ such that:
\begin{equation}\begin{split}\label{pgeq2lowbd3}
u(x_0) \geq c\int_0^{\infty}& \exp\Bigl( c \int_0^r\Bigl( \frac{\tilde\sigma (B(x_0, t/2))}{t^{n-p}}\Bigl)^{1/(p-1)}\frac{dt}{t}\Bigl)\Bigl(\frac{1}{r^{n-p}}\\
& \cdot \int_{B(x_0,r)}e^{c\mathbf{B}^r(\tilde\sigma)(z)}d\tilde\omega(z)\Bigl)^{1/(p-1)}\frac{dr}{r}.
\end{split}\end{equation}
But, using (\ref{ballest}), it is easy to see that, for any $z\in B(x_0,R)$:
$$\mathbf{B}^r(\tilde\sigma)(z) \gtrsim \mathbf{I}_p^r(d\tilde\sigma)(z) - 1.
$$
Hence (\ref{mainlowestpggeq2}) follows from (\ref{pgeq2lowbd3}) and recalling the definition of $\tilde\sigma$.
\subsection{The proof of Proposition \ref{lowbdprop} when $1<p<2$} We now move onto the case when $1< p < 2$.  Let $j_0$ be as in Lemma \ref{dyshift}, then we have the following lemma:
\begin{lem}\label{sleq2lowbdlemdy}  Let $1<p<2$, then there is a constant $C=C(n,p)>0$ such that for all $m\in \mathds{N}\cup\{0\}$, and $x\in B(x_0, R)$:
\begin{equation}\begin{split}\label{sleqlowbdlemdyst}
\int_{B(x,r)}&\int_0^{\infty}\Bigl( \frac{1}{s^{n-p}}\int_{B(y, s)}\Bigl(\sum_{\substack{z\in Q\\ Q\in \mathcal{Q}}}c_Q \abs{Q\cap B(y,s)}_{\tilde\sigma}^{1/(p-1)}\Bigl)^m(z)\\&\cdot d\tilde\omega(z)\Bigl)^{1/(p-1)}\frac{ds}{s}\Bigl\}^{p-1}d\tilde\sigma(y)\\
&  \geq \frac{C}{m+1} \int_{B(x,r)} \Bigl(\sum_{\substack{z\in Q\\Q\in \mathcal{Q}}}c_Q \abs{Q\cap B(x,r)}_{\tilde\sigma}^{1/(p-1)}\Bigl)^{m+1}(z)d\tilde\omega(z)
\end{split}\end{equation}
\end{lem}
\begin{proof}  First note that, for any ball $B(y, 2^j)$, there is a dyadic cube $Q$ such that $y\in Q$, $Q\subset B(y, 2^j)$, and $\ell(Q) = 2^{j-j_0}$.  Thus, the left hand side of (\ref{sleqlowbdlemdyst}) is greater than a constant multiple, depending on $n$, $p$ and $j_0$, of:
\begin{equation}\begin{split}\label{lowbdbeforedual}\int_{B(x,r)}\Bigl\{ & \sum_{\substack{ y\in P \\ P\in \mathcal{Q}}} c_P \Bigl( \int_{P\cap B(x,r)}\Bigl(\sum_{\substack{z\in Q\\\ Q\subset P}}c_Q \abs{Q\cap P}_{\tilde\sigma}^{1/(p-1)} \Bigl)^m(z) d\tilde\omega(z)\Bigl)^{1/(p-1)}\Bigl\}^{p-1}d\tilde\sigma(y)
\end{split}\end{equation}
By appealing to duality, in the form of Corollary \ref{lowbdduality}, we see that (\ref{lowbdbeforedual}) is greater than a constant multiple of:
$$ \sum_{P\in \mathcal{Q} } c_P \abs{P\cap B(x, r)}_{\tilde\sigma}^{1/(p-1)}\int_{P\cap B(x,r)}\Bigl(\sum_{\substack{z\in Q\\\ Q\subset P}}c_Q \abs{Q\cap B(x,r)}_{\tilde\sigma}^{1/(p-1)} \Bigl)^m(z) d\tilde\omega(z)
$$
An application of Fubini's theorem, followed by the summation by parts inequality (\ref{sumparts1}) proves the lemma.
\end{proof}
Let us now complete the proof of Proposition \ref{lowbdprop}.  Note that Lemmas \ref{lowerouterbuildup} and \ref{sleq2lowbdlemdy} combine with Lemma \ref{lowbdsum}, as in Section \ref{pgeq2lowsec}, to show that there is a positive constant $c>0$, with:
\begin{equation}\label{dylowbd1}\begin{split}
u(x_0) \geq c\int_0^{\infty}& \exp\Bigl( c \int_0^r\Bigl( \frac{\tilde\sigma (B(x_0, t/2))}{t^{n-p}}\Bigl)^{1/(p-1)}\frac{dt}{t}\Bigl)\Bigl(\frac{1}{r^{n-p}}\\
& \cdot \int_{B(x_0,r)}\exp\Bigl(c\sum_{\substack{z\in Q\\ Q\in \mathcal{Q}}}c_Q \abs{Q\cap B(x_0,r)}_{\tilde\sigma}^{1/(p-1)}\Bigl)d\tilde\omega(z)\Bigl)^{1/(p-1)}\frac{dr}{r}.
\end{split}\end{equation}
It remains to use a shifting argument to recover (\ref{mainlowestpleq2}).  
The bound (\ref{dylowbd1}) continues to hold if we shift the dyadic lattice for any $t\in \mathbf{R}^n$.  Averaging over all shifts $t\in B(0, 2^{k+j_0})$, it follows:
\begin{equation}\begin{split}\nonumber
u(x)\geq &\frac{C}{2^{n(k+j_0)}} \int_{B(0, 2^{k+j_0})}  \int_0^{2^{k}} \exp\Bigl( c \int_0^r\Bigl( \frac{\tilde\sigma (B(x_0, t/2))}{t^{n-p}}\Bigl)^{1/(p-1)}\frac{dt}{t}\Bigl)\Bigl(\frac{1}{r^{n-p}}\\
&\cdot \int_{B(x_0,r)}\exp\Bigl( c\sum_{\substack{z\in Q_t\\ Q_t \in \mathcal{Q}_t }}c_{Q_t} \abs{Q_t\cap B(x_0,r)}_{\tilde\sigma}^{1/(p-1)}\Bigl)d\tilde\omega(z)\Bigl)^{1/(p-1)}\frac{dr}{r}dt.
\end{split}\end{equation}
By Jensen's inequality, the right hand side of the previous display is greater than:
\begin{equation}\begin{split}\label{lowshiftstep2}
\int_0^{2^{k}}& \exp\Bigl( c \int_0^r\Bigl( \frac{\tilde\sigma (B(x_0, t/2))}{t^{n-p}}\Bigl)^{1/(p-1)}\frac{dt}{t}\Bigl)\Bigl(\frac{1}{r^{n-p}}\int_{B(x_0,r)}\exp\Bigl(c\frac{1}{2^{n(k+j_0)}}\\
&\cdot \int_{B(0, 2^{k+j_0})} \sum_{\substack{z\in Q_t\\ Q_t \in \mathcal{Q}_t}}c_{Q_t} \abs{Q_t\cap B(x_0,r)}_{\tilde\sigma}^{1/(p-1)}dt\Bigl)d\tilde\omega(z)\Bigl)^{1/(p-1)}\frac{dr}{r}.
\end{split}\end{equation}
Next, using Lemma \ref{dyshift} in (\ref{lowshiftstep2}), it follows that there exists a constant $c=c(n,p)$ so that:
\begin{equation}\begin{split}
u(x)\geq c\int_0^{2^{k}}& \exp\Bigl( c \int_0^r\Bigl( \frac{\tilde\sigma (B(x_0, t/2))}{t^{n-p}}\Bigl)^{1/(p-1)}\frac{dt}{t}\Bigl)\Bigl(\frac{1}{r^{n-p}}\int_{B(x_0,r)}\\
&\cdot \exp\Bigl(c\int_0^{2^k}\frac{\tilde\sigma(B(z,r))}{r^{n-p}}\frac{dr}{r}\Bigl)d\tilde\omega(z)\Bigl)^{1/(p-1)}\frac{dr}{r},
\end{split}\end{equation}
But now note that for any $z\in B(x_0,R)$, and any $k>\log_2 R + 1$, one can estimate:
$$\int_0^{2^k}\frac{\tilde\sigma(B(z,r))}{r^{n-p}}\frac{dr}{r} \geq c(n,p) \int_0^{\infty}\frac{\tilde\sigma(B(z,r))}{r^{n-p}}\frac{dr}{r}.
$$
Thus, letting $k\rightarrow \infty$, we obtain the required bound (\ref{mainlowestpleq2}).  This completes the proof of Proposition \ref{lowbdprop}.

\section{Existence of solutions to (\ref{quas2w})}\label{existence}
In this section, we conclude the proofs of our main results.  We will consider the case when $\Omega$ is a bounded domain, i.e. Theorem \ref{bddupbd}.   The case $\Omega= \mathbf{R}^n$ (Theorems \ref{introthm} and \ref{mainthmleq2entire}) is similar, using the weak continuity of the $p$-Laplacian operator from \cite{TW1}.  Indeed, see \cite{JVFund}, Section 7 for the argument in the case $\omega = \delta_{x_0}$. There is no problem in generalizing the argument found there to the more general measure, and so we will omit the details here.
\begin{proof}[Proof of Theorem \ref{bddupbd}]
Let $\omega$ be a finite nonnegative measure in $\Omega$, and suppose that, for some $c>0$, the right hand side of either (\ref{pleq2upintrobd}) or (\ref{pgequpbdst}) is finite for some $x\in \Omega$.  By Lemma \ref{lemfinitesuff}, we have that the function $v$ defined in (\ref{vdefn}) is finite almost everywhere with $\beta = c$.  We wish to apply Theorem \ref{superthm}, to assert the existence of a positive constant $C_0 = C_0(n,p,\beta)>0$ such that if $C(\sigma)<C_0$, then there exists $v$ (finite almost everywhere), satisfying:
\begin{equation}\label{existintsup}
v(x) \geq K \mathbf{W}_{1,p}(v^{p-1}d\sigma)(x) + \, K \mathbf{W}_{1,p}(d\omega)(x)
\end{equation}
here $K = \max (1, 2^{\frac{2-p}{p-1}})C_1$, with $C_1$ as in Theorem \ref{thmpotest}.  To this end, note from Theorem \ref{superthm} that there exists such a constant $\tilde C_0>0$ so that provided $C(\sigma)<\tilde C_0$ there exists $C_2>0$ so that:
$$\mathbf{W}_{1,p}(v^{p-1}d\sigma)(x) + \, C_2 \mathbf{W}_{1,p}(d\omega)(x) \leq C_2v(x).
$$
Letting $\displaystyle C_0 = \frac{\tilde C_0}{C_2 K}$, we arrive at (\ref{existintsup}).

We are now in a position to begin the iterative argument.  Let us first define $u_0$ to be a renormalized solution of:
\begin{equation}\nonumber
\begin{cases}
-\Delta_p u_0 = \omega\;  \text{ in } \Omega,\\
\; u_0 = 0\; \text{ on } \partial\Omega.
\end{cases}\end{equation}
Then, note that by Theorem \ref{thmpotest} and (\ref{existintsup}), it follows that that $u_0\leq v$.

We now inductively produce a sequence $\{u_j\}_{j\geq 1}$ such that:
\begin{itemize}
\item $u_j \in \textrm{L}^{p-1}(\Omega, d\sigma)$,
\item $u_j \geq u_{j-1}$ for all $j\geq 1$, 
\item $u_j \leq v$ for all $j\geq 1$, and
\item each $u_j$ is a renormalized solution of:
\begin{equation}\label{renormiterates}
\begin{cases}
-\Delta_p u_j = \sigma u_{j-1}^{p-1} +\omega \; \text{ in } \Omega,\\
\; u_j = 0 \; \text{ on } \partial\Omega.
\end{cases}
\end{equation}
\end{itemize}
To see this, suppose $u_1, \dots, u_{j-1}$ have been constructed.  Then, let $u_{j}$ be a solution of (\ref{renormiterates}), such that $u_j \geq u_{j-1}$.  The existence of such a function $u_j$ is ensured by Lemma 6.9 of \cite{PV}.  By Theorem \ref{thmpotest} and (\ref{existintsup}), it follows that:
\begin{equation}\begin{split}u_j(x) & \leq K \mathbf{W}_{1,p}(u_{j-1}^{p-1}d\sigma)(x) + K \mathbf{W}_{1,p}(\omega)(x) \\
& \leq K \mathbf{W}_{1,p}(v^{p-1}d\sigma)(x) + K \mathbf{W}_{1,p}(\omega)(x) \\
& \leq v(x),
\end{split}\end{equation}
hence the claim follows.  Note that $v\in L^{p-1}(\Omega, d\sigma)$ by construction since $\Omega$ is a bounded domain.

Since the sequence $\{u_j\}_{j\geq 0}$ is increasing, there exists a function $u$ such that $u_j$ converges to $u$.  We wish to conclude that $u$ is a renormalized solution of (\ref{quas2w}).  This will follow as in the stability result of \cite{DMMOP99} once as we have proved that: $\nabla T_k(u_j)\rightarrow \nabla T_k(u)$ in $L^p(\Omega)$, for any $k>0$.  Here $T_k(u) = \min(u, k)$, is the truncation operator.  However, since $\{u_j\}_j$ form an increasing sequence this is not difficult to prove.  Indeed, let us fix $k>0$, then it is well known \cite{DMMOP99} that if $v_j = T_k(u_j)$, $v=T_k(u)$, then $v_j, v  \in W^{1,p}_0(\Omega)$, and $v_j \rightarrow v$ weakly in $W^{1,p}(\Omega)$.  In addition, the truncates are supersolutions, i.e.:
\begin{equation}\label{truncsuper}\int_\Omega \abs{\nabla v_j}^{p-2}\nabla v_j\cdot \nabla \psi \geq 0, \text{ whenever } \psi \in W^{1,p}_0(\Omega) \cap L^{\infty}(\Omega), \text{ and } \psi\geq 0. 
\end{equation}
Under these assumptions, it is known that $v_j \rightarrow v$ strongly in $W^{1,p}_0(\Omega)$, see e.g. the proof of Theorem 3.75 from \cite{HKM}.  Thus, we conclude as in \cite{DMMOP99} that $u$ is a renormalized solution of:
\begin{equation}\begin{cases}
-\Delta_p u = \sigma u^{p-1}+\omega\; \text{ in } \Omega,\\
u = 0 \;\text{ on } \partial\Omega,
\end{cases}\end{equation}
and $u\leq v$.  The proof of the existence of solutions, along with the estimates (\ref{pleq2upintrobd}) and (\ref{pgequpbdst}) is complete.  This completes the proof.\end{proof}

\section{On the equations (\ref{1rhsentire}) and (\ref{riccatiintro})}\label{1bvsection}

In this section, we make a short study of the equation (\ref{1rhsentire}) along with the related equation (\ref{riccatiintro}).  Here we will also assert the pointwise bound (\ref{1bventirebd}) for (\ref{1bventire}) stated in the introduction, which plays an important role in our two weight problem.  The results here are of interest in their own right, as we indicated in the introduction.

Our first goal here is to note the equivalence between (\ref{1bventire}) and (\ref{1rhsentire}).  Indeed, suppose $u$ is a solution of (\ref{1bventire}).  Denoting $v=u+1$ it follows:
\begin{equation}\label{equiv1bvsigma}\min(2^{2-p}, 1) \sigma v^{p-1}\leq -\Delta_p v \leq \max(2^{2-p}, 1)\sigma v^{p-1} \text{ in } \mathbf{R}^n,
\end{equation}
and clearly $\inf_{x\in \mathbf{R}^n}u(x) =0.$  Note the same relationship holds in the converse direction too.  As a result, the bound (\ref{1bventirebd}) for solutions of (\ref{1bventire}) will follow from Proposition \ref{1bvexist} below.  Let us in general consider:
\begin{equation}\begin{cases}\label{1bv}
-\Delta_p u = \sigma u^{p-1} \text{ in } \Omega,\\
\; u=1 \text{ on } \partial\Omega.
\end{cases}\end{equation}  So that equation (\ref{1bv}) reads as (\ref{1rhsentire}) if $\Omega = \mathbf{R}^n$.  Recall the capacity condition (\ref{strongcapcond}) and the constant $C(\sigma)$ along with it.
\begin{prop} \label{1bvexist} (i).  Let $\Omega=\mathbf{R}^n$.  Suppose that $\sigma$ satisfies (\ref{strongcapcond}) for all compact sets $E\subset \mathbf{R}^n$.  Then there exists a constant $C_0 = C_0(n,p)>0$, such that if $C(\sigma)<C_0$, then there exists a solution $u\in W^{1,p}_{\text{loc}}(\mathbf{R}^n)$ of (\ref{1bv}) along with a constant $c>0$ such that:
\begin{equation}\begin{split}\label{1rhsentirebd}
\exp\Bigl[\frac{1}{c}\int_0^{\infty}&\Bigl(\frac{\sigma(B(x,r))}{r^{n-p}}\Bigl)^{1/(p-1)}\frac{dr}{r}\Bigl]\\
&\leq u(x) \leq \exp\Bigl[c\int_0^{\infty}\Bigl(\frac{\sigma(B(x,r))}{r^{n-p}}\Bigl)^{1/(p-1)}\frac{dr}{r}\Bigl].
\end{split}\end{equation}
(ii).  Let $\Omega$ be a bounded domain.  Suppose that $\sigma$ satisfies (\ref{strongcapcond}) for all compact sets $E\subset \mathbf{R}^n$.  Then there exists a constant $C_0 = C_0(n,p)>0$, such that if $C(\sigma)<C_0$, then there exists a solution $u\in W^{1,p}(\Omega)$ of (\ref{1bv}), along with a constant $c>0$ such that:
\begin{equation}\begin{split}\label{1rhsbddbd}
\exp\Bigl[\frac{1}{c}\int_0^{d(x)/2}&\Bigl(\frac{\sigma(B(x,r))}{r^{n-p}}\Bigl)^{1/(p-1)}\frac{dr}{r}\Bigl]\\
&\leq u(x) \leq \exp\Bigl[c\int_0^{2\text{diam}(\Omega)}\Bigl(\frac{\sigma(B(x,r))}{r^{n-p}}\Bigl)^{1/(p-1)}\frac{dr}{r}\Bigl]
\end{split}\end{equation}
\end{prop}
Here $d(x)$ is the distance to the boundary on $\Omega$.  
To prove the lower bounds in (\ref{1rhsentirebd}) and (\ref{1rhsbddbd}), it will be convenient to go through the equation:
\begin{equation}\begin{cases}\label{riccati}
-\Delta_p u = (p-1)\abs{\nabla u}^p + \sigma, \text{ in }\Omega,\\
\,u = 0\, \text{ on }\partial\Omega.
\end{cases}\end{equation}
The connection between (\ref{1bv}) and (\ref{riccati}) is the content of the following result:
\begin{prop} \label{logsubgen} Let $\Omega$ be a connected open set.  Let $\mu^s$ be a measure singular with respect to capacity, and suppose $\sigma$ is absolutely continuous with respect to capacity.  Suppose that $u\in L_{\text{loc}}^{p-1}(\Omega,d\sigma)$ is a $p$-superharmonic solution of:
\begin{equation}\begin{cases}\label{singularRHS}
-\Delta_p  u = \sigma u^{p-1} + \mu^s \text{ in }\Omega,\\
\,u>0 \text{ in }\Omega.
\end{cases}\end{equation}  Then $v=\log u \in \textrm{W}^{\, 1,p}_{\text{loc}}(\Omega)$, and $v$ satisfies (\ref{riccati}) in the sense of distributions.
\end{prop}

The proof of Proposition \ref{logsubgen} below makes use of a recent result of Kilepl\"{a}inen, Kuusi and Tuhola-Kujanp\"{a}\"{a}  \cite{KKT11}.  This replaces our original argument based on techniques of Dal Maso and Malusa \cite{DMM97}.  Similar results have recently been proved in \cite{AHBV}.  

\begin{proof}[Proof of Proposition \ref{logsubgen}] Since $u$ is $p$-superharmonic, from Theorem 3.2 of \cite{KKT11} it follows that $u$ is a local renormalized solution.  By definition, this means that for each $\phi\in C^{\infty}_0(\Omega)$ and $h\in W^{1,\infty}(\mathbf{R})$ with derivative $h'$ having compact support, we have:
\begin{equation}\label{localrenorm}
\int_{\Omega} |\nabla u|^{p-2}\nabla u \cdot \nabla (h(u)\phi) dx = \int_{\Omega} h(u)\phi u^{p-1}d\sigma + \int_{\Omega} h(u) \phi d\mu^s.
\end{equation}
Assume $\phi\geq 0$.  For $k>0$, let $h(u) = T_k(u)^{1-p}$, with $T_k(s) = \min(k,s)$.  Then:
\begin{equation}\label{localrenorm2}\begin{split}
\int_{\Omega} \frac{|\nabla u|^{p-1}\nabla u}{T_k(u)^{p-1}} \cdot \nabla \phi dx = & (p-1)\int_{\Omega\cap \{u\leq k\}}\frac{|\nabla u|^p}{u^p}\phi dx + \int_{\Omega} \phi \frac{u^{p-1}}{T_k(u)^{p-1}}d\sigma \\
&+ \int_{\Omega} \frac{\phi}{T_k(u)^{p-1}} d\mu^s.\end{split}
\end{equation}
We now need a few standard properties of $p$-superharmonic functions.  First, since $\mu^s$ is singular with respect to capacity, we have that $\mu^s(\{u<k\}) = 0$ for all $k>0$ (see e.g. Lemma 2.9 of \cite{KKT11}).  Therefore:
$$ \int_{\Omega} \frac{\phi}{T_k(u)^{p-1}} d\mu^s = \frac{1}{k^{p-1}}\int_{\Omega}\phi d\mu^s \rightarrow 0\text{ as }k\rightarrow\infty.
$$
Next, recall that (see e.g. Chapter 7 of \cite{HKM}) that $v = \log(u) \in W^{1,p}_{\text{loc}}(\Omega)$.  Therefore, from Lebesgue's dominated convergence theorem:
$$\int_{\Omega\cap \{u\leq k\}}\frac{|\nabla u|^p}{u^p}\phi dx \rightarrow \int_{\Omega}|\nabla v|^p hdx \text{ as }k\rightarrow \infty,
$$
and, also:
$$\int_{\Omega} \frac{|\nabla u|^{p-1}\nabla u}{T_k(u)^{p-1}} \cdot \nabla \phi dx\rightarrow \int_{\Omega}|\nabla v|^{p-2}\nabla v\cdot \nabla \phi dx \text{ as }k\rightarrow \infty.
$$
From a final application of the monotone convergence on the term involving the measure $\sigma$, we deduce from letting $k\rightarrow \infty$ in (\ref{localrenorm2}) that $v$ is a weak solution of (\ref{riccati}).
\end{proof}

Using Proposition \ref{logsubgen}, we readily conclude the lower bounds in Proposition \ref{1bvexist}.  Indeed, both follow from the following local lemma:

\begin{lem}\label{1bvlowbdprop} Let $\Omega$ be an open set, and suppose that $u$ is a $p$-superharmonic solution of (\ref{1bv}).  Then, there exists a constant $c=c(n,p)>0$:
\begin{equation}\label{1bvlowbd2}
u(x) \geq \exp\Bigl[c\int_0^{d(x)/2} \Bigl(\frac{\sigma(B(x,r))}{r^{n-p}}\Bigl)^{1/(p-1)}\frac{dr}{r} \Bigl],\text{ for all }x\in \Omega.
\end{equation}
\end{lem}

\begin{proof}  By Proposition \ref{logsubgen}, $v = \log u$ solves (\ref{riccati}) in the weak sense.  This clearly implies that $v$ solves $-\Delta_p v \geq \sigma$ in the $p$-superharmonic sense in $\Omega$.  Applying Theorem \ref{potest}, it follows there is a constant $c=c(n,p)>0$ so that:
$$v(x) \geq c\mathbf{W}_{1,p}^{\frac{d(x)}{2}}(d\sigma)(x) \; \text{ for all } x\in \Omega.
$$
Recalling that $v= \log u$, the proposition follows.
\end{proof}

Turning now to the existence of solutions along with the upper bounds in Proposition \ref{1bvexist}, the primary ingredient is the following lemma.  The proof also serves as a prototype of the kind of supersolutions developed in Section \ref{construct}.

\begin{lem} \label{super1bv} Suppose that $\sigma$ is a nonnegative measure satisfying (\ref{strongcapcond}), and define:
\begin{equation}\label{vrhs1}v(x) = \exp\Bigl[\beta\int_0^{\infty} \Bigl(\frac{\sigma(B(x,r))}{r^{n-p}}\Bigl)^{1/(p-1)}\frac{dr}{r} \Bigl]=  \exp\bigl[\beta\mathbf{W}_{1,p}(d\sigma)(x)\bigl]\end{equation}
For any $\beta>0$, there exists $C_0=C_0(n,p)>0$ such that if $C(\sigma)<C_0$, then there exists a constant $C=C(n,p)>0$ such that:
:\begin{equation}\label{super1bvint}
\mathbf{W}_{1,p}(v^{p-1}d\sigma)(x) \leq C (v(x)-1).
\end{equation}
Furthermore, $\inf_{x\in \mathbf{R}^n} v= 1$.
\end{lem}

\begin{proof}[Proof of Lemma \ref{super1bv}]  
Writing out $\mathbf{W}_{1,p}(v^{p-1} d\sigma)(x)$, and applying Lemma \ref{lemtailest}, we derive:
\begin{equation}\begin{split}\label{super1bvstep2}
\mathbf{W}_{1,p}(v^{p-1} & d\sigma)(x)   \lesssim \int_0^{\infty}  \Bigl[\frac{1}{r^{n-p}} \exp\Bigl((p-1)\beta\int_r^{\infty}\Bigl(\frac{\sigma(B(x,t))}{t^{n-p}}\Bigl)^{\frac{1}{p-1}}\frac{dt}{t}\Bigl)\\
& \cdot \int_{B(x,r)}\exp\Bigl((p-1)\beta\int_0^r\Bigl(\frac{\sigma(B(y,t))}{t^{n-p}}\Bigl)^{\frac{1}{p-1}}\frac{dt}{t}\Bigl)d\sigma(y)\Bigl]^{1/(p-1)}\frac{dr}{r}.
\end{split}\end{equation}
On the other hand, as a result of the exponential integrability lemma (Lemma \ref{expest}), there exists a positive constant $C_0 = C_0(n,p)>0$, so that, provided $\beta C(\sigma)<C_0$
\begin{equation}\begin{split}\nonumber
\int_{B(x,r)}& \exp\Bigl((p-1)\beta\int_0^r\Bigl(\frac{\sigma(B(y,t))}{t^{n-p}}\Bigl)^{\frac{1}{p-1}}\frac{dt}{t}\Bigl)d\sigma(y)
\lesssim \sigma(B(x,2r)).
\end{split}
\end{equation}
Note that also, using the estimate on balls (\ref{ballest}), it follows that there exists a positive constant $C=C(n,p, \beta, C(\sigma))$, so that:
\begin{equation}\nonumber \exp\Bigl(\beta\int_r^{2r}\Bigl(\frac{\sigma(B(y,t))}{t^{n-p}}\Bigl)^{\frac{1}{p-1}}\frac{dt}{t}\Bigl) \leq C
\end{equation}
Substituting these two estimates into (\ref{super1bvstep2}): 
\begin{equation}\nonumber\begin{split}
\mathbf{W}_{1,p}(v^{p-1}d\sigma)(x) & \lesssim  \int_0^{\infty} \Bigl(\frac{\sigma(B(x,2r))}{(2r)^{n-p}}\Bigr)^{\frac{1}{p-1}}\exp\Bigl(\beta\int_{2r}^{\infty}\Bigl(\frac{\sigma(B(x,t))}{t^{n-p}}\Bigl)^{\frac{1}{p-1}}\frac{dt}{t}\Bigl)\frac{dr}{r}\\
& \lesssim \int_0^{\infty} \Bigl(\frac{\sigma(B(x,r))}{r^{n-p}}\Bigr)^{\frac{1}{p-1}}\exp\Bigl(\beta \int_{r}^{\infty}\Bigl(\frac{\sigma(B(x,t))}{t^{n-p}}\Bigl)^{\frac{1}{p-1}}\frac{dt}{t}\Bigl)\frac{dr}{r}\\
&\lesssim (v(x)-1),
\end{split}\end{equation}
the last inequality in the sequence follows from integration by parts.  The statement that $\inf_{x\in \mathbf{R}^n} v(x) =1$ follows from:
$$\inf_{x\in \mathbf{R}^n} \mathbf{W}_{1,p}(d\sigma)(x) = 0,
$$
the latter assertion may be verified in a similar manner to the argument around display (6.20) of \cite{JVFund}.  This concludes the proof of the lemma.
\end{proof}
Note that, for any measure $\sigma$ supported on $\Omega$, it follows:
$$\mathbf{W}_{1,p}(d\sigma)(x) \lesssim \mathbf{W}_{1,p}^{2\text{diam}(\Omega)}(d\sigma)(x), \text{ for all } x\in \Omega,
$$
where the implied constant in the inequality depends only on $n$ and $p$.  Let us now complete the proof of Proposition \ref{1bvexist}:
\begin{proof}[Proof of the existence part of Proposition \ref{1bvexist}]  Part (i).  Let us denote:
\begin{equation}\label{super1rhsentire}v(x) = \exp \bigl(\mathbf{W}_{1,p}(d\sigma)(x)\bigl),
\end{equation}
and note that by virtue of Lemma \ref{super1bv}, we may choose $C(\sigma)$ small enough so that:
\begin{equation}\label{entire1bvsuper}C_1\mathbf{W}_{1,p}(v^{p-1})(x) + 1 \leq v(x)
\end{equation}
with $C_1$ as in Theorem \ref{thmpotest}.   We claim that one can choose $C(\sigma)$, depending on only $n$ and $p$, so that both:
\begin{equation}\label{localvint}
v\in L^{p}_{\text{loc}}(\mathbf{R}^n), \text{ and } v\in L^{p}_{\text{loc}}(\mathbf{R}^n, d\sigma).
\end{equation}
These two properties follow from applications of Theorem \ref{morremb} and Lemma \ref{corexpest} respectively.
Let us construct a sequence $(u_j)_j$ so that $u_0=1$ and:
\begin{equation}\begin{cases}\label{entire1bvapprox}
-\Delta_p u_j = \sigma u_{j-1}^{p-1}\text{ in }\mathbf{R}^n,\\
\,u_j \leq v,\, u_j \in W^{1,p}_{\text{loc}}(\mathbf{R}^n),\\
\,\inf_{\mathbf{R}^n}u_j =1.
\end{cases}\end{equation}
Suppose that $u_0,\dots u_{j-1}$ have been constructed.  Then let $B_{\ell} = B(0, 2^{\ell})$ and denote by $u_j^{\ell}$ the unique solution of:
$$-\Delta_p u_j^{\ell} = \sigma u^{p-1}_{j-1} \text{ in } B_{\ell}, \, u_j^{\ell}-1\in W^{1,p}_0(B_{\ell}).
$$ 
the existence of $u_j^{\ell}$ follows from monotone operator theory (see Proposition 5.1 in Chapter 2 of \cite{Sho97}), since $\sigma u^{p-1}_{j-1}\in W^{-1,p'}(B_{\ell})$ by hypothesis and Theorem \ref{multcond}.   Note that $u_j^{\ell}\leq v$, as follows from Theorem \ref{thmpotest} together with the assumption $u_{j-1}\leq v$ and display (\ref{entire1bvsuper}).

The sequence $(u_j^{\ell})_{\ell}$ form an increasing sequence by the classical comparison principle, and therefore using Theorem 1.17 of \cite{KM92} in combination with the weak continuity of the $p$-Laplacian (see \cite{TW1}), we can find a $p$-superharmonic function $u_j$ so that $u_j^{\ell}\rightarrow u_j$ almost everywhere and :
$$-\Delta_p u_j =\sigma u_{j-1}^{p-1} \text{ in }\mathcal{D}'(\mathbf{R}^n).
$$
Furthermore $1\leq u_j\leq v$, so $\inf_{\mathbf{R}^n}u_j = 1$.  To see that $u_j \in W^{1,p}_{\text{loc}}(\mathbf{R}^n)$, note that if $\phi\in C^{\infty}_0(B_{\ell})$:
$$\int_{\Omega} |\nabla u^{\ell}_j|^p \phi^p  \leq \int_{\Omega} u_{j-1}^{p-1}u_j^{\ell}\phi^p d\sigma + p\int_{\Omega} |\nabla u_{j}^{\ell}|^{p-1}u_j^{\ell} |\nabla \phi|\phi^{p-1} dx
$$
Using Young's inequality on the right hand side along with blunt estimates we obtain:
\begin{equation}\label{localsobest}\int_{\Omega} |\nabla u^{\ell}_j|^p \phi^p  \lesssim \int_{\Omega} v^{p}|\nabla\phi|^p dx + \int_{\Omega} v^{p}\phi^p d\sigma <\infty.
\end{equation}
The finiteness follows from (\ref{localvint}).  Using weak compactness along with the a.e. convergence it follows that $u_j\in W^{1,p}_{\text{loc}}(\mathbf{R}^n)$. The sequence (\ref{entire1bvapprox}) has been constructed.
Appealing once again to Theorem 1.17 of \cite{KM92} with the weak continuity of quasi-linear operators, we deduce the existence of a solution  $u$ of:
$$-\Delta_p u =\sigma u^{p-1} \text{ in }\mathbf{R}^n.
$$
Furthermore $1\leq u\leq v$ and so $\inf_{\mathbf{R}^n}u = 1$.  It remains to show that $u\in W^{1,p}_{\text{loc}}(\mathbf{R}^n)$, by this follows immediately from (\ref{localsobest}) and weak compactness.  Recalling the definition of $v$ from (\ref{super1rhsentire}), the proposition follows.

The proof of part (ii) is easier.  Indeed, first note that under the present assumptions in $\sigma$, it follows from Theorem \ref{multcond} and H\"{o}lder's inequality that $\sigma \in W^{-1, p'}(\Omega)$.  Let $u_0 = 1$.  Appealing to monotone operator theory (see e.g. \cite{Sho97}), we inductivly find a sequence $\{u_j\}$ such that:
\begin{equation}\label{sobapprox}\begin{cases}
-\Delta_p u_j = \sigma u_{j-1}^{p-1},\\
u_j - 1 \in W^{1,p}_0(\Omega).
\end{cases}
\end{equation}
Arguing inductively as in part (i), we can further suppose that $u_j\leq v$ for all $j$.  In addition, we may suppose $v\in L^{p}(\Omega)\cap L^p(\Omega, d\sigma)$, again as in part (i) (recall $\Omega$ is bounded).

Furthermore, the sequence $(u_j)_j$ is increasing, and hence testing $u_j-1$ in (\ref{sobapprox}) and using Minkowski's inequality:
\begin{equation}\label{standardsobest1}\Bigl(\int_{\Omega} |\nabla u_j|^p dx\Bigl)^{1/p} \leq \Bigl(\int_{\Omega} (u_j -1)^p d\sigma \Bigl)^{1/p} + \sigma(\Omega)^{1/p} \lesssim \Bigl(\int_{\Omega} v^p d\sigma\Bigl)^{1/p}
\end{equation}
From (\ref{standardsobest1}), we find that there exists $u\in W^{1,p}(\Omega)$, so that $u_j \rightarrow u$ weakly in $W^{1,p}(\Omega)$.  Since $u_j$ is an increasing sequence of supsersolutions, it is standard (see e.g \cite{HKM}) that $u_j\rightarrow u$ strongly in $W^{1,p}(\Omega)$\footnote{This has already been seen in Section \ref{existence}}.  It readily follows that $u$ is a solution of (\ref{1bv}).  Furthermore $u\leq v$.\end{proof}



Since it may be of interest, we collect the results proved for the equation (\ref{ricthm}) into a single result:
\begin{thm}\label{ricthm} Let $\Omega$ be a bounded domain, or $\Omega = \mathbf{R}^n$.\\
\indent a). Suppose that $u\in W^{1,p}_{\text{loc}}(\Omega)$, is a positive weak solution of (\ref{riccati}).  Then $\sigma$ satisfies (\ref{capintro}), with constant $C=1$, and furthermore:
$$u(x) \geq \mathbf{W}_{1,p}^{d(x)/2}(d\sigma)(x), \text{ for all } x\in \Omega.
$$
\indent b).  Conversely, suppose that $\sigma$ satisfies (\ref{strongcapcond}), then there exists $C_0 = C_0(n,p)>0$ so that if $C(\sigma)<C_0$, then there is a positive weak solution of $u\in W^{1,p}_{\text{loc}}(\Omega)$ of (\ref{riccati}), such that:
$$u(x) \leq \mathbf{W}_{1,p}^{2\text{diam}(\Omega)}(d\sigma)(x), \text{ for all } x\in \Omega.
$$
Furthermore, if $\Omega$ is bounded, then $u$ can be chosen to lie in the class $W^{1,p}_0(\Omega)$.
\end{thm}
\begin{proof}
The proof of {\rm (a)} follows immediately from Lemma 4.2 of \cite{JVFund}, and Theorem \ref{thmpotest} above.  On the other hand, part {\rm (b)}  is an immediate corollary of Propositions \ref{1bvexist} and \ref{logsubgen}.
\end{proof}

\section{Examples}\label{examples}

In this section we will discuss conditions on $\omega$ so that our theorems guarantee the existence of a solution to (\ref{quas2w}).

\subsection{$L^q$ data}  Our first example concerns the case when $\omega \in L^q(\Omega)$ for some $q>1$.  Here we will let $\Omega$ be a bounded domain.  We will always extend measures and functions by zero outside of $\Omega$ so they are defined on $\mathbf{R}^n$.
\begin{prop} \label{lqexist} Let $1<p<n$, and let $\Omega$ be a bounded domain.  Furthermore, let $\omega \in L^q(\Omega)$ for $q>1$, and suppose that $\sigma$ satisfies (\ref{strongcapcond}).  Then, there exists $C_0=C(\sigma, q)$ so that if $C(\sigma)<C_0$ then there a positive constant $c>0$, depending on $n$ and $p$, together with a renormalized solution $u$ of (\ref{quas2w}) so that:
\begin{equation}\label{lpbound}\begin{split}
u(x) \leq C(n,p,q,\Omega) \Bigl[\int_0^{\text{diam}(\Omega)}& \Bigl(\frac{1}{r^{n-p}}\int_{B(x,r)} \omega^q(z) dz \Bigl)^{1/(p-1)}\frac{dr}{r}\Bigl]^{1/q}\\
&\cdot \exp \Bigl[\frac{c}{q'}\int_0^{\text{diam}(\Omega)} \Bigl(\frac{\sigma(B(x,r))}{r^{n-p}}\Bigl)^{1/(p-1)}\frac{dr}{r}\Bigl].
\end{split}\end{equation}
Furthermore, for all $r<\displaystyle \frac{nq}{n-p}$, one can choose $C_0$ to in addition depend on $r$ so that $u\in L^r(\Omega)$.
\end{prop}

\begin{proof}  Our first aim is to show that the quantities appearing in (\ref{pgequpbdst}) and (\ref{pleq2upintrobd}) are finite almost everywhere for a choice of $c>0$.  To this end, let us momentarily fix $x\in \Omega$, and denote by $\mathbf{T}$:
\begin{equation}\nonumber\mathbf{T}^r(\sigma)(z) = \begin{cases} \mathbf{I}_p(\chi_{B(x,r)}d\sigma)(z) \text{ if }1<p<2,\\
\mathbf{W}_{1,p}(\chi_{B(x,r)}d\sigma)(z) \text{ if } p\geq 2.
\end{cases}
\end{equation}
This allows us to deal with both cases $1<p<2$ and $p\geq 2$ simultaneously.  Note by H\"{o}lder's inequality:
\begin{equation}\label{lqholder}\int_{B(x,r)} e^{c\mathbf{T}^r(\sigma)(z)}d\omega \leq \Bigl(\int_{B(x,r)} e^{cq'\mathbf{T}^r(\sigma)(z)}dx\Bigl)^{1/q'}\Bigl(\int_{B(x,r)} \omega^q dx\Bigl)^{1/q}.
\end{equation}
Since we are assuming that $\sigma$ satisfies (\ref{strongcapcond}), it follows that the hypothesis of Theorem \ref{morremb} are valid (in the case of Lebesgue measure).  Therefore there exists $c_1=c_1(n,p)>0$ so that provided $cq'<c_1$, there is a constant $C=C(n,p,q)$ so that:
\begin{equation}\label{lqball}\int_{B(x,r)} e^{c\mathbf{T}^r(\sigma)(z)}d\omega \leq C(n,p,q)r^{n/q'}\Bigl(\int_{B(x,r)} \omega^q dx\Bigl)^{1/q}.
\end{equation}
Let us now form the right hand sides appearing in  (\ref{pgequpbdst}) and (\ref{pleq2upintrobd}):
\begin{equation}\label{lqrhs}\int_0^{d_{\Omega}} \Bigl[\frac{e^{c\mathbf{W}^r_{1,p}(\chi_{\Omega})(x)}}{r^{n-p}} \int_{B(x,r)}e^{c\mathbf{T}^r(\sigma)(z)} d\omega\Bigl]^{1/(p-1)}\frac{dr}{r}.
\end{equation}
Substituting (\ref{lqball}) into (\ref{lqrhs}) and appealing to H\"{o}lder's inequality a second time, we derive that (\ref{lqrhs}) is less than a constant multiple of:
\begin{equation}\begin{split}\label{lqupbd1}\Bigl(\int_0^{d_{\Omega}} & \Bigl[\frac{e^{c\mathbf{W}^r_{1,p}(\chi_{\Omega})(x)}}{r^{n-p}} \cdot r^n\Bigl]^{1/(p-1)}\frac{dr}{r}\Bigl)^{1/q'}\\
&\cdot\Bigl(\int_0^{d_{\Omega}} \Bigl[\frac{e^{c\mathbf{W}^r_{1,p}(\chi_{\Omega})(x)}}{r^{n-p}} \int_{B(x,r)}\omega^q dz \Bigl]^{1/(p-1)}\frac{dr}{r}\Bigl)^{1/q}
\end{split}\end{equation}
Since (\ref{lqupbd1}) is bounded by a constant multiple of the right hand side of (\ref{lpbound}), it follows that both (\ref{pgequpbdst}) and (\ref{pleq2upintrobd}) are finite almost everywhere for a choice of $c>0$.  Hence by Theorem
 \ref{bddupbd}, there exists a solution $u$ of (\ref{quas2w}) so that (\ref{lpbound}) holds.

To see the regularity property, note that, for any $\omega \in L^q$, we have:
$$\int_0^{\text{diam}(\Omega)} \Bigl(\frac{1}{r^{n-p}}\int_{B(x,r)} \omega^q(z) dz \Bigl)^{1/(p-1)}\frac{dr}{r}\, \in L^s(\Omega), \text{ for all } s<n/(n-p).
$$
This is a standard regularity property of the nonlinear potential (and remains true if $\omega^q$ is replace by a finite measure).  For any $r<qn/n-p$, choose $s$ so that:
$$\frac{r}{q} < s< \frac{n}{n-p}.
$$
We see by applying H\"{o}lder's inequality in (\ref{lpbound}):
\begin{equation}\label{lpbound2}\begin{split}
\int_{\Omega}u^rdx \leq C &\Bigl[\int_{\Omega} \Bigl(\int_0^{\text{diam}(\Omega)} \Bigl(\frac{1}{r^{n-p}}\int_{B(x,r)} \omega^q(z) dz \Bigl)^{1/(p-1)}\frac{dr}{r}\Bigl)^{s}dx\Bigl]^{r/qs}\\
& \cdot \Bigl[\int_{\Omega} \exp\Bigl(c\frac{sq}{s-q}\int_0^{d_{\Omega}}\Bigl(\frac{\sigma(B(x,r))}{r^{n-p}}\Bigl)^{1/(p-1)}\frac{dr}{r} \Bigl)dx \Bigl]^{\frac{1-sq}{sq}}.
\end{split}\end{equation}
From Theorem \ref{morremb} together with a covering of $\Omega$ by balls, we can pick $c>0$ (or $C(\sigma)>0$) so that:
$$\int_{\Omega} \exp\Bigl(c\frac{sq}{s-q}\int_0^{d_{\Omega}}\Bigl(\frac{\sigma(B(x,r))}{r^{n-p}}\Bigl)^{1/(p-1)}\frac{dr}{r} \Bigl)dx \leq C(\Omega),
$$
the proposition follows.
\end{proof}

\subsection{Weak-$A_{\infty}$ data} The second class we consider is when $\omega$ is a \textit{weak-$A_{\infty}$ weight}.  We recall that nonnegative function $w$ is a weak $A_{\infty}$ weight if there are constants $C_w>0$ and $\theta>0$ so that, for all balls $B$, and measurable subsets $E\subset B$:
\begin{equation}\label{weakainfdef2}
\frac{|E|_w}{|2B|_w} \leq C_w\Big(\frac{|E|}{|Q|}\Bigl)^{\theta}.
\end{equation}
All locally integrable power weights $\omega = |x|^{q}$ for $q>-n$ are included in this class, and these form an important subclass of right hand sides for the equation (\ref{quas2w}).  In this subclass, the pointwise bound presented in our main theorems collapses to a much simpler expression.  The point here is that we continue to \textit{obtain a sharp bilateral pointwise bound} for such $\omega$, in contrast to the $L^q$ case considered above.  

\begin{prop}\label{ainfcoincide}  Let $\Omega = \mathbf{R}^n$ and $1<p<n$.  Suppose that $\sigma$ satisfies (\ref{strongcapcond}), and suppose that $\omega$ is a weak $A_{\infty}$ weight.  There exist constants $c_1, c_2, C_0>0$, depending on $n, p, \theta$ and $C_w$, so that if $C(\sigma)<C_0$, then there exists a $p$-superharmonic solution of (\ref{quas2w}) so that:
\begin{equation}\begin{split}\label{ainfcoincidest}
c_1\int_{0}^{\infty}\Bigl(&\frac{e^{c_1\mathbf{W}^r_{1,p}( d\sigma)(x)}}{r^{n-p}} |B(x,r)|_{\omega}\Bigl)^{1/(p-1)}\frac{dr}{r}\\
&\leq u(x) \leq c_2\int_{0}^{\infty}\Bigl(\frac{e^{c_2\mathbf{W}^r_{1,p}( d\sigma)(x)}}{r^{n-p}} |B(x,r)|_{\omega}\Bigl)^{1/(p-1)}\frac{dr}{r}
\end{split}\end{equation}
\end{prop}

\begin{proof}[Proof of Proposition \ref{ainfcoincide}]  In light of Theorems \ref{introthm} and \ref{mainthmleq2entire}, the lower bound in display (\ref{ainfcoincidest}) is clear, so we only need to show the upper bound.   To this end we need to assert the existence of a constant $c>0$ so that the bounds (\ref{pgeq2entireupbd}) and (\ref{pleq2entireup}) are finite, and are bounded above by the expression appearing in the right hand side (\ref{ainfcoincidest}).  By inspection, this will follow if we show (with $q=\min(1,1/(p-1))$):
\begin{equation}\label{corofmoremb}\int_{B(x,r)}\exp\Bigl[c\int_0^{r}\Bigl(\frac{\sigma(B(z,s)}{s^{n-p}}\Bigl)^q\frac{ds}{s}\Bigl] d\omega(z) \leq C|B(x,4r)|_{\omega}.
\end{equation}
Provided $C(\sigma)$ or $c>0$ are small enough, display (\ref{corofmoremb}) follows directly from Theorem \ref{morremb}.
\end{proof}

\subsection{Morrey space data} Our third example of $\omega$ is when it lies in a suitable Morrey space.  We will again discuss the case when $\Omega$ is a bounded domain.  The precise condition we use is the following:  there exists $\epsilon>0$ and a constant $C(\omega)>0$ so that:
\begin{equation}\label{omemor}\omega(B(x,r))\leq C(\omega)r^{n-p+\epsilon} \text{ for all balls }B(x,r)\subset \mathbf{R}^n.
\end{equation}
Note that this condition is not contained in the $L^q$ space condition, since no higher integrability is assumed.  With this condition, the following holds:
\begin{prop}\label{morexist}  Let $\Omega$ be a bounded domain.  Suppose that $\omega$ satisfies (\ref{omemor}), and that $\sigma$ satisfies (\ref{strongcapcond}).  Then, there exists $C_0=C(n,p,C(\omega))$ so that if $C(\sigma)<C_0$ then there is a positive constant $c>0$, depending on $n, p$ and $C(\omega)$, together with a renormalized solution $u$ of (\ref{quas2w}) so that:
\begin{equation}
u(x) \leq C(C(\omega), \Omega, n, p)\exp \Bigl[c\int_0^{\text{diam}(\Omega)} \Bigl(\frac{\sigma(B(x,r))}{r^{n-p}}\Bigl)^{1/(p-1)}\frac{dr}{r}\Bigl].
\end{equation}
\end{prop}

The proposition follows in precisely a very similar fashion to Proposition \ref{lqexist} above once we have asserted a suitable integrability result.  We therefore will prove the following lemma, and leave the remainder of the argument to the reader:

\begin{lem}  Suppose that $\sigma$ satisfies the weak ball estimate (\ref{ballest}) with constant $C(\sigma)>0$, and that $\omega$ satisfies (\ref{omemor}), then there exists a constant $c>0$, depending on $n$, $p$, $C(\sigma)$ and $C(\omega)$ so that, for any ball $B(x,r)$:
\begin{equation}\label{expinteps}
\int_{B(x,r)} \exp\Bigl(c\int_0^r \frac{\sigma(B(z,r))}{r^{n-p}}\Bigl)^{1/(p-1)}\frac{dr}{r}\Bigl)d\omega \leq \frac{1}{c} r^{n-p+\epsilon}.
\end{equation}\end{lem}

\begin{proof}  The proof we sketch is an adaptation of the dyadic shifting argument from \cite{JVFund}, see the argument from Lemma 4.8 through Corollary 4.11.  We reduce matters to a dyadic version.
Recall the definition of the dyadic Wolff potential from (4.12) of \cite{JVFund}.  From the argument of Lemma 4.9 of \cite{JVFund}, in particular (4.19), it follows that in order to prove (\ref{expinteps}) it suffices to prove that:
\begin{equation}\begin{split}\label{mfoldlebest}
I & = \int_{B(x,r)} \sum_{z\in Q_1\in \mathcal{D}}  \Bigl(\frac{|Q_1\cap B(x,r)|_{\sigma}}{\ell(Q_1)^{n-p}}\Bigl)^{1/(p-1)}\cdot \sum_{z\in Q_2\subset Q_1}\Bigl(\frac{|Q_2\cap B(x,r)|_{\sigma}}{\ell(Q_2)^{n-p}}\Bigl)^{1/(p-1)}\\
&\cdots \sum_{z\in Q_m\subset Q_{m-1}} \Bigl(\frac{|Q_m\cap B(x,r)|_{\sigma}}{\ell(Q_m)^{n-p}}\Bigl)^{1/(p-1)}d\omega \\
& \leq A(n,p, C(\omega), C(\sigma)) (B(n,p) C(\omega)C(\sigma))^m r^{n-p+\epsilon},
\end{split}\end{equation}
where all sums are taken over dyadic cubes, see Section \ref{disccarlsec} above for notation.
Let $\tilde\sigma$ be the measure $\tilde{\sigma}(E) = \sigma(E\cap B(x,r))$ and $\tilde\omega(E) = \omega(E\cap B(x,r))$.  Note that any ball $B(x,r)$ is contained in the union of at most $2^n$ dyadic cubes $\{P^j\}_{j=1, \dots 2^n}$ of radius $\lesssim r$.  The argument to show (\ref{mfoldlebest}) splits into two cases, depending on whether $1<p\leq2$ or $p>2$.

First suppose that $1<p\leq 2$, then note that, for any dyadic cube $P\in \mathcal{D}$, it follows:
\begin{equation}\begin{split}\label{1stitstepleq2}\sum_{Q\subset P}& \Bigl(\frac{|Q|_{\tilde{\sigma}}}{\ell(Q_1)^{n-p}}\Bigl)^{1/(p-1)}\omega(Q\cap B(x,r))\\
& = \int_{P\cap B(x,r)} \sum_{z\in Q\subset P} \Bigl(\frac{|Q|_{\tilde{\sigma}}}{\ell(Q)^{n-p}}\Bigl)^{(2-p)/(p-1)}\frac{\omega(Q\cap B(x,r))}{\ell(Q)^{n-p}}d\sigma\end{split}\end{equation}
Using the property on $\sigma$, we bound the right hand side of (\ref{1stitstepleq2}) by:
\begin{equation}\begin{split} C(\sigma)^{\frac{2-p}{p-1}}\int_{P\cap B(x,r)} \sum_{z\in Q\subset P} \frac{\omega(Q\cap B(x,r))}{\ell(Q)^{n-p}}d\sigma\lesssim C(\omega)C(\sigma)^{\frac{2-p}{p-1}}|P|_{\tilde{\sigma}}r^{\epsilon},
\end{split}\end{equation}
where we have estimated the sum in the integrand in the following way:  for each $z\in B(x,r)$, $z\in P^j$ for some $j$, with $P^j$ the dyadic covering of $B(x,r)$ as above.  Then:
\begin{equation}\begin{split}\label{pleqintegrandest}
 \sum_{z\in Q\subset P} \frac{\omega(Q\cap B(x,r))}{\ell(Q)^{n-p}} & = \sum_{z\in Q \subset P^j} \cdots + \sum_{P^j \subset Q\subset P} \cdots\\
 & \leq C(\omega)\sum_{z\in Q\subset P^j} \ell(Q)^{\epsilon} + \omega(B(x,r)) \sum_{Q\supset P^j} \ell(Q)^{p-n} \\
 & \lesssim C(\omega)\ell(P^j)^{\epsilon} + r^{n-p+\epsilon}C(\omega) \ell(P^j)^{p-n} \lesssim r^{\epsilon}.
\end{split}\end{equation}

Applying (\ref{1stitstepleq2}) into the left hand side of (\ref{mfoldlebest}), it follows that:
\begin{equation}\begin{split}\nonumber
I\lesssim  & C(\sigma)^{\frac{2-p}{p-1}}  r^{\epsilon} \sum_{Q_1\in \mathcal{Q}}\Bigl(\frac{|Q_1\cap B(x,r)|_{\sigma}}{\ell(Q_1)^{n-p}}\Bigl)^{1/(p-1)} \sum_{Q_2\subset Q_1}\Bigl(\frac{|Q_2\cap B(x,r)|_{\sigma}}{\ell(Q_2)^{n-p}}\Bigl)^{1/(p-1)}\\
&\cdots \sum_{Q_{m-1}\subset Q_{m-2}} \Bigl(\frac{|Q_{m-1}\cap B(x,r)|_{\sigma}}{\ell(Q_{m-1})^{n-p}}\Bigl)^{1/(p-1)} |Q_{m-1}|_{\tilde\sigma},
\end{split}\end{equation}
Now, following (4.19) of \cite{JVFund}, we conclude that:
\begin{equation}\begin{split}\nonumber
I \leq B^m & r^{\epsilon} C(\sigma)^{\frac{2-p}{p-1}} C(\sigma)^{m-1}\sigma(B(x,r)) \leq r^{n-p+\epsilon}B^mC(\sigma)^mC(\sigma)^{\frac{2-p}{p-1}},
\end{split}\end{equation}
for a suitable constant $B=B(n,p)>0$.  We conclude that (\ref{mfoldlebest}) holds if $p\leq 2$.

Let us now consider the case when $p\geq 2$.  This is slightly more involved.  In this case note:
\begin{equation}\label{omegaobs}|Q|_{\tilde\omega} \leq C(\omega)^{1/(p-1)}\bigl(\min(\ell(Q),r)^{n-p+\epsilon}\big)^{1/(p-1)}|Q|_{\tilde\omega}^{(p-2)/(p-1)}.\end{equation}
Hence, if we set:
\begin{equation}\begin{split}\label{1stitstepgeq2}II=& \sum_{Q\subset P} \Bigl(\frac{|Q|_{\tilde{\sigma}}}{\ell(Q)^{n-p}}\Bigl)^{1/(p-1)}\omega(Q\cap B(x,r)),
\end{split}\end{equation}
then, using H\"{o}lder's inequality in (\ref{1stitstepgeq2}), it follows:
\begin{equation}\begin{split}\nonumber
II& \leq C(\omega)^{1/(p-1)} \Bigl(\sum_{Q\subset P} |Q|_{\tilde\sigma}\Bigl(\frac{\min(\ell(Q),r)^{n-p+\epsilon}}{\ell(Q)^{n-p}}\Big)^{1/(p-1)}\Bigl)^{(2-p)/(p-1)}\\
&\;\;\;\;\;\;\;\;\;\;\;\;\;\cdot \Bigl(\sum_{Q\subset P}|Q\cap B(x,r)|_{\omega}\Bigl(\frac{\min(\ell(Q),r)^{n-p+\epsilon}}{\ell(Q)^{n-p}}\Big)^{1/(p-1)}\Bigl)^{1/(p-1)}.\\
\end{split}\end{equation}
Subsequently, by Fubini's theorem:
\begin{equation}\begin{split}\nonumber
II& \leq \Bigl(\int_{P\cap B(x,r)}\sum_{z\in Q\subset P}\Bigl(\frac{\min(\ell(Q),r)^{n-p+\epsilon}}{\ell(Q)^{n-p}}\Big)^{1/(p-1)}d\sigma \Bigl)^{\frac{p-2}{p-1}} \\
&\;\;\;\; \;\;\;\;\;\;\;\;\;\cdot\Bigl(\int_{P\cap B(x,r)}\sum_{z\in Q\subset P}\Bigl(\frac{\min(\ell(Q),r)^{n-p+\epsilon}}{\ell(Q)^{n-p}}\Big)^{1/(p-1)}dx \Bigl)^{\frac{1}{p-1}}\\
& \lesssim C(\omega)^{1/(p-1)}C(\sigma)^{1/(p-1)}\min(\ell(P),r)^{\frac{n-p+\epsilon}{p-1}}\omega(P\cap B(x,r))^{(p-2)/(p-1)}.
\end{split}\end{equation}
The implicit constants here depend on $n$, $p$ and $\epsilon$. In the previous calculation the integrands have been estimated in a similar manner as in (\ref{pleqintegrandest}) above, splitting the sums over the small cubes and large cubes (compared with $r$).


One can then iterate this calculation to estimate:
\begin{equation}\begin{split}\nonumber
I & 
\leq (C(p,n, \epsilon))^m(C(\omega)C(\sigma))^{m/(p-1)}r^{\frac{n-p+\epsilon}{p-1}}\omega(B(x,r))^{(p-2)/(p-1)}\\
&\leq (C(p,n, \epsilon))^m(C(\omega)C(\sigma))^{m/(p-1)}r^{n-p+\epsilon}.
\end{split}\end{equation}
In conclusion, we have asserted that (\ref{mfoldlebest}) holds, and therefore the lemma is proved.
\end{proof}

\section{A  fully nonlinear analogue: the $k$-Hessian}\label{Hessian}  In this section we briefly describe how one can obtain a $k$-Hessian analogue of our primary results.  Let $1\leq k < n/2$, then we will be interested in the problem:
\begin{equation}\label{hessentire}
\begin{cases}
\textrm{F}_k(-u) = \sigma u^k + \omega \text{ in } \mathbf{R}^n,\\
\; \inf_{x\in \mathbf{R}^n}u(x) = 0.
\end{cases}
\end{equation}
Here, $\textrm{F}_k$ is the $k$-Hessian operator (see \cite{CNS1}), defined for smooth functions $u$ by:
$$F_k(u) = \sum_{1\leq i_1 < \cdots < i_k \leq n} \lambda_{i_1}\dots \lambda_{i_k},$$
with $\lambda_1, \dots \lambda_n$ denoting the eigenvalues of the Hessian matrix $D^2u$.  We will use the notion of $k$-convex  functions, introduced by Trudinger and Wang \cite{TW2}  to state our results.  We say that $u$ is $k$-convex in $\Omega$ if $u:\Omega \rightarrow [-\infty, \infty)$ is upper semicontinuous and satisfies $F_k(u)\geq 0$ in the viscosity sense, i.e. for any $x\in \Omega$, $F_k(q)(x)\geq 0$ for any quadratic polynomial $q$ so that $u-q$ has a local finite maximum at $x$.  We will seek solutions $u \ge 0$ of (\ref{hessentire}) so that $-u$ is $k$-convex. This convention allows us to state  our results in the form  analogous  to the quasilinear case. 
Equivalently (see \cite{TW2}), we may define $k$-convex functions by a comparison principle: an upper semicontinuous function $u : \Omega \rightarrow [-\infty, \infty)$ is $k$-convex in $\Omega$ if for every open set $D \subset\subset \Omega$, and $v\in C^2_{\text{loc}}(D)\cap C(\bar{D})$ with $F_k(v)\geq 0$ in $D$, then 
%

The necessary condition on $\sigma$ is now considered in terms of the $k$-Hessian capacity, introduced again in \cite{TW2}.  The $k$-Hessian capacity of a compact set $E$ is defined by:
\begin{equation}\label{Hesscap}
\text{cap}_k (E) = \sup \{\mu_k[u](E) \, : \, u \text{ is } k\text{-convex in } \mathbf{R}^n, \; -1<u<0 \}.
\end{equation}
Here $\mu_k[u]$ is the $k$-Hessian measure of $u$, see Theorem \ref{Hessweakcont} below.  The two local potentials that will be relevant for our bounds will be the fractional linear potential $\mathbf{I}_{2k}^r(d\sigma)$ and $\mathbf{W}_{\frac{2k}{k+1}, k+1}^r(d\sigma)$, as defined in (\ref{Riesz}) and (\ref{Wolff}) respectively.
\begin{thm}\label{mainthmHess} Let $1\leq k <n/2$, and let $\alpha = 2k/(k+1)$.  Suppose that $u\ge 0$ is a solution of (\ref{hessentire}) so that 
$-u$ is $k$-convex. Then $\sigma$ satisfies:
\begin{equation}\label{hesscapnec}\sigma(E) \leq C(\sigma) \text{cap}_k(E) \; \text{ for all compact sets } E\subset \mathbf{R}^n.
\end{equation}
 In addition there is a constant $c=c(n,k)>0$ so that $u$ satisfies:
\begin{equation}\begin{split}\label{hesslow}
u(x) \geq c & \int_{0}^{\infty}\Bigl(\frac{e^{c\mathbf{W}^r_{\alpha,k+1}(d\sigma)(x)}}{r^{n-\alpha s}} \int_{B(x,r)} e^{c\mathbf{I}^r_{2k}(d\sigma)(z)}d\omega(z)\Bigl)^{\frac{1}{k}}\frac{dr}{r},
\end{split}
\end{equation}
for all $x\in \mathbf{R}^n$.\\
Conversely, assuming the right hand side of (\ref{hessup}) below is finite for some $x\in \Omega$ and $c>0$, there exists a positive constant $C_0 = C_0(n,k,c)>0$, such that if $\sigma$ satisfies (\ref{hesscapnec}), with constant $C(\sigma)<C_0$, then there exists a solution $u\ge 0$ of (\ref{hessentire}) so that $-u$ is $k$-convex, and 
\begin{equation}
\label{hessup}
u(x) \leq c_1 \int_{0}^{\infty}\Bigl(\frac{e^{c\mathbf{W}^r_{\alpha,k+1}(d\sigma)(x)}}{r^{n-\alpha s}} \int_{B(x,r)}e^{c\mathbf{W}^r_{\alpha, k+1}(d\sigma)(z)}d\omega(z)\Bigl)^{\frac{1}{2k}}\frac{dr}{r},
\end{equation}
for all $x\in \mathbf{R}^n$. Here $c_1 = c_1(n,k,c)>0$.
\end{thm}

Let $\Phi ^k(\Omega)$ be the set of $k$-convex functions such that $u$ is not identically infinite in each component of $\Omega$.  The following weak continuity result is key to us. 
\begin{thm}\cite{TW2}\label{Hessweakcont} Let $u\in \Phi^k(\Omega)$. Then there is a  nonnegative Borel measure $\mu_k[u]$ in $\Omega$ such that 
\begin{itemize}
\item $\mu_k[u] = F_k(u)$ whenever $u\in C^2(\Omega)$, and
\item If $\{u_m\}_m$ is a sequence in $\Phi^k(\Omega)$ converging in $L^1_{\text{loc}}(\Omega)$ to a function $u$, then the sequence of measures $\{\mu_k[u_m]\}_m$ converges weakly to $\mu_k[u]$. 
\end{itemize}
\end{thm}

The measure $\mu_k[u]$ associated to $u\in \Phi^k(\Omega)$ is called the \textit{Hessian measure} of $u$.  Hessian measures were used by Labutin \cite{Lab1} to deduce Wolff's  potential estimates for a $k$-convex function in terms of its Hessian measure.  The following global version of Labutin's estimate is deduced from his result in \cite{PV}:
\begin{thm}\cite{PV} \label{Hesswolff}
Let $1\leq k \leq n$, and suppose that $u\geq 0$ is such that $-u\in \Phi^k(\Omega)$ and $\inf_{x\in \mathbf{R}^n}u(x) =0$.  Then, if $\mu=\mu_k[u]$, there is a positive constant $K$, depending on $n$ and $k$, such that:
$$c_1\mathbf{W}_{\frac{2k}{k+1}, k+1}\mu(x) \leq u(x) \leq c_2 \mathbf{W}_{\frac{2k}{k+1}, k+1}\mu(x), \qquad x \in \mathbf{R}^n.
$$
\end{thm}

Using Theorem \ref{Hesswolff}, we see that there is a positive constant $C>0$, so that any solution $u$ of (\ref{hessentire}) satisfies the following estimate:
\begin{equation}\label{intineqHess}
u(x)\geq C \int_0^{\infty} \Bigl(\frac{1}{r^{n-2k}}\int_{B(x,r)} u^k(z)d\sigma(z)\Bigl)^{1/k}\frac{dr}{r} +C \int_0^{\infty} \Bigl(\frac{\omega(B(x,r))}{r^{n-2k}}\Bigl)^{1/k}\frac{dr}{r} 
\end{equation}
By iterating (\ref{intineqHess}), one obtains a lower bound for $u$ in terms of a formal Neumann series of iterated operators, cf. Lemma \ref{lowbdsum}.  The iterates can be estimated using the techniques of Lemmas \ref{lowerouterbuildup} and \ref{sgeq2lowbdlem}.  One then arrives at (\ref{hesslow}).  The calculation was carried out in full for the fundamental solution in \cite{JVFund}.

Regarding (\ref{hessup}), one can construct a supersolution (\ref{intineqHess}) using precisely the same techniques as Lemma \ref{sgeq2upbdlem2}.  In other words, the integral inequalities that govern the $k$-Hessian have the same behaviour as the integral inequalities for the $p$-Laplacian when $p\geq 2$.  From the integral supersolution, we conclude using Theorem \ref{Hessweakcont} that there exists a solution of (\ref{hessentire}) so that (\ref{hessup}) holds.  This mimics the iterative argument of \cite{JVFund}, and is similar in nature to the arguments spelled out in Section \ref{existence} above.

\appendix 
\section{A tail estimate for nonlinear potentials: \\Proof of Lemma \ref{lemtailest}} \label{appendix1}
\begin{lem} Let $\sigma$ be a Borel measure satisfying:
\begin{equation}\label{potestap}\sigma (B(x,r)) \leq C r^{n-\alpha s}, \text{ for all balls } B(x,r).
\end{equation} Then there is a positive constant $C=C(n, \alpha, s, \sigma)>0$ so that for all $x \in \mathbb{R}^n$ and  $y \in B(x,t)$, $t>0$, it follows: 
\begin{equation} \abs{\int_t^{\infty} \left [ \Bigl(\frac{\sigma(B(x,r))}{r^{n-\alpha s}}\Bigr)^{\frac{1}{s-1}} - \Bigl(\frac{\sigma(B(y,r))}{r^{n-\alpha s}}\Bigr)^{\frac{1}{s-1}} \right ] \frac{dr}{r}} \leq C . 
\end{equation}
\end{lem}

\begin{proof}
Without loss of generality, suppose that:
\begin{equation}\nonumber
\int_t^{\infty} \left [ \Bigl(\frac{\sigma(B(x,r))}{r^{n-\alpha s}}\Bigr)^{\frac{1}{s-1}} - \Bigl(\frac{\sigma(B(y,r))}{r^{n-\alpha s}}\Bigr)^{\frac{1}{s-1}} \right ] \frac{dr}{r} > 0 . 
\end{equation}
We want to rearrange the integrand so it is nonnegative.  To this end, we define two sets:
\begin{align}\nonumber A = \{z \in \mathbf{R}^n : \abs{x-z} \leq \abs{y-z} \},\,\, \text{ and } \,\,\, B = \{z \in \mathbf{R}^n : \abs{y-z} < \abs{x-z} \} . 
\end{align}
Then if $z \in B$ and $\abs{z-x} < r$, we have that $\abs{y-z} < r$, and thus $B(x,r) \cap B \subset B(y,r) \cap B$ so that:
\begin{equation}\label{elemobs1}
\sigma(B(x,r) \cap B) \leq \sigma(B(y,r) \cap B), \text{ and:}
\end{equation}
\begin{equation}\label{elemobs2}
\sigma(B(y,r) \cap A) \leq \sigma(B(x,r) \cap A) . 
\end{equation}
Using (\ref{elemobs1}) gives:
\begin{equation}\begin{split}\nonumber
\int_t^{\infty}& \left [ \Bigl(\frac{\sigma(B(x,r))}{r^{n-\alpha s}}\Bigr)^{\frac{1}{s-1}} - \Bigl(\frac{\sigma(B(y,r))}{r^{n-\alpha s}}\Bigr)^{\frac{1}{s-1}}  \right ] \frac{dr}{r} \\
& \leq \int_t^{\infty} \Bigl(\frac{\sigma(B(x,r)\cap A) + \sigma(B(y,r)\cap B)}{r^{n-\alpha s}}\Bigr)^{\frac{1}{s-1}} \\ 
&- \Bigl(\frac{\sigma(B(y,r)\cap A) + \sigma(B(x,r)\cap B)}{r^{n-\alpha s}}\Bigr)^{\frac{1}{s-1}} \frac{dr}{r} = \int_t^{\infty} \left [ I^{\frac{1}{s-1}} -  II^{\frac{1}{s-1}} \right] \frac{dr}{r} . 
\end{split}
\end{equation}
From (\ref{elemobs1}) and (\ref{elemobs2}) it immediately follows that the integrand in nonnegative, i.e. that $I \geq II$.

The proof now splits into two cases, when $1<s<2$ and when $s\geq 2$.  First suppose $1<s < 2$, then note the following elementary inequality; for $a,b \in (0, \infty)$ with $a>b$, and $\gamma \geq 1$:
\begin{equation}\label{elemab1} a^{\gamma} - b^{\gamma} \leq \gamma a^{\gamma-1}(a-b) \end{equation}
Plugging $I$ and $II$ into (\ref{elemab1}) yields: $ I^{\frac{1}{s-1}} -  II^{\frac{1}{s-1}} \leq \frac{1}{s-1} (I -II)I^{\frac{2-s}{s-1}} \leq C( I - II ).$ Here we have used the estimate (\ref{potestap}) in the last inequality, noting that $2-s >0$.

As a result (in case $1<s<2$), the lemma will follow from the following inequality: 
\begin{equation}\label{redpleq2}\begin{split}
\int_t^{\infty} & \frac{\sigma(B(x,r)\cap A) + \sigma(B(y,r)\cap B)}{r^{n-\alpha s}} \\& - \frac{\sigma(B(y,r)\cap A) + \sigma(B(x,r)\cap B)}{r^{n-\alpha s}}  \frac{dr}{r} \leq C . 
\end{split}\end{equation}
Let us now split $\sigma$ into $\sigma_1 = \sigma \cdot \chi_{ \mathbf{R}^n \backslash B(x,2t)}$ and $\sigma_2 = \sigma \cdot \chi_{B(x,2t)}$ and if we can control the left hand side of (\ref{redpleq2}) with either $\sigma_1$ or $\sigma_2$ in place of $\sigma$ then we are done.\\
The estimate for $\sigma_2$ is a straightforward application of  (\ref{potestap}):
\begin{equation}\begin{split}\nonumber
 \int_t^{\infty} & \frac{\sigma_2(B(x,r)\cap A) + \sigma_2(B(y,r)\cap B)}{r^{n-\alpha s}} - \frac{\sigma_2(B(y,r)\cap A) + \sigma_2(B(x,r)\cap B)}{r^{n-\alpha s}} \frac{dr}{r}\\
&\leq C \sigma(B(x,2t)) \int_{t}^{\infty}\frac{1}{r^{n-\alpha s}}\frac{dr}{r} \leq C \frac{\sigma(B(x,2t)) }{(2t)^{n-\alpha s}} \leq C
\end{split}\end{equation}
where (\ref{potestap}) has been used in this last inequality.

We now move onto the estimate for $\sigma_1$.  First note that if $r<t$ and $y\in B(x,t)$, then $B(y,r) \subset B(x,2t)$ and so $\sigma_1(B(y,r))=0$. This allows us to extend the integration to over the half line:
\begin{equation}\begin{split}\nonumber
\int_t^{\infty}& \frac{\sigma_1(B(x,r)\cap A) + \sigma_1(B(y,r)\cap B)}{r^{n-\alpha s}}   - \frac{\sigma_1(B(y,r)\cap A) + \sigma_1(B(x,r)\cap B)}{r^{n-\alpha s}} \frac{dr}{r}\\
& = \frac{1}{n-\alpha s} \int _{\mathbf{R}^n}  \left [ \frac{\chi_A(z)}{\abs{x-z}^{n-\alpha s}} - \frac{\chi_A(z)}{\abs{y-z}^{n-\alpha s}}  +   \frac{\chi_B(z)}{\abs{y-z}^{n-p}} - \frac{\chi_B(z)}{\abs{x-z}^{n-\alpha s}} \right ] d \sigma_1(z)\\
& = \frac{1}{n-\alpha s} \int _{\mathbb{R}^n\backslash B(x,2t)}  \abs{\frac{1}{\abs{x-z}^{n-\alpha s}} - \frac{1}{\abs{y-z}^{n-\alpha s}} } d \sigma(z)
\end{split}\end{equation}
Let $z \notin B(x,2t)$, then whenever $y \in B(x,t)$, it is easy to see that:
\begin{equation}\label{ballequiv} \frac{1}{2}\abs{y-z} \leq \abs{x-z} \leq 2 \abs{y-z} . 
\end{equation}
Note the following elementary inequality. For $a,b \in (0, \infty)$ with $a>b$, and $\gamma \geq 0$:
\begin{equation}\label{elemab2} a^{\gamma} - b^{\gamma} \leq \gamma (a^{\gamma-1} + b^{\gamma-1})(a-b) .\end{equation}
Due to (\ref{ballequiv}) and (\ref{elemab2}), and since $y\in B(x,t)$, it follows:
\begin{equation}\nonumber\begin{split}
\left \vert \frac{1}{\abs{x-z}^{n-\alpha s}} - \frac{1}{\abs{y-z}^{n-\alpha s}}\right \vert  
 \leq C \frac{\abs{x-y}}{\abs{x-z}^{n-\alpha s+1}} \leq C \frac{t}{\abs{x-z}^{n-\alpha s+1}} . 
\end{split}
\end{equation}
Combining these observations, we obtain:
\begin{equation}\nonumber\begin{split}
\int_{\mathbb{R}^n \backslash B(x,2t)} \abs{\frac{1}{\abs{x-z}^{n-\alpha s}} - \frac{1}{\abs{y-z}^{n-\alpha s}}}d \sigma(z) & \leq C \int_{\mathbb{R}^n \backslash B(x,2t)}  \frac{t}{\abs{x-z}^{n-\alpha s+1}} d \sigma (z)\\
& \leq C t \int_{2t}^{\infty} \frac{\sigma(B_{r}(x))}{r^{n-\alpha s}} \frac{dr}{r^2} \leq C t\int_{2t}^{\infty} \frac{dr}{r^2}\leq C . 
\end{split}
\end{equation}
As a result, the lemma is proved in the case $1<s \leq 2$. 

We now move onto the $s\geq 2$ case.  First recall that with $I$ and $II$ as before, we have $I \geq II$, and hence $I^{\frac{1}{s-1}} - II^{\frac{1}{s-1}} \leq (I - II)^{\frac{1}{s-1}}$.  This implies that:
\begin{equation}\begin{split}\nonumber
\int_t^{\infty}& \left [ \Bigl(\frac{\sigma(B(x,r))}{r^{n-\alpha s}}\Bigr)^{\frac{1}{s-1}} - \Bigl(\frac{\sigma(B(y,r))}{r^{n-\alpha s}}\Bigr)^{\frac{1}{s-1}} \right ] \frac{dr}{r} \\
& \leq \int_t^{\infty} \Bigl(\frac{\sigma(B(x,r)\cap A) + \sigma(B(y,r)\cap B)}{r^{n-\alpha s}} - \frac{\sigma(B(y,r)\cap A) + \sigma(B(x,r)\cap B)}{r^{n-\alpha s}}\Bigr)^{\frac{1}{s-1}} \frac{dr}{r} . 
\end{split}
\end{equation}
Let $\epsilon > 0 $ small enough so that $\epsilon(s-2) < \min(n-\alpha s, 1)$ . Then, by H\"{o}lder's inequality:
\begin{equation}\begin{split}\nonumber
\int_t^{\infty} (I -  II)^{\frac{1}{s-1}}  \frac{dr}{r} & \leq C t^{\epsilon (\frac{1}{s-1}-1)}  \Bigl(\int_t^{\infty}  (I -  II)r^{\epsilon(s-1)}  \frac{dr}{r^{1+\epsilon}}\Bigr)^{\frac{1}{s-1}}\\
&= C t^{\epsilon (\frac{1}{s-1}-1)}  \Bigl(\int_t^{\infty} \Bigl(\frac{\sigma(B(x,r)\cap A) + \sigma(B(y,r)\cap B)}{r^{n-\alpha s}}\\
& - \frac{\sigma(B(y,r)\cap A) + \sigma(B(x,r)\cap B)}{r^{n-\alpha s}}\Bigr)r^{\epsilon(s-1)}  \frac{dr}{r^{1+\epsilon}}\Bigr)^{\frac{1}{s-1}} . 
\end{split}
\end{equation}

We wish to bound the right hand side by a constant.  To this end we will split the measure $\sigma$ as before into $\sigma_1$ and $\sigma_2$.  The following estimate for $\sigma_2$ follows easily using (\ref{potestap}):
\begin{equation}\begin{split}\nonumber
t^{\epsilon (\frac{1}{s-1}-1)} & \Bigl(\int_t^{\infty} \Bigl(\frac{\sigma_2(B(x,r)\cap A) + \sigma_2(B(y,r)\cap B)}{r^{n-\alpha s}}\\&- \frac{\sigma_2(B(y,r)\cap A) + \sigma_2(B(x,r)\cap B)}{r^{n-\alpha s}}\Bigr)r^{\epsilon(s-1)}  \frac{dr}{r^{1+\epsilon}}\Bigr)^{\frac{1}{s-1}} \leq C,  
\end{split}
\end{equation}

We now concentrate on the $\sigma_1$ estimate. First we note that we may extend the domain of integration over the whole half line and use Fubini's theorem  as in the $1<s\leq 2$ case to find that 
\begin{equation}\begin{split}\nonumber
t^{\epsilon (\frac{1}{s-1}-1)} &  \Bigl(\int_t^{\infty} \Bigl(\frac{\sigma(B(x,r)\cap A) + \sigma(B(y,r)\cap B)}{r^{n-\alpha s}}\\
& - \frac{\sigma(B(y,r)\cap A) + \sigma(B(x,r)\cap B)}{r^{n-\alpha s}}\Bigr)r^{\epsilon(s-1)}  \frac{dr}{r^{1+\epsilon}}\Bigr)^{\frac{1}{s-1}}\\
& \leq C t^{\epsilon (\frac{1}{s-1}-1)}\Bigl(\int_{\mathbb{R}^n\backslash B(x,2t)}\abs{\frac{1}{\abs{x-z}^{n-\alpha s-\epsilon(s-1)+\epsilon}} - \frac{1}{\abs{y-z}^{n-\alpha s-\epsilon(s-1)+\epsilon}}}d\sigma(z)\Bigr)^{\frac{1}{s-1}}\\
& = III . 
\end{split}
\end{equation}

Now by adapting the previous argument in the $s\leq 2$ case, we have: 
\begin{equation}\nonumber
\abs{\frac{1}{\abs{x-z}^{n-\alpha s-\epsilon(s-1)+\epsilon}} - \frac{1}{\abs{y-z}^{n-\alpha s-\epsilon(s-1)+\epsilon}}} \leq C \frac{t}{\abs{x-z}^{n-\alpha s-\epsilon(s-1)+\epsilon + 1}} . 
\end{equation}
Hence 
\begin{equation}\begin{split}\nonumber
III & \leq C t^{\epsilon (\frac{1}{s-1}-1)}\Bigl(\int_{\mathbb{R}^n\backslash B(x,2t)}\frac{t}{\abs{x-z}^{n-\alpha s-\epsilon(s-1)+\epsilon + 1}}d\sigma(z)\Bigr)^{\frac{1}{s-1}}\\
& \leq C t^{\epsilon (\frac{1}{s-1}-1)}\Bigl(\int_{2t}^{\infty}\frac{t\, \sigma(B(x,r))}{r^{n-\alpha s-\epsilon(s-1)+\epsilon + 1}}\frac{dr}{r}\Bigr)^{\frac{1}{s-1}} \leq C , 
\end{split}
\end{equation}
where in the last inequality we have used (\ref{potestap}), then we are left with a convergent integral by choice of $\epsilon$.   This completes the proof in the case $s \geq 2$. 
\end{proof}

\section{Duality in discrete Littlewood--Paley spaces}\label{duality}

\begin{lem}\label{lpdual}  Fix a dyadic cube $P\in \mathcal{Q}$, and let $s>1$.   There is a constant $c = c(s)>0$ such that:
\begin{equation}\begin{split}\label{2sided}
\frac{1}{c} \sup_{\{\mu_Q\}_{Q\subset P}} \sum_{Q\subset P} \lambda_Q \mu_Q |Q|_{\sigma}  \leq \int_P \Big[\sum_{x\in Q\subset P}&\lambda_Q^{s}\Bigl]^{1/s}d\sigma(x) \\
&\leq c \sup_{\{\mu_Q\}_{Q\subset P}}\sum_{Q\subset P}\lambda_Q \mu_Q|Q|_{\sigma},
\end{split}
\end{equation}
where the supremum is taken over all sequences $\{\mu_Q\}_{Q\subset P}$ satisfying:
\begin{equation}\label{appendspace}\sup_{Q\subset P}\frac{1}{\abs{Q}_{\sigma}}\sum_{R\subset Q, \, R\in \mathcal{Q}}\mu_R^{s'}\abs{R}_{\sigma} \leq 1.
\end{equation}
\end{lem}

\begin{proof}  The lower estimate is known to be a (nontrivial) consequence of the Carleson measure theorem, Theorem \ref{Carlemb}.  A proof can be found in Theorem 4, part (b) of \cite{Verb96}, see also Theorem 3.1 of \cite{NTV03}. To prove the upper estimate, write:
\begin{equation}\begin{split}\nonumber\int_P\Big\{\sum_{x\in Q\subset P} & \lambda_Q^{s}\Bigl\}^{1/s}d\sigma(x) = \sum_{Q\subset P} \lambda_Q^{s}\int_Q\Big\{\sum_{x\in R\subset P} \lambda_R^{s}\Bigl\}^{\frac{1}{s}-1}d\sigma(x)\\
& =\sum_{Q\subset P}\lambda_Q\abs{Q}_{\sigma}\frac{1}{\abs{Q}_{\sigma}}\int_Q\Big\{\frac{\lambda_Q^{s}}{\sum_{x\in R\subset P} \lambda_R^{s}}\Bigl\}^{1-1/s}d\sigma(x)
\end{split}\end{equation}
Let us set:
$$\mu_Q = \frac{1}{\abs{Q}_{\sigma}}\int_Q\Big\{\frac{\lambda_Q^{s}}{\sum_{x\in R\subset P} \lambda_R^{s}}\Bigl\}^{1/s'}d\sigma(x).
$$
It remains to see that 
$\{\mu_Q\}$ is admissible for (\ref{appendspace}).  This is a simple consequence of H\"{o}lder's inequality and interchanging summation and integration.  Indeed, for any $Q\subset P$, and with this choice of $\mu_Q$, it follows:
\begin{equation}\begin{split}\nonumber
\sum_{R\subset Q, \, R\in \mathcal{Q}}\mu_R^{s'}\abs{R}_{\sigma} &\leq \sum_{R\subset Q, \, R\in \mathcal{Q}} \int_R\frac{\lambda_R^{s}}{\sum_{x\in S\subset P} \lambda_S^{s}}d\sigma(x)\\
& = \int_Q \frac{\sum_{x\in R\subset Q}\lambda_R^{s}}{\sum_{x\in S\subset P} \lambda_S^{s}}d\sigma(x)\leq \abs{Q}_{\sigma},
\end{split}\end{equation}
as required.  
\end{proof}


\begin{thebibliography}{DMMOP99}

\bibitem[ADP06]{ADP06} B. Abdellaoui, A. Dall'Aglio and I. Peral, \emph{Some remarks on elliptic problems with critical growth in the gradient}, J. Diff. Equations, \textbf{222} (2006), 21--62.

\bibitem[AHBV09]{AHBV}
H.~Abdul-Hamid and M-F. Bidaut-V\'{e}ron, \emph{On the connection between two
  quasilinear elliptic problems with source terms of order 0 or 1}, Commun. Contemp. Math. \textbf{12} (2010),
  727--788. 


\bibitem[AH96]{AH}
D.~R. Adams and L.~I. Hedberg, \emph{Function Spaces and Potential Theory},
  Grundlehren der mathematischen Wissenschaften \textbf{314}, Springer, Berlin, 1996.


\bibitem[BBGPV]{BBGVP95}
P.~B\'{e}nilan, L.~Boccardo, R.~Gariepy, M.~Pierre, and J.~Vazquez, \emph{An $L^1$ theory of exsitence and uniqueness of solutions of nonlinear elliptic equations}, Ann. Scuola Norm. Super. Pisa \textbf{22}  (1995), 
241--273. 


\bibitem[CNS85]{CNS1}
L.~Caffarelli, L.~Nirenberg, and J.~Spruck, \emph{The {D}irichlet problem for
  nonlinear second-order elliptic equations {III}. {F}unctions of the
  eigenvalues of the {H}essian}, Acta Math. \textbf{155} (1985), 
  261--301.

\bibitem[COV00]{COV1}C.~Cascante, J.~M. Ortega, and I.~E. Verbitsky, \emph{Trace inequalities of Sobolev type in the upper triangle case}, 
Proc. London Math. Soc.  \textbf{80} (2000),  391--414. 

\bibitem[COV04]{COV}
C.~Cascante, J.~M. Ortega, and I.~E. Verbitsky, \emph{Nonlinear potentials and
  two weight trace inequalities for general dyadic and radial kernels}, Indiana
  Univ. Math. J. \textbf{53} (2004), 845--882.


\bibitem[CZ95]{CZ95} K. L. Chung and Z. Zhao,  
\emph{From Brownian Motion to Schr\"odinger's Equation}.  Grundlehren der
math. Wissenschaften \textbf{312}, Springer,
Berin--Heidelberg, 1995.

\bibitem[DMM97]{DMM97}
G. Dal Maso and A. Malusa, \emph{Some properties of reachable solutions of nonlinear elliptic equations with measure data,}  
Ann. Scuola Norm. Super. Pisa \textbf{25} (1997),  375--396.

\bibitem[DMMOP]{DMMOP99}
G.~Dal Maso, F.~Murat, L.~Orsina, and A.~Prignet, \emph{Renormalized solutions of elliptic equations with general measure data},
Ann. Scuola Norm. Super. Pisa  \textbf{28} (1999), 741--808. 

\bibitem[DM10]{DM10}  F. Duzaar and  G. Mingione, \emph{Gradient estimates via linear and nonlinear potentials}, J. Funct. Anal. \textbf{259} (2010), 2961--2998. 

\bibitem[DM11]{DM11}  F. Duzaar and  G. Mingione, \emph{Gradient estimates via nonlinear potentials}, Amer. J. Math. (to appear). 



\bibitem[FS71]{FS}
C.~Fefferman and E.~M. Stein, \emph{Some maximal inequalities}, Amer. J. Math. 
  \textbf{93} (1971), 107--115.


\bibitem[FM00]{FM1}
V.~Ferone and F.~Murat, \emph{Nonlinear problems having natural growth in the
  gradient: an existence result when the source terms are small}, Nonlinear
  Analysis \textbf{42} (2000), 1309--1326.



\bibitem[FV09]{FV2}
M.~Frazier and I.~E. Verbitsky, \emph{Solvability conditions for a discrete
  model of {S}chr{\"{o}}dinger's equation}, Analysis, PDE and Appl.,
The Vladimir Maz'ya Anniversary Volume.
Operator Theory: Adv. Appl. \textbf{179},   
Birk\"{a}user, 2010.


\bibitem[FV10]{FV1}
M.~Frazier and I.~E. Verbitsky, \emph{Global {G}reen{'}s function estimates}, 
Around the Research of Vladimir Maz'ya III, Analysis and Applications, Ed. Ari Laptev,  Intern. Math. Series \textbf{13}, Springer,  2010, 105--152.


\bibitem[FNV10]{FNV} 
M.~Frazier, F.~Nazarov, and I.~E. Verbitsky,
\emph{Global estimates for kernels of Neumann series, Green's
functions, and the conditional gauge}, Preprint (2011). 



\bibitem[GJ82]{GJ}
J.~Garnett and P.~Jones, \emph{{BMO} from dyadic {BMO}}, Pacific J. Math. \textbf{99} (1982), 351--371.

\bibitem[Gre02]{Gre02}
N.~Grenon, \textit{Existence results for semilinear elliptic equations with small measure data.}
Ann. Inst. H. Poincar\'{e}, Anal. Non Lin\'{e}aire \textbf{19} (2002), 1--11.

\bibitem[GT03]{GT03}
N. Grenon and C. Trombetti, \emph{Existence results for a class of nonlinear elliptic problems with p-growth in the gradient.} Nonlinear Anal. \textbf{52} (2003), 931--942.

\bibitem[HMV99]{HMV}
K.~Hansson, V.~G. Maz'ya, and I.~E. Verbitsky, \emph{Criteria of solvability for
  multidimensional {R}iccati equations}, Ark. Mat. \textbf{37} (1999), 87--120.


\bibitem[HKM06]{HKM}
J.~Heinonen, T.~Kilpel\"{a}inen, and O.~Martio, \emph{Nonlinear Potential
  Theory of Degenerate Elliptic Equations}, Dover Publ., 2006 (unabridged republication of 1993 edition, Oxford Universiy Press).

\bibitem[Hi48]{Hi48}
E.~Hille, \emph{Non-oscillation theorems,} 
Trans. Amer. Math. Soc. \textbf{64}, (1948), 234--252.

\bibitem[HJ11]{HJ11}
P. Honz\'\i k and B. J. Jaye, \emph{On the good-$\lambda$ inequality for nonlinear potentials}, Proc. Amer. Math. Soc. (to appear), arXiv:1105.6152.

\bibitem[JV10]{JVFund}
B.~J.~Jaye and I.~E.~Verbitsky, \emph{The fundamental solution of nonlinear operators with natural growth terms}, Ann. Scuola Norm. Super. Pisa  (to appear) arXiv:1002.4664. 

\bibitem[KK79]{KK79}J.~Kazdan and R.~Kramer, \emph{Invariant criteria for existence of second-order
  quasi-linear elliptic equations}, Commun. Pure. Appl. Math. \textbf{31} (1978),
  619--645.

\bibitem[Kil02]{K02}
T.~Kilpel\"{a}inen, \emph{$p$-Laplacian type equations involving measures}, Proc. ICM, Vol. III (Beijing, 2002), 167--176, Higher Ed. Press, Beijing, 2002.

\bibitem[KKT11]{KKT11} T. Kilpel\"{a}inen, T. Kuusi, and A. Tuhola-Kujanp\"{a}\"{a},  \emph{Superharmonic functions are locally renormalized solutions}, Ann. Inst. Henri Poincar\'{e} (to appear).


\bibitem[KM92]{KM92}
T.~Kilpel\"{a}inen and J.~Mal\'y, \emph{Degenerate elliptic equations with
  measure data and nonlinear potentials}, Ann. Scuola Norm. Super. Pisa, 
  \textbf{19} (1992), 591--613.
  

\bibitem[KM94]{KM1}
T.~Kilpel\"{a}inen and J.~Mal\'y, \emph{The {W}iener test and potential estimates for
  quasilinear elliptic equations}, Acta Math.  \textbf{172}  (1994), 137--161.

\bibitem[Lab02]{Lab1}
D.~Labutin, \emph{Potential estimates for a class of fully nonlinear elliptic
  equations}, Duke Math. J. \textbf{111} (2002), 1--49.
  
  \bibitem[LU68]{LU68} O. A. Ladyzhenskaya and N. N. Ural'tseva, \emph{Linear and Quasilinear Elliptic Equations,}  Academic Press, New York-London, 1968. 


\bibitem[Maz85]{MazSob}
V.~Maz'ya, \emph{Sobolev Spaces}, Springer Series in Soviet Math.,
  Springer, Berlin, 1985 (new edition in press).



\bibitem[Min07]{Min07}
G.~Mingione, \emph{The Calder\'{o}n-Zygmund theory for elliptic problems with measure data}.
Ann. Scuola Norm. Super. Pisa \textbf{6} (2007),  195--261. 

\bibitem[Min11]{Min11}
G.~Mingione, 
 \emph{Gradient potential estimates}, J. Europ. Math. Soc. (JEMS) 
 \textbf{13} (2011),  459--486.
 
 \bibitem[MP02]{MP02}
F. Murat and A. Porretta, \emph{Stability properties, existence, and nonexistence of renormalized solutions for elliptic equations with measure data,} Commun. PDE \textbf{27} (2002),  2267--2310. 


\bibitem[NTV99]{NTV99}
F.~Nazarov,  S.~Treil and A.~Volberg, \emph{The {B}ellman functions and
 two-weight inequalities for {H}aar multipliers,}
J. Amer. Math. Soc. \textbf{12} (1999), 909--928.

\bibitem[NTV03]{NTV03}
F.~Nazarov,  S.~Treil and A.~Volberg, \emph{The Tb-theorem on non-homogeneous spaces,} Acta Math. \textbf{190} (2003),  151--239.


\bibitem[PV06]{PV0}
N.~C. Phuc and I.~E. Verbitsky, \emph{Local integral estimates and removable singularities for 
quasilinear and {H}essian equations with nonlinear source 
terms}, Commun. PDE \textbf{31} (2006),  1779--1791. 

\bibitem[PV08]{PV}
N.~C. Phuc and I.~E. Verbitsky, \emph{Quasilinear and {H}essian equations of
  {L}ane-{E}mden type}, Ann. Math. \textbf{168} (2008), 859--914.


\bibitem[PV09]{PV2}
N.~C. Phuc and I.~E. Verbitsky, \emph{Singular quasilinear and Hessian equations and inequalities},
  J. Funct. Anal. \textbf{256} (2009),  1875--1905.


\bibitem[Por02]{Por02} A. Porretta, \emph{Nonlinear equations with natural growth terms and measure data,} Proc. 2002 Fez Conf. PDE, Electron. J. Differ. Equ. Conf., \textbf{9}, Southwest Texas State Univ., San Marcos, TX (2002), 183--202. 


\bibitem[PS06]{PS06} A. Porretta and S. Segura de Le\'{o}n, \emph{Nonlinear elliptic equations having a gradient term with natural growth,} 
J. Math. Pures Appl. \textbf{85} (2006),  465--492. 

\bibitem[Ser64]{Ser2}
J.~Serrin, \emph{Local behaviour of solutions to quasi-linear equations}, Acta
  Math. \textbf{111} (1964), 247--301.



\bibitem[Sho97]{Sho97} R. E. Showalter, \emph{Monotone Operators in Banach Space and Nonlinear Partial Differential Equations,} Math.  Surveys and Monographs, \textbf{49}, Amer.  Math.  Soc., Providence, RI, 1997. 

\bibitem[Tru67]{T67}
N.~Trudinger, \emph{On Harnack type inequalities and their application to quasilinear elliptic equations}, Commun. Pure Appl. Math. \textbf{20} (1967) 721--747.

\bibitem[TW99]{TW2}
N.~Trudinger and X.-J. Wang, \emph{Hessian measures {II}}, Ann. Math.
  \textbf{150} (1999),  579--604.



\bibitem[TW02b]{TW1}
N.~Trudinger and X.-J. Wang, \emph{On the weak continuity of elliptic operators and applications to
  potential theory}, \textbf{124} Amer. J. Math. (2002),  369--410.

\bibitem[Ver96]{Verb96}
I. E.~Verbitsky, \emph{Imbedding and multiplier theorems for discrete Littlewood-Paley spaces,} Pacific J. Math. \textbf{176} (1996), 529--556. 

\end{thebibliography}
\end{document}